\documentclass[11pt]{article}%
\usepackage{amsmath}
\usepackage{amsfonts}
\usepackage{amssymb}
\usepackage{graphicx}
\usepackage{acronym}
\usepackage{amsthm}%
\newtheorem{theorem}{Theorem}[section]

\theoremstyle{definition}
\newtheorem{definition}[theorem]{Definition}
\theoremstyle{plain}

\newtheorem{lemma}[theorem]{Lemma}
\newtheorem{proposition}[theorem]{Proposition}
\theoremstyle{remark}
\newtheorem{remark}[theorem]{Remark}
\numberwithin{equation}{section}
\begin{document}

\title{The Euler Characteristic and Finiteness Obstruction of Manifolds with Periodic Ends}
\author{John G. Miller\\Indiana U.-Purdue U. at Indianapolis}
\date{}
\maketitle

\begin{abstract}
Let $M$ be a complete orientable manifold of bounded geometry. Suppose that
$M$ has finitely many ends, each having a neighborhood quasi-isometric to a
neighborhood of an end of an infinite cyclic covering of a compact manifold.
We consider a class of exponentially weighted inner products $\left(
\cdot,\cdot\right)  _{k}$ on forms, indexed by $k>0$. Let $\delta_{k}$ be the
formal adjoint of $d$ for $\left(  \cdot,\cdot\right)  _{k}.$ It is shown that
if $M$ has finitely generated rational homology, $d+\delta_{k}$ is Fredholm on
the weighted spaces for all sufficiently large $k.$ The index of its
restriction to even forms is the Euler characteristic of $M.$

This result is generalized as follows. Let $\pi=\pi_{1}\left(  M\right)  .$
Take $d+\delta_{k}$ with coefficients in the canonical $C^{\ast}\left(
\pi\right)  $-bundle $\psi$ over $M.$ If the chains of $M$ with coefficients
in $\psi$ are $C^{\ast}\left(  \pi\right)  $-finitely dominated, then
$d+\delta_{k}$ is Fredholm in the sense of Mi\u{s}\u{c}enko and Fomenko for
all sufficiently large $k.$ The index in $\tilde{K}_{0}\left(  C^{\ast}\left(
\pi\right)  \right)  $ is related to Wall's finiteness obstruction. Examples
are given where it is nonzero.

\end{abstract}

\setcounter{section}{-1}

\section{\label{s0}Introduction}

The analytic index of the operator $d+\delta$ on a compact orientable
Riemannian manifold $M^{n}$ is the Euler characteristic of $M,$ $\chi\left(
M\right)  .$ This paper extends this result to a class of complete noncompact
manifolds, those with finitely generated rational homology and finitely many
quasi-periodic ends. The latter term means that there is a neighborhood of
each end which is quasi-isometric to a neighborhood of an end of an infinite
cyclic covering of a smooth compact manifold. One reason for interest in such
manifolds is a result stated by Siebenmann \cite{sie} and proved by Hughes and
Ranicki \cite{hura}: if $M$ is a manifold of dimension greater than 5 with
finitely many ends satisfying a certain tameness condition, then each end has
a neighborhood homeomorphic to a neighborhood of an end of an infinite cyclic
covering of a compact topological manifold.

$d+\delta$ acting on $\mathcal{L}^{2}$ forms is a Fredholm operator only in
special circumstances. We consider more generally weighted $\mathcal{L}^{2}$
spaces. These were first used in index theory on manifolds with asymptotically
cylindrical ends by Lockhart and McOwen \cite{lomc} and Melrose and Mendoza.
Let $\rho\left(  x\right)  $ be a smooth nonnegative function on $M$ with
bounded gradient which tends to $\infty$ at $\infty.$ Let $k>0.$ The weighted
inner product on compactly supported smooth forms is $\left(  u,v\right)
_{k}=(k^{\rho\left(  x\right)  }u,k^{\rho\left(  x\right)  }v),$ where
$\left(  \cdot,\cdot\right)  $ is the $\mathcal{L}^{2}$ inner product. The
weighted forms are obtained by completion. In other words, they are the
$\mathcal{L}^{2}$ space of the measure $k^{2\rho\left(  x\right)  }dx,$ where
$dx$ is the Riemannian measure. In the quasi-periodic case $\rho\left(
x\right)  $ is chosen to change approximately linearly under iterated covering
translations. We consider the operator $D_{k}=d+\delta_{k},$ where $\delta
_{k}$ is the formal adjoint of $d$ for the weighted inner product. $D_{k}$ is
essentially self-adjoint. We denote by $\bar{D}_{k}$ the closure of $D_{k}.$
Let $\bar{D}_{k}^{even}$ be its restriction to even forms. Let $\chi$ and
$\chi^{\ell f}$ be the Euler characteristic of the homology and locally finite
homology of $M.$ The first main result follows.

\begin{theorem}
\label{th0.1}Let $M^{n}$ be a complete connected Riemannian manifold of
bounded geometry. $\bar{D}_{k}$ is Fredholm if and only if $\bar{D}_{1/k}$ is,
and the indexes satisfy $Ind$ $\bar{D}_{1/k}^{even}=\left(  -\right)  ^{n}Ind$
$\bar{D}_{k}^{even}.$ If $M$ has finitely generated rational homology and
finitely many quasi-periodic ends, $\bar{D}_{k}$ is Fredholm for all $k>0$
which are sufficiently large or small. The index of $\bar{D}_{k}^{even}$ is
\[
\left\{
\begin{array}
[c]{l}%
\left(  -\right)  ^{n}\chi\\
\left(  -\right)  ^{n}\chi^{\ell f}=\chi
\end{array}
\right.  \text{for all }k>0\text{ which are sufficiently }\left\{
\begin{array}
[c]{l}%
\text{large.}\\
\text{small.}%
\end{array}
\right.
\]

\end{theorem}

The factors of $\left(  -\right)  ^{n}$ and the relation $\chi^{\ell
f}=\left(  -\right)  ^{n}\chi$ come from Poincar\'{e} duality. This is a
special case of a more general theorem involving an analytical version of
Wall's finiteness obstruction. For a ring $R,$ a complex of $R$-modules is
said to be $R$-finitely dominated if it is equivalent to a finite dimensional
complex of finitely generated projective $R$-modules. Then $\chi_{R}\in
K_{0}\left(  R\right)  $ is the Euler characteristic, and $\tilde{\chi}_{R}%
\in\tilde{K}_{0}\left(  R\right)  $ is its reduction. Let $X$ be a $CW$
complex, $\tilde{X} $ its universal covering, and $\pi$ the group of covering
transformations. If $X$ is dominated by a finite complex, or equivalently
$\pi$ is finitely presented and the cellular chains $C_{\ast}\left(  \tilde
{X}\right)  $ are $\mathbb{Z}\left(  \pi\right)  $-finitely dominated, then
Wall's obstruction $o_{M}\in\tilde{K}_{0}\left(  \mathbb{Z}\left(  \pi\right)
\right)  $ is defined. It is the Euler characteristic of $C_{\ast}\left(
\tilde{X}\right)  . $ Its vanishing is necessary and sufficient for $X$ to
have the homotopy type of a finite complex.

Let $\pi$ be the group of a regular covering of $M,$ and $C^{*}\left(
\pi\right)  $ be the group $C^{*}$-algebra. There is a canonical bundle $\psi$
with fiber $C^{*}\left(  \pi\right)  $ over $M.$ If the local coefficient
chains of $M$ with coefficients in $\psi$ are $C^{*}\left(  \pi\right)
$-finitely dominated, then $\chi_{C^{*}\left(  \pi\right)  }$ is defined. For
the trivial group and $R$ a field of characteristic $0,$ finite domination is
the same as finitely generated rational homology. The augmentation
$K_{0}\left(  C^{*}\left(  \pi\right)  \right)  \rightarrow K_{0}\left(
\mathbb{C}\right)  =\mathbb{Z}$ takes $\chi_{C^{*}\left(  \pi\right)  }$ to
$\chi.$ If $M$ is dominated by a finite complex and $\pi$ is the group of the
universal covering, $\mathbb{Z}\left(  \pi\right)  \rightarrow C^{*}\left(
\pi\right)  $ takes $o_{M}$ to $\tilde{\chi}_{C^{*}\left(  \pi\right)  }.$

A locally finite Euler characteristic $\chi_{C^{\ast}\left(  \pi\right)
}^{\ell f}$ is defined similarly if the locally finite chains of $M$ with
coefficients in $\psi$ are $C^{\ast}\left(  \pi\right)  $-finitely dominated.
It reduces to $\chi^{\ell f}$ for the trivial group. We replace $D_{k}$ by the
same operator with coefficients in $\psi$ without changing notation. By
\textquotedblleft Fredholm\textquotedblright\ in the context of operators over
$C^{\ast}$-algebras we mean Fredholm in the sense of Mi\u{s}\u{c}enko and Fomenko.

\begin{theorem}
\label{th0.2}Theorem \ref{th0.1} holds with the following changes: in place of
finitely generated rational homology we assume that the local coefficient
chains of $M$ with coefficients in $\psi$ are $C^{\ast}\left(  \pi\right)
$-finitely dominated. $\chi$ and $\chi^{\ell f}$ are replaced by
$\chi_{C^{\ast}\left(  \pi\right)  }$ and $\chi_{C^{\ast}\left(  \pi\right)
}^{\ell f}.$
\end{theorem}

This is actually proved with a fundamental group hypothesis. Let $\bar
{N}\rightarrow N$ be the model infinite cyclic covering for an end of $M.$ We
assume that $\pi_{1}\left(  N\right)  =$ $\pi_{1}\left(  \bar{N}\right)
\times\mathbb{Z}.$ This is to avoid dealing with twisted group rings.

It seems very possible that the homomorphism $\tilde{K}_{0}\left(
\mathbb{Z}\left(  \pi\right)  \right)  \rightarrow\tilde{K}_{0}\left(
C^{\ast}\left(  \pi\right)  \right)  $ is always $0.$ This is the case if
$C^{\ast}\left(  \pi\right)  $ is replaced by the group von Neumann algebra
\cite{sch}. However, a manifold may be $C^{\ast}\left(  \pi\right)  $-finitely
dominated without being finitely dominated. In this case $\tilde{\chi
}_{C^{\ast}\left(  \pi\right)  }$ is still a finiteness obstruction, since it
vanishes if $M$ has the homotopy type of a finite complex. We give examples of
manifolds with finite fundamental group for which the above indexes are
nontrivial. The index is just the $\pi$-equivariant Euler characteristic.
Examples with infinite fundamental group are obtained using free products and
semidirect products.

The proofs are based on a connection between exponential weights and boundedly
controlled topology. A translation of Euclidean space induces a bounded
operator on exponentially weighted spaces. In general, we say that an operator
is spatially bounded if, roughly speaking, it moves things a bounded distance.
This is the boundedness of bounded topology. It is related to, but different
from, the finite propagation of Roe and Higson \cite[Chs. 3, 4]{roe}. The
underlying principle is that, frequently, a spatially bounded operator is
analytically bounded on exponentially weighted spaces.

The main point is to show that weighted complexes of forms are chain
equivalent to standard cochain complexes. Let $\Omega_{c}$ be the forms with
coefficients in $\psi$ with compact supports. Let $\bar{\Omega}_{d,k}$ be the
domain of the closure of $d$ acting on $\Omega_{c}$ in the $k$-norm. We make
the same fundamental group hypothesis as for Theorem \ref{th0.2}.

\begin{theorem}
\label{th0.3}Under the conditions of Theorem \ref{th0.2}, $\bar{\Omega}_{d,k}
$ is equivalent to the compactly supported simplicial cochains $C_{c}^{\ast
}\left(  M;\psi\right)  $ for $k$ large, and to the simplicial cochains
$C^{\ast}\left(  M;\psi\right)  $ for $k>0$ small.
\end{theorem}

The idea for this is as follows. Suppose that the complement of some compact
set in $M$ is isometric to $V\times\lbrack0,\infty),$ with $V$ of dimension
$n-1.$ Let $u$ be any smooth form. Pushing in along the normal rays induces a
form from $u$ which satisfies the $k$-growth condition for any $0<k<1.$ This
gives an equivalence of the two spaces. There is a related argument in the
other case. More details can be found in \cite[6.4]{mel}. We will carry out a
controlled pushing operation in some cases where $M$ doesn't admit a boundary.

We proceed by several reductions. The first is from weighted forms to weighted
simplicial cochains. This uses a de Rham-type theorem extending one of Pansu
for the $\mathcal{L}^{2}$ cohomology of manifolds of bounded geometry. The
theorem incorporates both weights and spatial boundedness. The problem is then
transferred to an algebraic complex for the infinite cyclic cover modelling an
end. This is a direct translation into analysis of the framework of Hughes and
Ranicki. The complex has the structure of a doubly infinite algebraic mapping
telescope, which may be pushed either off one of its ends or to infinity.
Analytically, this amounts to the invertibility of a standard weighted shift
operator. This is an analog of Ranicki's result on the vanishing of homology
with Novikov ring coefficients \cite{ran2}.

There are a number of further connections with other work. Among these are
Taubes' study of analysis on manifolds with periodic ends, and a conjecture of
Bueler on weighted $\mathcal{L}^{2}$ cohomology. A discussion is given at the
end of the paper.

We make use of the standard material on Hilbert $C^{*}$-modules, which may be
found in \cite[Ch. 15]{w-o}. $A$ will always denote a unital $C^{*}$-algebra.
All modules will be separable. The compact operators on an $A$-module $P$ are
$\mathcal{K}_{A}\left(  P\right)  .$ The distinction between the adjointable
operators $\mathcal{L}_{A}\left(  P\right)  $ and the bounded ones
$\mathcal{B}_{A}\left(  P\right)  $ is crucial at some points. All chain
complexes will be finite dimensional. The proofs in the references are often
for $A=\mathbb{C}.$ They have been chosen so that little or no change is
required to make them valid for general $A.$

The contents are as follows: Section \ref{s1} contains background material,
and accomplishes the proof of Theorem \ref{th0.2} using results from later
sections. Section \ref{s2} introduces spatial boundedness and contains the
proof of the de Rham theorem. Section \ref{s3} is about algebraic versions of
infinite cyclic covers and mapping telescopes. It completes the proof of
Theorem \ref{th0.3}. Section \ref{s4} contains background on finiteness
obstructions and examples for Theorems \ref{th0.1} and \ref{th0.2}. Section
\ref{s5} is the analytic basis for the paper. It shows that the differential
operators we use have the expected properties. We prove a mild extension of a
theorem of Kasparov, which he stated with only a brief sketch of proof.
Section \ref{s6} is the discussion.

To a large extent, this paper is an analytical version of parts of the book of
Hughes and Ranicki. The text doesn't acknowledge all of my borrowings. I wish
to thank Jonathan Rosenberg for suggestions and encouragement at the beginning
of this project.

\section{\label{s1}Forms and weights}

This section contains preliminaries and the proof of Theorem \ref{th0.2},
assuming the results of the remainder of the paper.

\subsection{\label{ss1.1}}

Let $M^{n}$ be a complete, oriented, connected Riemannian manifold. Let
$\Lambda$ be the complexified exterior algebra bundle of the cotangent bundle.
The forms on $M$ with compact support $\Omega_{c}$ are the compactly supported
smooth sections of $\Lambda.$ Let $\ast$ be the Hodge operator. For
$u,v\in\Omega_{c}^{p},$ a pointwise inner product is defined by $\left\langle
u,v\right\rangle \left(  x\right)  =\ast\left(  \bar{u}\left(  x\right)
\wedge\ast v\left(  x\right)  \right)  .$ The bar denotes conjugation, so this
is conjugate-linear in the first variable. A global inner product is defined
by $\left(  u,v\right)  =\int_{M}\left\langle u,v\right\rangle dx.$ Let $A$ be
a unital $C^{\ast}$-algebra. We consider forms with coefficients in a flat
bundle of $A$-modules. This is a bundle $V=\tilde{M}\times_{\pi}P\rightarrow
M,$ with $\tilde{M}$ a regular covering of $M,$ $\pi$ its group, and $P$ a
finitely generated (so projective) Hilbert $A$-module with a unitary
representation of $\pi$. The relation is $\left(  x,p\right)  \sim\left(
gx,gp\right)  .$ The most important case is the canonical bundle $\psi$, where
$P=C^{\ast}\left(  \pi\right)  $ and the regular representation is used. $V$
has a natural flat connection. Let $\Omega_{V,c}$ be the compactly supported
smooth sections of $\Lambda\otimes V.$ Let $d_{V}$ be the exterior derivative
with coefficients in $V.$ Since the connection is flat, $\left(  d_{V}\right)
^{2}=0.$ Thus we have a de Rham complex with coefficients in $V.$

An $A$-valued inner product is determined as follows: If $u,v\in\Omega_{V,c}$
can be written as $s\otimes k,$ $t\otimes\ell,$ with $s,t\in\Omega_{c}^{p}$
and $k,\ell$ sections of $V,$ let $\left\langle u\left(  x\right)  ,v\left(
x\right)  \right\rangle _{V}=\left\langle s\left(  x\right)  ,t\left(
x\right)  \right\rangle \left\langle k\left(  x\right)  ,\ell\left(  x\right)
\right\rangle .$ Then $\left(  u,v\right)  _{V}=\int_{M}\left\langle
u,v\right\rangle _{V}dx.$ All the integrals in this paper are Riemann. This
makes $\Omega_{V,c}$ into a complex of pre-Hilbert $A$-modules. Henceforth we
will usually drop $V$ from the notation and just write $\Omega_{c}$. There is
a star operator given by $\ast\left(  s\otimes k\right)  =\ast s\otimes k.$

We will define weighted inner products on $\Omega_{c},$ generalizing the
$\mathcal{L}^{2}$ inner products defined above. See \cite[Section 2]{bue} for
more details. Let $h\left(  x\right)  $ be a smooth real function on M. Let
$d\mu=e^{2h\left(  x\right)  }dx,$ and $\left(  u,v\right)  _{\mu}=\int
_{M}\left\langle u,v\right\rangle d\mu=\left(  e^{h}u,e^{h}v\right)  $ The
weights that will be used in this paper are much more special. Let
$\rho\left(  x\right)  $ be a smooth real function on $M$ with bounded
gradient. Let $h\left(  x\right)  =\rho\left(  x\right)  \log k$ for some
$k>0.$ Then $d\mu=k^{2\rho\left(  x\right)  }dx.$ In this situation we will
write $\left(  \cdot,\cdot\right)  _{\mu}=\left(  \cdot,\cdot\right)  _{k}.$
The case $k=1$ is the $\mathcal{L}^{2}$ inner product, in which case we will
often simply write $\left(  \cdot,\cdot\right)  .$ $\Omega_{c}$ with such an
inner product will be denoted by $\Omega_{\mu},$ or by $\Omega_{k}$ when using
the $k$-inner products. The completions are $\bar{\Omega}_{\mu}$ and
$\bar{\Omega}_{k}. $ The inner products extend by continuity.

Let $\Omega_{d,\mu}$ be $\Omega_{c}$ with the graph inner product $\left(
u,v\right)  _{d,\mu}=\left(  u,v\right)  _{\mu}+\left(  du,dv\right)  _{\mu}.$
The main space of forms we will use is the domain of $\bar{d},$ the closure of
$d $ in the $\mu$-norm. This may be described as the completion of
$\Omega_{d,\mu}.$ We denote it by $\bar{\Omega}_{d,\mu}$ or $\bar{\Omega
}_{d,k}.$ $\bar{d}:\bar{\Omega}_{d,\mu}^{j}\rightarrow\bar{\Omega}_{d,\mu
}^{j+1}$ is bounded.

Let $\delta$ be the $\mathcal{L}^{2}$ formal adjoint of $d$ on $\Omega_{c}. $
One computes that the formal adjoint of $d$ with respect to $\left(
\cdot,\cdot\right)  _{\mu}$ is $\delta_{\mu}=e^{-2h}\delta e^{2h}%
=\delta-2dh\llcorner$ , where $\llcorner$ denotes interior multiplication. Let
$D_{\mu}=d+\delta_{\mu},$ which is formally self-adjoint. Multiplication by
$e^{h}$ induces a unitary between the $\mu$-inner product and the
$\mathcal{L}^{2}$-inner product on $\Omega_{c}$. Then $D_{\mu}$ is unitarily
equivalent to $d+\delta-\left(  dh\wedge+dh\llcorner\right)  .$

Let $C_{b}^{\infty,1}\left(  M\right)  $ be the space of smooth functions
which are bounded and whose differentials are bounded. It has the norm
$\sup_{x\in M}\left\vert \phi\left(  x\right)  \right\vert +\sup_{x\in
M}\left\Vert d\phi\left(  x\right)  \right\Vert .$ The following Lemma is a
standard fact for forms with values in $\mathbb{C}.$ Additional care is
required for coefficients in a $C^{\ast}$-algebra.

\begin{lemma}
\label{le1.1}$C_{b}^{\infty,1}\left(  M\right)  $ acts continuously on
$\bar{\Omega}_{d,\mu}.$
\end{lemma}

%

%TCIMACRO{\TeXButton{Proof}{\proof} }%
%BeginExpansion
\proof
%EndExpansion
For $\phi\in C_{b}^{\infty,1}\left(  M\right)  ,$ $u\in$ $\Omega_{c},$%
\begin{align*}
\left\|  \phi u\right\|  _{d,\mu}^{2}  &  =\left\|  \left(  \phi u,\phi
u\right)  _{\mu}+\left(  d\left(  \phi u\right)  ,d\left(  \phi u\right)
\right)  _{\mu}\right\|  \le\left\|  \left(  \phi u,\phi u\right)  _{\mu
}\right\|  +\left\|  \left(  d\left(  \phi u\right)  ,d\left(  \phi u\right)
\right)  _{\mu}\right\| \\
&  =\left\|  \phi u\right\|  _{\mu}^{2}+\left\|  d\left(  \phi u\right)
\right\|  _{\mu}^{2}=\left\|  \phi u\right\|  _{\mu}^{2}+\left\|  \phi
du+d\phi\wedge u\right\|  _{\mu}^{2}\\
&  \le\left\|  \phi u\right\|  _{\mu}^{2}+2\left\|  \phi du\right\|  _{\mu
}^{2}+2\left\|  d\phi\wedge u\right\|  _{\mu}^{2}\\
&  =\left\|  \left(  \phi u,\phi u\right)  _{\mu}\right\|  +2\left\|  \left(
\phi du,\phi du\right)  _{\mu}\right\|  +2\left\|  \left(  d\phi\wedge
u,d\phi\wedge u\right)  _{\mu}\right\|  .
\end{align*}

\noindent The terms are easily estimated. For example,
\begin{align*}
\left(  d\phi\wedge u,d\phi\wedge u\right)  _{\mu}  &  =\int_{M}\left\langle
d\phi\wedge u,d\phi\wedge u\right\rangle d\mu\\
&  \leq\sup_{x\in M}\left\Vert d\phi\left(  x\right)  \right\Vert ^{2}\int
_{M}\left\langle u,u\right\rangle d\mu=K^{2}\left(  u,u\right)  _{\mu}.
\end{align*}

\noindent Then $\left\|  d\phi\wedge u\right\|  _{\mu}\le K\left\|  u\right\|
_{\mu.}.$%
%TCIMACRO{\TeXButton{End Proof}{\endproof}}%
%BeginExpansion
\endproof
%EndExpansion

\subsection{\textbf{\label{ss1.2}}}

We need some definitions concerning the bounded geometry (BG) category. For
more information see \cite[Appendix 1]{shu}.

\begin{definition}
\label{de1.2}Riemannian metrics $\left\langle \cdot,\cdot\right\rangle $ and
$\left\langle \cdot,\cdot\right\rangle ^{\prime}$ on $M$ are
\textit{{quasi-isometric }}if there exists $C>1$ such that for all $x\in M$
and $X\in TM_{x},$
\end{definition}

$\dfrac1C\left\langle X,X\right\rangle <\left\langle X,X\right\rangle
^{\prime}<C\left\langle X,X\right\rangle .$

It follows that there is $K>1$ such that for all $u\in\Omega_{c},$
$\dfrac1K\left(  u,u\right)  <\left(  u,u\right)  ^{\prime}<K\left(
u,u\right)  .$ A similar statement then holds for the weighted $d$-inner
products, so that the complexes $\bar{\Omega}_{d,\mu}$ are the same, with
equivalent norms.

A manifold of bounded geometry is a Riemannian manifold with certain
uniformity properties. They are of two different types.

\begin{description}
\item[(I)] The injectivity radii at points of $M$ are bounded below by a
constant $r_{0}.$
\end{description}

This condition implies that $M$ is complete. The statement of the second
condition requires the notion of canonical coordinates at a point $x\in M.$
Choose an orthonormal basis in $T_{x}M,$ thus identifying it with
$\mathbb{R}^{n}.$ Choose some $r<r_{0}.$ Then a canonical coordinate
neighborhood of $x$ is given by the exponential map at $x$ restricted to the
open ball of radius $r$ in $\mathbb{R}^{n}.$

\begin{description}
\item[(B1)] For some fixed $r,$ there exists a covering of $M$ by canonical
coordinate neighborhoods such that the differentials of the exponential maps
and their inverses are uniformly bounded.
\end{description}

Examples include compact manifolds and covering spaces of compact manifolds.
Uniform boundedness of some higher derivatives of the transition functions is
often required. These conditions are implied by conditions on the curvature
tensor and its covariant derivatives. In \cite{shu}, all higher derivatives
are assumed uniformly bounded. The statements in the present paper using only
(I) and (B1) come from examining the proofs. With these definitions, it is not
the case that a manifold which is quasi-isometric to a BG manifold is BG.

\subsection{\textbf{\label{ss1.3}}}

Recall that $\mathcal{B}_{A}$ (resp. $\mathcal{L}_{A}$) is the category of
Hilbert $A$-modules and bounded (resp. adjointable) homomorphisms. In the
following discussion \textquotedblleft complex" means \textquotedblleft
cochain complex". Analogous statements hold for complexes. Let $\left(
C,\beta\right)  $ be an $A$-finitely dominated complex in $\mathcal{B}_{A}$.
This means that $C$ is equivalent in $\mathcal{B}_{A}$ to a complex of
finitely generated modules. We may define its Euler characteristic as
$\chi\left(  C\right)  =\sum\left(  -\right)  ^{i}\left[  F_{i}\right]  \in
K_{0}\left(  A\right)  ,$ where $F$ is an equivalent complex of finitely
generated modules. This is independent of the choice of $F,$ since $\chi$ is a
chain homotopy invariant of finitely generated complexes.

We will make use of the theory of Fredholm complexes, introduced by Segal
\cite{seg}. A complex $\left(  C,\beta\right)  $ in $\mathcal{L}_{A}$ is said
to be $A$-Fredholm if there exists a parametrix, a homomorphism $g\in$
$\mathcal{L}\left(  C\right)  $ of degree 1 satisfying $\beta g+g\beta=I+c,$
with $c\in\mathcal{K}\left(  C\right)  .$ A Fredholm operator is a Fredholm
complex $\beta:C_{0}\rightarrow C_{1}$ which is invertible modulo
$\mathcal{K}\left(  C\right)  $. A complex in $\mathcal{L}_{A}$ is Fredholm if
and only if it is finitely dominated in $\mathcal{L}_{A},$ by \cite{kami}
Propositions 3.2 and 3.9. Therefore $\chi\left(  C\right)  $ is defined for a
Fredholm complex. For a Fredholm operator it is called the index of $\beta,$
$Ind\,\beta.$ It has the the expected properties \cite[Ch. 17]{w-o}. The
following Lemma improves on the stated relationship between finite domination
and Fredholm complexes. It is necessary because the equivalences involving
$\bar{\Omega}_{d,k}$ will only be established in $\mathcal{B}_{A}.$

\begin{lemma}
\label{le1.3}A cochain complex $C$ in $\mathcal{L}_{A}$ is Fredholm if and
only if it is finitely dominated in $\mathcal{B}_{A}.$
\end{lemma}

%

%TCIMACRO{\TeXButton{Proof}{\proof}}%
%BeginExpansion
\proof
%EndExpansion
A Fredholm complex is finitely dominated in $\mathcal{L}_{A}$ and thus in
$\mathcal{B}_{A}.$ Let $C$ be equivalent in $\mathcal{B}_{A}$ to the finitely
generated complex $F.$ Since homomorphisms with domain a finitely generated
module are in $\mathcal{K}_{A}$, $F$ is a complex in $\mathcal{L}_{A},$ and
the map $f:F\rightarrow C$ is in $\mathcal{L}_{A}.$ Since $f$ induces an
isomorphism of homology, it has a homotopy inverse in $\mathcal{L}_{A}$
\cite[Prop. 2.7]{kami}. Therefore $C$ is finitely dominated in $\mathcal{L}%
_{A}.$%
%TCIMACRO{\TeXButton{End Proof}{\endproof}}%
%BeginExpansion
\endproof
%EndExpansion

We consider $\tau$-complexes $\left(  E,\beta\right)  $ in the sense of
\cite[Section 1]{mil1}. These are simplified notation for complexes of
differential forms. They are $n$-dimensional cochain complexes $E$ in
$\mathcal{L}_{A}$ with differential $\beta$ and self-adjoint involution
$\tau:E\rightarrow E^{n-\ast}$ satisfying $\beta^{\ast}=\tau\beta\tau$. Let
the dual complex $\left(  E^{\prime},\beta^{\prime}\right)  $ be defined by
$\left(  E^{\prime}\right)  ^{j}=\left(  E^{n-j}\right)  ^{\prime}$ and
$\left(  \beta^{\prime}\right)  ^{j}=\left(  \beta^{n-j-1}\right)  ^{\prime}.$
The map $\phi:E\rightarrow E^{\prime}$ defined by $\phi\left(  u\right)
\left(  v\right)  =\left(  u,\tau v\right)  $ is an isomorphism. It is shown
in \cite[Th. 3.3]{kami} that for a Fredholm $\tau$-complex, the signature
operator $S=-i\left(  d-\tau d\tau\right)  $ is an $A$-Fredholm operator. It
is self-adjoint. It follows that $S^{even}:E^{even}\rightarrow E^{odd}$ is
Fredholm. The adjoint of $S^{even}$ is $S^{odd}.$ The following replaces a
standard Hodge theory argument for $A=\mathbb{C}.$ The first part of the proof
is taken from Segal \cite[Section 5]{seg}. We use the notation $\approx$ for
congruence modulo $\mathcal{K}_{A}.$

\begin{proposition}
\label{pr1.4}If $\left(  E,\beta\right)  $ is a cochain complex in
$\mathcal{L}_{A}$ such that $S=-i\left(  \beta-\beta^{\ast}\right)  $ is
Fredholm, then $E$ is a Fredholm complex and $Ind$ $S^{even}=\chi\left(
E\right)  \in K_{0}\left(  C^{\ast}\left(  \pi\right)  \right)  $.
\end{proposition}

%

%TCIMACRO{\TeXButton{Proof}{\proof}}%
%BeginExpansion
\proof
%EndExpansion
Let $E$ be any Fredholm complex. A parametrix $g$ may be chosen so that
$g^{2}\simeq0.$ In fact, if $g$ is any parametrix, then $g\beta g$ has this
property. For any such $g,$ $\beta+g:E^{even}\rightarrow E^{odd}$ is a
Fredholm operator, since $\left(  \beta+g\right)  ^{2}\approx I.$ We claim
that $Ind$ $\left(  \beta+g\right)  $ is independent of the choice of such a
$g.$ If $g_{0}$ and $g_{1}$ are parametrices for $E,$ $g_{t}=\left(
1-t\right)  g_{0}+tg_{1}$ is a norm-continuous family of parametrices. The
same is true of $g_{t}\beta g_{t}.$ Thus $Ind$ $\left(  \beta+g_{0}\beta
g_{0}\right)  =Ind$ $\left(  \beta+g_{1}\beta g_{1}\right)  .$ Now suppose
that $g_{0}^{2}\approx0$ and $g_{1}^{2}\approx0.$ Then $g_{0}-g_{0}\beta
g_{0}=g_{0}\left(  1-\beta g_{0}\right)  \approx g_{0}^{2}\beta\approx0.$
Therefore $\beta+g_{0}$ is Fredholm and has the same index as $\beta
+g_{0}\beta g_{0}.$ Similarly for $g_{1}.$ We conclude that $Ind$ $\left(
\beta+g_{0}\right)  =Ind$ $\left(  \beta+g_{1}\right)  .$ We can thus refer to
$Ind$ $E.$

Suppose that $E$ is contractible. Then there exists $g$ such that $\beta
g+g\beta=I.$ $\beta\left(  g\beta g\right)  +\left(  g\beta g\right)
\beta=\beta\left(  I-\beta g\right)  g+g\left(  I-g\beta\right)  \beta
=g\beta+\beta g=I,$ so $g\beta g$ is again a contraction. Therefore
$\beta+g\beta g $ is an isomorphism, so has index 0. It is shown in
\cite{kami}, proof of Proposition 2.9, that for any Fredholm complex $E,$
There exist a finitely generated complex $F$ and contractible complexes $M$
and $N$ such that $E\oplus M\cong F\oplus N.$ By additivity, $Ind$%
~$E=Ind$~$F=\left[  F^{even}\right]  -\left[  F^{odd}\right]  =\sum\left(
-\right)  ^{i}\left[  F^{i}\right]  =\chi\left(  E\right)  .$

Now let $E$ be such that $S$ is Fredholm. Let $\Delta=S^{odd}S^{even}.$ This
is self adjoint Fredholm, so $Ind$ $S^{even}=-Ind$ $S^{odd}.$ Let
$\Delta^{\prime}$ be an inverse for $\Delta$ $\operatorname{mod}\mathcal{K}.$
Then $\Delta^{\prime}$ is self adjoint $\operatorname{mod}$ $\mathcal{K}.$ For
$\left(  \Delta^{\prime}\right)  ^{\ast}\Delta=\left(  \Delta\Delta^{\prime
}\right)  ^{\ast}\approx I,$ and similarly $\Delta\left(  \Delta^{\prime
}\right)  ^{\ast}\approx I.$ But $\Delta^{\prime}$ is unique
$\operatorname{mod}\mathcal{K}$, so $\left(  \Delta^{\prime}\right)  ^{\ast
}\approx\Delta^{\prime}.$

$\Delta$ commutes with $\beta$ and $\beta^{\ast}.$ It follows that
$\Delta^{\prime}$ commutes $\operatorname{mod}$ $\mathcal{K}$ with $\beta$ and
$\beta^{\ast}.$ For if $T$ is an operator such that $\Delta T\approx T\Delta,$
then $\Delta^{\prime}T\approx\Delta^{\prime}T\Delta\Delta^{\prime}%
\approx\Delta^{\prime}\Delta T\Delta^{\prime}\approx T\Delta^{\prime}.$ Let
$g=\beta^{\ast}\Delta^{\prime}$. Then $g$ is a parametrix for $E,$ since
$\beta\beta^{\ast}\Delta^{\prime}+\beta^{\ast}\Delta^{\prime}\beta\approx
\beta\beta^{\ast}\Delta^{\prime}+\beta^{\ast}\beta\Delta^{\prime}=\Delta
\Delta^{\prime}\simeq I.$ Thus $\beta+\beta^{\ast}\Delta^{\prime}$ is
Fredholm. Also $\left(  \beta^{\ast}\Delta^{\prime}\right)  \left(
\beta^{\ast}\Delta^{\prime}\right)  \approx\left(  \beta^{\ast}\right)
^{2}\left(  \Delta^{\prime}\right)  ^{2}=0.$ Therefore $Ind$ $E=Ind\left(
\beta+\beta^{\ast}\Delta^{\prime}\right)  .$ But $\left(  \beta+\beta^{\ast
}\Delta^{\prime}\right)  S^{odd}\approx-i\left(  \beta^{\ast}\beta
\Delta^{\prime}-\beta\beta^{\ast}\right)  $ is skew-adjoint
$\operatorname{mod}\mathcal{K},$ so has index 0. Thus $Ind$ $E=-Ind$
$S^{odd}=Ind$ $S^{even}.$ It follows that $Ind$ $S^{even}=\chi\left(
E\right)  $.%
%TCIMACRO{\TeXButton{End Proof}{\endproof}}%
%BeginExpansion
\endproof
%EndExpansion

\subsection{\textbf{\label{ss1.4}}}

In Section \ref{s5}, extending a theorem of Kasparov \cite{kas1}, we show that
if $M$ is of bounded geometry and $V=\psi,$ then $\bar{D}_{\mu}$ and $\bar
{D}_{\mu}^{2}$ are symmetric with real spectrum, that of the latter lying in
$[0,\infty)$. In particular, this allows us to construct operators like
$\left(  \bar{D}_{\mu}^{2}+I\right)  ^{-1/2}$. We also will use $d_{\mu}%
^{\ast}\bar{d},$ where the adjoint is taken with respect to the $\mu$-inner
product. It is symmetric with nonnegative spectrum.

Let $E_{\mu}$ be the complex with $E_{\mu}^{j}$ $=$ $\bar{\Omega}_{\mu}^{j}$
and differential $d_{E_{\mu}}=\bar{d}\left(  \bar{D}_{\mu}^{2}+I\right)
^{-1/2}.$ Proofs of the following statements are in Section \ref{ss5.3}:
$d_{E_{\mu}}$ is bounded with adjoint $d_{\mu}^{\ast}\left(  \bar{D}_{\mu}%
^{2}+I\right)  ^{-1/2}$; $\left(  d_{\mu}^{\ast}\bar{d}+I\right)  ^{1/2}%
:\bar{\Omega}_{d,\mu}\rightarrow\bar{\Omega}_{\mu}$ is a degree-preserving
unitary. It is shown that this is a cochain isomorphism $\left(  \bar{\Omega
}_{d,\mu},\bar{d}\right)  \rightarrow\left(  E_{\mu},d_{E_{\mu}}\right)  .$

It is emphasized by Bueler \cite{bue} that the reason why weighted spaces are
interesting with respect to cohomology is that they do \textit{not }satisfy
the self-duality implied by the definition of $\tau$-complex. If
$d\mu=e^{2h\left(  x\right)  }dx,$ let $d\mu^{-}=e^{-2h\left(  x\right)  }dx.$
Let
\[
\beta_{\mu}^{j}=\left\{
\begin{array}
[c]{rl}%
id_{E_{\mu}}^{j} & \text{ }j\text{ even }\\
d_{E_{\mu}}^{j} & j\text{ odd}%
\end{array}
\right.  ,\qquad\tau_{\mu}^{j}=\left\{
\begin{array}
[c]{rl}%
ie^{2h}*^{j} & \text{ }n\text{ even and }j\text{ odd}\\
e^{2h}*^{j} & \text{otherwise}%
\end{array}
\right.  .
\]

\noindent$\tau_{\mu}$ is a unitary $E_{\mu}^{\ast}\rightarrow E_{\mu^{-}%
}^{n-\ast}$ with $\tau_{\mu}^{\ast}=\tau_{\mu^{-}}.$ By Lemma \ref{le5.8},
$\tau_{\mu}\beta_{\mu}\tau_{\mu^{-}}=\beta_{\mu^{-}}^{\ast}.$ The map
$\phi:\left(  E_{\mu},\beta_{\mu}\right)  \rightarrow\left(  E_{\mu^{-}%
}^{^{\prime}},\beta_{\mu^{-}}^{^{\prime}}\right)  $ defined by $\phi\left(
u\right)  \left(  v\right)  =\left(  u,\tau_{\mu^{-}}v\right)  _{\mu}$ is an
isomorphism. We define a $\tau$-complex structure on $E_{\mu}\oplus E_{\mu
^{-}}$. Let $\beta=\beta_{\mu}\oplus\beta_{\mu^{-}},$ and
\[
\tau=\left(
\begin{array}
[c]{rr}%
0 & \tau_{\mu^{-}}\\
\tau_{\mu} & 0
\end{array}
\right)  .
\]

\noindent$\tau$ is a self-adjoint unitary. $\beta^{*}=\tau\beta\tau,$ so we
have a $\tau$-complex. The signature operator is $S_{\mu}\oplus S_{\mu^{-}%
}=-i\left(  \beta_{\mu}-\beta_{\mu}^{*}\right)  \oplus-i\left(  \beta_{\mu
^{-}}-\beta_{\mu^{-}}^{*}\right)  .$ We find that $\tau_{\mu}S_{\mu}\tau
_{\mu^{-}}=-S_{\mu^{-},}$ so one is Fredholm if and only if the other is, and
$S_{\mu}\oplus S_{\mu^{-}}$ is Fredholm if and only if either is. If $n$ is
even, $\tau_{\mu}S_{\mu}^{even}\tau_{\mu^{-}}=-S_{\mu^{-}}^{even}$ and $Ind$
$S_{\mu^{-}}^{even}=Ind$ $S_{\mu}^{even}$. If $n$ is odd, $\tau_{\mu}S_{\mu
}^{even}\tau_{\mu^{-}}=-S_{\mu^{-}}^{odd}=\left(  -S_{\mu^{-}}^{even}\right)
^{*},$ so $Ind$ $S_{\mu^{-}}^{even}=-Ind$ $S_{\mu}^{even}.$

Thus to get Theorem \ref{th0.2} it is sufficient that either half of Theorem
\ref{th0.3} holds. If $A=\mathbb{C},$ a straightforward application of Hodge
theory shows that the two halves of Theorem \ref{th0.3} are equivalent.
However, there doesn't seem to be a direct argument in general. Therefore we
will continue with the two cases in parallel.

The analog of the usual signature operator on weighted spaces is $\bar{D}%
_{\mu}=\bar{d}+d_{\mu}^{*}.$ The standard bounded operator on $\bar{\Omega
}_{\mu}$ corresponding to this is $\bar{D}_{\mu}\left(  \bar{D}_{\mu}%
^{2}+I\right)  ^{-1/2}.$ The latter is unitarily equivalent to $S_{\mu}. $ For
let $\alpha$ act on $\bar{\Omega}_{\mu}^{j}$ by $i^{\left[  j/2\right]  } $
(the greatest integer function). Then $\alpha S_{\mu}\alpha^{*}=\bar{D}_{\mu
}\left(  \bar{D}_{\mu}^{2}+I\right)  ^{-1/2}.$ Therefore we may refer to $Ind$
$S_{\mu}^{even}$ as $Ind$ $\bar{D}_{\mu}^{even}.$

\subsection{\textbf{\label{ss1.5}}}

We complete the proof of Theorem \ref{th0.2}. From now on we use $k$-inner
products. By above discussion, we are interested in the operators $S_{k}.$ If
$M$ is of bounded geometry $S_{k}$ exists. $S_{k}$ is Fredholm if and only if
$S_{1/k}$ is, in which case $Ind$ $S_{k}=\left(  -\right)  ^{n}Ind$ $S_{1/k}.$

Let $M$ have finitely many quasi-periodic ends. Assume that $C_{\ast}\left(
M;\psi\right)  $ is $A$-finitely dominated. By Theorem \ref{th0.3},
$\bar{\Omega}_{d,k}$ is equivalent to $C_{c}^{\ast}\left(  M;\psi\right)  $
for $k$ large and to $C^{\ast}\left(  M;\psi\right)  $ for $k>0$ small. By
Poincar\'{e} duality, these are equivalent (up to sign) to $C_{n-\ast}\left(
M;\psi\right)  $ and $C_{n-\ast}^{\ell f}\left(  M;\psi\right)  .$ By Lemma
\ref{le4.2}, $C_{\ast}^{\ell f}\left(  M;\psi\right)  $ is finitely dominated
and $\chi_{C^{\ast}\left(  \pi\right)  }^{\ell f}=\left(  -\right)  ^{n}%
\chi_{C^{\ast}\left(  \pi\right)  }$. Thus, under the conditions on $k,$
$\bar{\Omega}_{d,k}$ is finitely dominated and $\chi\left(  \bar{\Omega}%
_{d,k}\right)  =\left(  -\right)  ^{n}\chi_{C^{\ast}\left(  \pi\right)  },$
and $\chi\left(  \bar{\Omega}_{d,k}\right)  =\left(  -\right)  ^{n}%
\chi_{C^{\ast}\left(  \pi\right)  }^{\ell f}=\chi_{C^{\ast}\left(  \pi\right)
}.$ $\bar{\Omega}_{d,k}$ is equivalent to $\left(  E_{k},d_{E_{k}}\right)  .$
The factors of $i$ in the definition of $\beta_{k}$ don't affect finite
domination or Euler characteristic. (Do the same to an equivalent finitely
generated complex.) Therefore $\left(  E_{k},\beta_{k}\right)  $ is finitely
dominated with the same Euler characteristic. By Lemma \ref{le1.3}, the $\tau
$-complex $E_{k}\oplus E_{1/k}$ is Fredholm, since it is the sum of two
finitely dominated complexes. Then its signature operator $S_{k}\oplus
S_{1/k}$ is Fredholm, so $S_{k}$ is Fredholm. By Proposition \ref{pr1.4} its
index is $\left(  -\right)  ^{n}\chi_{C^{\ast}\left(  \pi\right)  }$ or
$\left(  -\right)  ^{n}\chi_{C^{\ast}\left(  \pi\right)  }^{\ell f}.$

\section{\label{s2}de Rham theory}

We discuss a de Rham-type theorem for the $\mathcal{L}^{2}$ cochains of
manifolds of bounded geometry. The forms and cochains take values in a bundle
of modules over a $C^{\ast}$-algebra. This builds on a theorem of Pierre Pansu
\cite{pan1}, \cite[Ch. 4]{pan2}, in which the usual conclusion is strengthened
to bounded equivalence of the complexes. This means that both the maps and
homotopies involved are bounded in suitable norms. In essence, he shows that
the usual double complex proof \cite[Ch. II]{botu} works under suitable
bounded geometry assumptions. Key features of our generalization are that it
applies to weighted spaces, and that the resulting cochain equivalences are
spatially bounded in a sense to be defined below. Some knowledge of Pansu's
proof is necessary in order to understand the remainder of this section.

\subsection{\label{ss2.1}}

\begin{definition}
\label{de2.1}An open covering $U$ $=\left\{  U_{\alpha}|\alpha\in I\right\}  $
of a metric space $X$ is \textit{{uniform} }if\ \textit{ }

\begin{enumerate}
\item for some $\epsilon>0$ the sets $U_{\alpha}^{\epsilon}=\left\{  x\in
U_{\alpha}|d\left(  x,X-U_{\alpha}\right)  >\epsilon\right\}  $ cover $X;$\ 

\item each $U_{\alpha}$ intersects a bounded number of others;

\item the diameters of the $U_{\alpha}$ are bounded.
\end{enumerate}
\end{definition}

A uniform covering of a separable space is countable. In what follows we will
use only uniform coverings. A BG manifold has uniform covers by open metric
balls of arbitrarily small fixed radius \cite[Lemma 1.1.2]{shu}. The version
of the Poincar\'{e} lemma used by Pansu is valid for such coverings with
sufficiently small radius. This condition will sometimes be abbreviated
``small balls''.

Let $M$ be a BG Riemannian manifold. As in Section 1, let $V$ be a unitary
flat bundle of $A$-modules over $M,$ and $\Omega_{c}$ and $\Omega_{d}$ be the
unweighted compactly supported forms with values in $V.$ $\Omega_{d}$ has the
inner product $\left(  u,v\right)  _{d}=\left(  u,v\right)  +\left(
du,dv\right)  .$

Spaces of smooth forms define presheaves. For an open set $U\subset$ $M,$ let
$\Omega_{d}\left(  U\right)  $ be the space of restrictions of elements of
$\Omega_{d}$ to $U,$ and similarly for other spaces. If $W\subset U$, the
restriction map is $r_{UW}.$ We will sometimes write $u|_{W}$ for $r_{UW}u.$

Let $\mathcal{F}$ be a presheaf on $M.$ For an open cover $\mathcal{U=}%
\left\{  U_{\alpha}\right\}  ,$ a compactly supported \v{C}ech $j$-cochain
with coefficients in $\mathcal{F}$ is an antisymmetric function $c_{\beta}%
\in\mathcal{F}\left(  U_{\beta}\right)  $ of nonempty $\left(  j+1\right)
$-fold intersections $U_{\beta}=U_{\alpha_{0}}\cap\cdots\cap U_{\alpha_{j}},$
such that $\overline{\cup_{\beta}U_{\beta}}$ with $c_{\beta}\neq0$ is compact.
The group of $j$-cochains is $\check{C}_{c}^{j}\left(  \mathcal{U}%
;\mathcal{F}\right)  .$

Pansu's proof requires some small modifications to work in the context of
Hilbert modules. Norms must be derived from inner products. Let $\mathcal{U}$
be a uniform cover of $M.$ For $c,d\in\check{C}_{c}^{j}\left(  \mathcal{U}%
;\Omega_{d}\right)  ,$ let $\left(  c,d\right)  _{d}=\sum_{\beta}\left(
c_{\beta},d_{\beta}\right)  _{d},$ with norm $\left\Vert \left(  c,c\right)
_{d}\right\Vert _{C^{\ast}}^{1/2}.$ The $\mathcal{L}^{2}$ \v{C}ech cochains
with coefficients in $\Omega_{d},$ $\check{C}_{1}^{j}\left(  \mathcal{U}%
;\Omega_{d}\right)  $ are the completion of the compactly supported cochains
in this norm. (The subscript means $k=1.)$ These form a double complex with
bounded differentials.

A locally constant section $c$ of $V$ on an open set $U\subset M$ is one for
which $dc=0$. Therefore $\left(  c,e\right)  _{d}=\left(  c,e\right)  $ for
any section $e.$ We denote (by abuse of notation) the compactly supported
cochains with values in the locally constant sections by $\check{C}_{c}^{\ast
}\left(  \mathcal{U};V\right)  .$ These are exactly the kernel of the
differential $\check{C}_{c}^{\ast}\left(  \mathcal{U};\Omega_{d}^{0}\right)
\rightarrow\check{C}_{c}^{\ast}\left(  \mathcal{U};\Omega_{d}^{1}\right)  .$
The completion is $\check{C}_{1}^{\ast}\left(  \mathcal{U};V\right)  .$
Generalizing the result of Pansu,

\begin{theorem}
\label{th2.2}If $\mathcal{U}$ is a uniform cover by open balls of sufficiently
small radius, the inclusions of $\Omega_{d}$ and $\check{C}_{c}^{\ast}\left(
\mathcal{U};V\right)  $ into $\check{C}_{c}^{\ast}\left(  \mathcal{U}%
;\Omega_{d}\right)  $ are bounded homotopy equivalences. Therefore
$\bar{\Omega}_{d}$ is boundedly equivalent to $\check{C}_{1}^{\ast}\left(
\mathcal{U};V\right)  .$%
%TCIMACRO{\TeXButton{End Proof}{\endproof}}%
%BeginExpansion
\endproof
%EndExpansion

\end{theorem}

We will give some refinements of this theorem after formalizing several
aspects of the proof. The first is the notion of a global inner product
derived from a pointwise inner product. In the following Definition, one could
take integrability in the strong sense. However, the Riemann integral suffices
for our purposes. Functions differing on sets of measure 0 are identified.

\begin{definition}
\label{de2.3} An $A$-\textit{{Hilbert presheaf }}consists of the following: a
presheaf $\mathcal{E}$ of pre-Hilbert $A$-modules over $M$ with all
restriction maps surjective; a positive Borel measure $\mu$ on $M$; a family
of Hermitian pairings $\left\langle .,.\right\rangle _{U}$ on $\mathcal{E}%
\left(  U\right)  $ for $\ U\subset M$ open, with values integrable $A$-valued
functions on $M.$ We assume these properties:
\end{definition}

\begin{enumerate}
\item If $u,v\in\mathcal{E}\left(  U\right)  ,$ $\left(  u,v\right)  _{U}%
=\int_{M}\left\langle u,v\right\rangle _{U}d\mu.$

\item $\left\langle u,u\right\rangle _{U}\geq0$.

\item If $W\subset U$, $\left\langle u|_{W},v|_{W}\right\rangle _{W}=\chi
_{W}\left\langle u,v\right\rangle _{U}.$ $\chi_{W}$ is the characteristic
function of $W.$
\end{enumerate}

\noindent We will sometimes write $\left\langle .,.\right\rangle $ for
$\left\langle .,.\right\rangle _{M}.$ For $\mathcal{E}=\Omega_{d}$ we use
$\left\langle u,v\right\rangle _{d,U}\left(  x\right)  =\left\langle u\left(
x\right)  ,v\left(  x\right)  \right\rangle +\left\langle du\left(  x\right)
,dv\left(  x\right)  \right\rangle $ for $x\in U,$ and $0$ for $x\notin U.$

\v{C}ech cochains form presheaves. The restrictions are restrictions of
cochains to open sets with the induced coverings. For $\mathcal{E}=\check
{C}^{\ast}\left(  \mathcal{U};\Omega_{d}\right)  ,$ $\left\langle
c,d\right\rangle _{U}=\sum_{\beta}\left\langle c_{\beta},d_{\beta
}\right\rangle _{d,U},$ and similarly for other groups of \v{C}ech cochains.
In these examples $\mu$ is the Riemannian measure. We will also use weighted
measures. For simplicial cochains, to be introduced below, the measure is discrete.

The restrictions are bounded with norm $\leq1$, since for $W\subset U,$
$u\in\mathcal{E}\left(  U\right)  ,$
\begin{align*}
\left(  u|_{W},u|_{W}\right)  _{W}  &  =\int_{M}\left\langle u|_{W}%
,u|_{W}\right\rangle _{W}d\mu=\int_{M}\chi_{W}\left\langle u,u\right\rangle
_{U}d\mu\\
&  \leq\int_{M}\left\langle u,u\right\rangle _{U}d\mu=\left(  u,u\right)
_{U}.
\end{align*}

$\mathcal{E}$ satisfies the following half of the sheaf axiom.\smallskip\ 

\begin{enumerate}
\item[$\mathcal{S}$:] Let an open set $U=\cup_{\alpha}U_{\alpha}$ with the
$U_{\alpha}$ open. If $u\in\mathcal{E}\left(  U\right)  $ is such that the
restrictions $u|_{U_{\alpha}}=0$ for all $\alpha,$ then $u=0.$
\end{enumerate}

\noindent For\label{2.1}%

\begin{align}
\left(  u,u\right)  _{U}  &  =\int_{M}\left\langle u,u\right\rangle _{U}%
d\mu\leq\sum_{\alpha}\int_{M}\chi_{U_{\alpha}}\left\langle u,u\right\rangle
_{U}d\mu\\
&  =\sum_{\alpha}\int_{M}\left\langle u|_{U_{\alpha}},u|_{U_{\alpha}%
}\right\rangle _{U_{\alpha}}d\mu=\sum_{\alpha}\left(  u|_{U\alpha
},u|_{U_{\alpha}}\right)  _{U_{\alpha}}.\nonumber
\end{align}
\noindent Therefore $\left(  u,u\right)  _{U}=0.$

\noindent The $\mathcal{L}^{2}$-type spaces we are using don't satisfy the
existence clause.

\subsection{\textbf{\label{ss2.2}}}

The idea of a spatially bounded operator is implicit in the proof. This is
related to, but rather different from, the concept of finite propagation
developed by Higson and Roe \cite[Chs. 3, 4]{roe}. It is introduced here to
allow a uniform treatment of several different situations. Let $\mathcal{E}$
be any presheaf satisfying $\mathcal{S},$ and $u\in\mathcal{E}\left(
M\right)  .$ There is a largest open set $V$ on which $u$ restricts to 0. By
$\mathcal{S}$ it is the union of all open sets on which $u$ restricts to 0.
The support\textit{\ }of $u,$ $Supp\left(  u\right)  ,$ is the complement of
$V.$

\begin{lemma}
\label{le2.4}Let $\mathcal{E}$ be a Hilbert presheaf. Elements of
$\mathcal{E}\left(  M\right)  $ with disjoint supports are orthogonal.
\end{lemma}

%

%TCIMACRO{\TeXButton{Proof}{\proof}}%
%BeginExpansion
\proof
%EndExpansion
For an open set $U,$ let $J_{U}=\left\{  u\in\mathcal{E}\left(  M\right)
:\left\langle u,u\right\rangle _{U}=0\,\right\}  .$ We claim that
\[
0\rightarrow J_{U}\rightarrow\mathcal{E}\left(  M\right)  \overset{r_{MU}%
}{\rightarrow}\mathcal{E}\left(  U\right)  \rightarrow0
\]

is exact. $r_{MU}$ is surjective by hypothesis. If $u\in J_{U},$ $\left(
u|_{U},u|_{U}\right)  =\int_{M}\left\langle u,u\right\rangle _{U}d\mu=0,$ so
$u|_{U}=0.$ If $u|_{U}=0,$ $\int_{M}\left\langle u,u\right\rangle _{U}d\mu=0,$
so $\left\langle u,u\right\rangle _{U}=0.$

Suppose that $u$ and $v$ have disjoint supports. Write \textquotedblleft$^{c}%
$" for complements. $u|_{Supp\left(  u\right)  ^{c}}=0,$ so $\chi_{Supp\left(
u\right)  ^{c}}\left\langle u,u\right\rangle =0.$ Therefore $\left\langle
u,u\right\rangle =0$ on $Supp\left(  u\right)  ^{c}.$Similarly, $\left\langle
v,v\right\rangle =0$ on $Supp(v)^{c}.$ $Supp\left(  u\right)  ^{c}\cup
Supp\left(  v\right)  ^{c}=M,$ so%

\[
\left\Vert \left\langle u,v\right\rangle \right\Vert \leq\left\Vert
\left\langle u,u\right\rangle \right\Vert \left\Vert \left\langle
v,v\right\rangle \right\Vert =0,
\]

\noindent and $\left(  u,v\right)  =0$.%
%TCIMACRO{\TeXButton{Endproof}{\endproof}}%
%BeginExpansion
\endproof
%EndExpansion

We will denote by $\mathcal{B}\left(  \mathcal{E},\mathcal{F}\right)  $ the
space $\mathcal{B}\left(  \mathcal{E}\left(  M\right)  ,\mathcal{F}\left(
M\right)  \right)  $ of bounded $A$-module homomorphisms. These are
not\textit{ }necessarily presheaf homomorphisms.

\begin{definition}
\label{de2.5} Let $\mathcal{E},$ $\mathcal{F}$ be two presheaves of Hilbert
modules satisfying $S.$ $T\in B\left(  \mathcal{E},\mathcal{F}\right)  $ is
\textit{{spatially bounded}\ }if there exists $R>0$ such that for all
$u\in\mathcal{E}\left(  M\right)  ,$ $Supp\left(  Tu\right)  \subset
N_{R}\left(  Supp\left(  u\right)  \right)  $ (the closed $R$-neighborhood).
The infimum of such $R$ is the \textit{spatial bound\ of }$T,$ $SB\left(
T\right)  .$
\end{definition}

Presheaf homomorphisms have spatial bound 0. Some elementary facts:

\label{2.2}%
\begin{align}
SB\left(  ST\right)   &  \leq SB\left(  S\right)  +SB\left(  T\right)  ,\\
SB\left(  S+T\right)   &  \leq\max\left\{  SB\left(  S\right)  ,SB\left(
T\right)  \right\}  .\nonumber
\end{align}
The completion $\overline{\mathcal{E}}$ of a Hilbert presheaf $\mathcal{E}$ is
formed by completing all the $\overline{\mathcal{E}}\left(  U\right)  .$ The
restrictions extend by continuity. $\overline{\mathcal{E}}$ is a presheaf of
Hilbert modules, but not a Hilbert presheaf in general. The restrictions may
not be surjective. There are difficulties involved in extending the pairing
$\left\langle \cdot,\cdot\right\rangle .$ $\overline{\mathcal{E}}$ satisfies
$\mathcal{S}$ because $\left(  \text{\ref{2.1}}\right)  $ holds in
$\overline{\mathcal{E}}$ by continuity. To relate completion and spatial
boundedness we must make an assumption.

\begin{enumerate}
\item[ ] $\mathcal{A}$: Any $u\in\overline{\mathcal{E}\left(  M\right)  }$ is
the limit of elements of $\mathcal{E}\left(  M\right)  $ with support in
$N_{\epsilon}\left(  Supp\left(  u\right)  \right)  $ for any $\epsilon
>0.$\smallskip\ 
\end{enumerate}

This condition holds for the relevant examples. For $\Omega_{d}$ we prove a
relative version. Let $U\subset M$ be open and $u\in\bar{\Omega}_{d}\left(
U\right)  .$ By definition, $u$ is the limit of a sequence $\left(
u_{n}\right)  $ of restrictions of elements $v_{n}$ of $\Omega_{d}$ to $U.$
Let $\psi\in C_{b}^{\infty,1}\left(  M\right)  $ be 1 on $Supp\left(
u\right)  $ and 0 on $M-N_{\epsilon}\left(  Supp\left(  u\right)  \right)  .$
Then $\psi u_{n}\in\Omega_{d}\left(  U\right)  $ since it is the restriction
of $\psi v_{n}.$ By Lemma \ref{le1.1}, $\psi u_{n}\rightarrow\psi u=u$ in
$\bar{\Omega}_{d}\left(  U\right)  .$

Let $c\in\check{C}_{1}^{j}\left(  \mathcal{U};V\right)  ,$ and $c_{n}\in
\check{C}_{c}^{j}\left(  \mathcal{U};V\right)  $ such that $c_{n}\rightarrow
c.$ For each $\beta,$ $c_{n\beta}\longrightarrow c_{\beta}.$ Since
$dc_{n\beta}=0,$ $dc_{\beta}=0,$ so $c_{\beta}$ is smooth. If $\beta_{i}$ are
an enumeration of the $\beta,$ $\sum_{i=i}^{N}c_{\beta_{i}}\longrightarrow
c$\ on $Supp\left(  c\right)  .$

Let $c\in\check{C}_{1}^{j}\left(  \mathcal{U};\Omega_{d}\right)  $, which is
the Hilbert sum $\bigoplus_{\beta}\bar{\Omega}_{d}\left(  U_{\beta}\right)  .$
For any $\epsilon>0$ and each $\beta$ there is a sequence $c_{n\beta}$ in
$\Omega_{d}\left(  U_{\beta}\right)  $ with supports in $N_{\epsilon}\left(
Supp\left(  c_{\beta}\right)  \right)  $ such that $c_{n\beta}\rightarrow
c_{\beta}.$ By passing to subsequences we obtain $c_{n}^{\prime}$ with support
in $N_{\epsilon}\left(  Supp\left(  c\right)  \right)  $ such that
$c_{n}^{\prime}\rightarrow c.$ Let $c_{n}^{^{\prime\prime}}$ be some
truncation of $c_{n}^{^{\prime}}$ with finitely many nonzero $c_{n\beta
}^{^{\prime\prime}}$ such that $\left\Vert c_{n}^{\prime\prime}-c_{n}^{\prime
}\right\Vert <1/2^{n}.$ Then the $c_{n}^{\prime\prime}\in\check{C}_{c}%
^{j}\left(  \mathcal{U};\Omega_{d}\right)  ,$ have supports in $N_{\epsilon
}\left(  Supp\left(  c\right)  \right)  ,$ and converge to $c.$

\begin{lemma}
\label{le2.6}Let $\mathcal{E},$ $\mathcal{F}$ be Hilbert presheaves satisfying
condition $\mathcal{A},$ and $T\in\mathcal{B}\left(  \mathcal{E}%
,\mathcal{F}\right)  $ have spatial bound $R.$ Then $T$ extends to an element
$\bar{T}$ of $\mathcal{B}\left(  \overline{\mathcal{E}},\mathcal{\bar{F}%
}\right)  $ with spatial bound $R.$
\end{lemma}

%

%TCIMACRO{\TeXButton{Proof}{\proof}}%
%BeginExpansion
\proof
%EndExpansion
Choose $u_{n}$ in $\mathcal{E}\left(  M\right)  $ converging to $u$ in some
$N_{\epsilon}\left(  Supp\left(  u\right)  \right)  .$ Then $Supp\left(
Tu_{n}\right)  \allowbreak\subset N_{R}\left(  Supp\left(  u_{n}\right)
\right)  \subset N_{R+\epsilon}\left(  Supp\left(  u\right)  \right)  .$
Therefore $Tu_{n}$ restricts to 0 on the complement of $N_{R+\epsilon}\left(
Supp\left(  u\right)  \right)  .$ By continuity of the restrictions, the same
is true of $\bar{T}u.$ Therefore $Supp\left(  \bar{T}u\right)  \subset
N_{R+\epsilon}\left(  Supp\left(  u\right)  \right)  .$ Since $\epsilon$ is
arbitrary, $Supp\left(  \bar{T}u\right)  \subset N_{R}\left(  Supp\left(
u\right)  \right)  .$%
%TCIMACRO{\TeXButton{End Proof}{\endproof}}%
%BeginExpansion
\endproof
%EndExpansion

For example, the exterior derivative and multiplication by a smooth function
on $\bar{\Omega}_{d}$ have spatial bound 0, since this is evidently the case
on $\Omega_{c}.$

We will also need a fineness assumption. The support of a set of elements is
defined to be the union of their supports. We assume that there exists a
sequence $\left\{  S_{i}\right\}  \subset B\left(  \mathcal{E}\right)  $ of
operators with spatial bound 0 such that each $Supp\left(  \operatorname{Im}%
\left(  S_{i}\right)  \right)  $ is compact and $\sum S_{i}$ converges
strongly to the identity. It will be seen at the end of Section 2.3 that this
is satisfied by the relevant examples. Let $\mathcal{E}$ be a Hilbert presheaf
satisfying this and $\mathcal{A}.$

\begin{lemma}
\label{le2.7}Elements $u,v$ of $\overline{\mathcal{E}}\left(  M\right)  $ with
disjoint supports are orthogonal.
\end{lemma}

%

%TCIMACRO{\TeXButton{Proof}{\proof} }%
%BeginExpansion
\proof
%EndExpansion
Suppose first that $u$ and $v$ have compact supports. For some $\epsilon>0$
there are disjoint $\epsilon$-neighborhoods $U$ and $V$ of $Supp(u)$ and
$Supp\left(  v\right)  .$ Choose elements $u_{n}$ of $\mathcal{E}$ with
supports in $U$ converging to $u,$ and similarly $v_{n}$ converging to $v$ in
$V$. Then $\left(  u,v\right)  =$ $\lim\left(  u_{n},v_{n}\right)  =\lim0=0. $

For the general case, by Lemma \ref{le2.6}, the $S_{i}$ extend to $\bar{S}%
_{i}\in\mathcal{B}\left(  \overline{\mathcal{E}}\right)  $ with spatial bound
$0.$ Therefore the elements $\bar{S}_{i}u$ and $\bar{S}_{i}v$ have compact
supports , and $\left(  \bar{S}_{i}u,\bar{S}_{j}v\right)  =0$ for all $i$ and
$j.$ Then $\left(  u,v\right)  =\lim_{k\rightarrow\infty}\left(  \sum
_{i=1}^{k}\bar{S}_{i}u,\sum_{i=1}^{k}\bar{S}_{i}v\right)  =0.$%
%TCIMACRO{\TeXButton{End Proof}{\endproof}}%
%BeginExpansion
\endproof
%EndExpansion

\subsection{\textbf{\label{ss2.3}}}

We now discuss the algebraic basis for applications of spatial boundedness.
Let $P$ and $Q$ be pre-Hilbert modules. Let $I$ be a countable index set. We
make the following assumptions:

\begin{enumerate}
\item For $i\in I$ there are operators $S_{i}\in\mathcal{B}\left(  P\right)  $
such that

\begin{enumerate}
\item The number of $k$ such that for a given $i,$ $\operatorname{Im}S_{i}$ is
not orthogonal to $\operatorname{Im}S_{k}$ is uniformly bounded.

\item For all $u,$ $S_{i}u=0$ except for finitely many $i$.

\item For any subset $J\subset I,$ the operator $\sum_{j\in J}S_{j}$ is bounded.
\end{enumerate}

\item There are uniformly bounded operators $T_{ji}$ with domains
$\operatorname{Im}S_{i}$ and ranges in $Q$ such that

\begin{enumerate}
\item The number of $T_{ji}$ for a given $i$ is uniformly bounded.

\item The number of pairs $\left(  \ell,k\right)  $ such that for a given
$\left(  j,i\right)  $, $\operatorname{Im}T_{ji}$ is not orthogonal to
$\operatorname{Im}T_{\ell k}$ is uniformly bounded.
\end{enumerate}
\end{enumerate}

\noindent In 1c the operator is a finite sum for each element of $P,$ so order
is irrelevant and the sum converges strongly.

The prototypical case is when $P=\bigoplus_{i}P_{i}$ and $Q=\bigoplus_{j}%
Q_{j}$ are orthogonal sums. Let $\left[  R_{ji}\right]  $ be a uniformly
bounded matrix of operators such that the number of nonzero elements in any
row or column is bounded. Let $p_{i}$ and $q_{j}$ be the projections and
inclusions. Then the matrix operator is $\sum_{i,j}T_{ji}S_{i}$ with
$S_{i}=p_{i}$ and $T_{ji}=q_{j}R_{ji}.$ This case is due to Higson and Roe.
The general case is needed to deal with partitions of unity.

We will make use of the following theorem of Paschke \cite[Theorem 2.8]{pas}:
a $\mathbb{C}$-linear mapping $T$ between pre-Hilbert modules is a bounded
$A$-module homomorphism if and only if there exists $K>0$ such that $\left(
Tu,Tu\right)  <K^{2}\left(  u,u\right)  $ for all $u,$ in which case $\left\|
T\right\|  \le K.$

\begin{proposition}
\label{pr2.8}$\sum_{i,j}T_{ji}S_{i}$ extends to an element of $\mathcal{B}%
\left(  \bar{P},\bar{Q}\right)  .$
\end{proposition}

%

%TCIMACRO{\TeXButton{Proof}{\proof} }%
%BeginExpansion
\proof
%EndExpansion
Let $T_{i}=\sum_{j}T_{ji}.$ Then the $\left\Vert T_{i}\right\Vert $ are
uniformly bounded, say by $K,$ and the number of $k$ such that for a given
$i,$ $\operatorname{Im}T_{i}$ and $\operatorname{Im}T_{k}$ are not orthogonal
is uniformly bounded. We may construct inductively a partition of $I$ into
finitely many disjoint sets $I_{\ell}$ such that if $i,j\in I_{\ell},$ $i\neq
j,$ then $\operatorname{Im}S_{i}\perp\operatorname{Im}S_{j}$ and
$\operatorname{Im}T_{i}\perp\operatorname{Im}T_{j}.$ It then suffices to show
that $\sum_{i\in I_{\ell}}T_{i}S_{i}$ is bounded for each $\ell.$ Taking all
summations over $I_{\ell},$%
\begin{align*}
\left(  \left(  \sum T_{i}S_{i}\right)  u,\left(  \sum T_{i}S_{i}\right)
u\right)   &  =\sum\left(  T_{i}S_{i}u,T_{i}S_{i}u\right)  \leq\sum\left\Vert
T_{i}\right\Vert ^{2}\left(  S_{i}u,S_{i}u\right) \\
&  \leq K^{2}\sum\left(  S_{i}u,S_{i}u\right)  =K^{2}\left(  \left(  \sum
S_{i}\right)  u,\left(  \sum S_{i}\right)  u\right) \\
&  \leq K^{2}L^{2}\left(  u,u\right)
\end{align*}

for some $L,$ by assumption.%
%TCIMACRO{\TeXButton{End Proof}{\endproof}}%
%BeginExpansion
\endproof
%EndExpansion

\noindent In the matrix case passage to subsets isn't required.

The following is a geometrical version of the previous proposition. Let
$\mathcal{E}$ and $\mathcal{F}$ be Hilbert presheaves. We assume that
$\mathcal{E}$ satisfies condition $\mathcal{A}$ as well as the following.

\begin{enumerate}
\item[(I)] $\mathcal{E}\left(  M\right)  $ consists of elements with compact support.

\item[(II)] There is a countable set $\left\{  S_{i}\right\}  _{i\in I}%
\subset\mathcal{B}\left(  \mathcal{E}\right)  $ such that

\begin{enumerate}
\item The $S_{i}$ have spatial bound 0.

\item The diameters of the $Supp\left(  \operatorname{Im}S_{i}\right)  $ are
uniformly bounded.

\item The set $\left\{  Supp\left(  \operatorname{Im}S_{i}\right)  \right\}  $
is uniformly locally finite. This means that for any $r>0$ there is an $n_{r}$
such that every ball of radius $r$ intersects no more than $n_{r}$ elements.

\item For any subset $J\subset I,$ $\sum_{j\in J}S_{j}\in\mathcal{B}\left(
\mathcal{E}\right)  .$
\end{enumerate}

\item[(III)] There are uniformly bounded operators $T_{ji}$ with domains
$\operatorname{Im}S_{i}$ such that

\begin{enumerate}
\item The number of $T_{ji}$ for a given $i$ is uniformly bounded.

\item Each $T_{ji}$ has spatial bound $\le R.$
\end{enumerate}
\end{enumerate}

\begin{proposition}
\label{pr2.9}$\sum_{i,j}T_{ji}S_{i}$ has an extension $\bar{T}\in
\mathcal{B}\left(  \overline{\mathcal{E}},\mathcal{\bar{F}}\right)  .$ If in
addition $\sum_{i}S_{i}=I,$ $\bar{T}$ has spatial bound $\leq R$.
\end{proposition}

%

%TCIMACRO{\TeXButton{Proof}{\proof}}%
%BeginExpansion
\proof
%EndExpansion
We check the hypotheses of Proposition \ref{pr2.8}. (1a) follows from (IIb,c)
since elements with disjoint supports are orthogonal. (1b) follows from (I)
and (IIa,c); (1c) from (IId) and (2a) from (IIIa). (IIb,c) and (IIIa,b) imply
that the diameters of the $Supp\left(  \operatorname{Im}T_{ji}\right)  $ are
uniformly bounded, and that the $Supp\left(  \operatorname{Im}T_{ji}\right)  $
are uniformly locally finite. Thus (2b) holds. Therefore $\sum_{j,i}%
T_{ji}S_{i}$ extends to $\overline{\mathcal{E}}$. By $\left(  \text{\ref{2.2}%
}\right)  $, each $T_{ji}S_{i}$ has spatial bound $\leq R,$ so that
$\sum_{j,i}T_{i}S_{i}$ has spatial bound $\leq R.$ Spatial boundedness of
$\bar{T}$ follows from Lemma \ref{le2.6}.
%TCIMACRO{\TeXButton{End Proof}{\endproof}}%
%BeginExpansion
\endproof
%EndExpansion

We now apply the above material to sharpen Theorem \ref{th2.2}. It is first
necessary to establish the boundedness and spatial boundedness of the maps and
homotopies occurring in the proof, at the level of compactly supported
cochains or smooth forms. This requires applications of Proposition
\ref{pr2.9} in several different contexts depending on $\mathcal{E}.$

Let $\mathcal{E=}$ $\Omega_{d}.$ Any uniform cover admits a uniformly bounded
partition of unity $\left\{  \phi_{i}\right\}  \subset C_{b}^{\infty,1}\left(
M\right)  $ \cite[Lemma 1.1.3]{shu}. We take $S_{i}=\phi_{i}.$ The conditions
on the $S_{i}$ are then clear$.$ As an example, the map $r:\Omega
_{d}\rightarrow\check{C}_{c}^{0}\left(  \mathcal{U};\Omega_{d}\right)  $ is
given by $\sum_{\beta}r_{MU_{\beta}}.$ Let $T_{\beta i}=$ $r_{MU_{\beta}%
}|\operatorname{Im}\phi_{i}.$ Since the $r_{MU_{\beta}}$ and $\phi_{i}$ have
spatial bounds 0, these do too. Since the $r_{MU_{\beta}}$ have norm 1, they
are uniformly bounded. Then $r=\sum_{i,\beta}T_{\beta i}\phi_{i}$ extends to
$\bar{\Omega}_{d}$ with spatial bound $0.$

The \v{C}ech groups $\check{C}_{c}^{\ast}\left(  \mathcal{U};\Omega
_{d}\right)  $ and $\check{C}_{c}^{\ast}\left(  \mathcal{U};V\right)  $ are
orthogonal sums by definition. The $S_{\beta}$ are the projections on the
$\Omega_{d}\left(  U_{\beta}\right)  .$ The boundedness and spatial
boundedness of maps with source a \v{C}ech group can be established as in the
example above from the corresponding facts about their components. The latter
are evident for the maps involved in the de Rham equivalence.

The additional hypothesis in Proposition 2.9 is satisfied in our examples.
$\sum_{i=1}^{n}S_{i}$ is the identity on elements with support in any compact
set for large enough $n.$

We conclude the following. Let $\mathcal{U}$ be a uniform covering by small balls.

\begin{theorem}
\label{th2.10}The de Rham equivalence between $\Omega_{d}$ and $\check{C}%
_{c}^{\ast}\left(  \mathcal{U};V\right)  $ is bounded and spatially bounded.
It therefore extends to an equivalence between $\bar{\Omega}_{d}$ and
$\check{C}_{1}^{\ast}\left(  \mathcal{U};V\right)  $ with the same properties.%
%TCIMACRO{\TeXButton{End Proof}{\endproof}}%
%BeginExpansion
\endproof
%EndExpansion

\end{theorem}

\subsection{\label{ss2.4}}

We will show that, under the assumption of spatial boundedness, operators on
elements with compact support give rise to operators between weighted spaces.
The analytic weighted spaces of forms have already been defined using the
weight functions $\tau\left(  x\right)  =k^{\rho\left(  x\right)  }$. The
definition extends immediately to define $\mathcal{E}_{k}$ for any Hilbert
presheaf $\mathcal{E}$. Let $\mathcal{E}$ and $\mathcal{F}$ be Hilbert presheaves.

\begin{lemma}
\label{le2.11}Let $T\in\mathcal{B}\left(  \mathcal{E},\mathcal{F}\right)  $
have spatial bound $R.$ For any $r>0,$ $T$ is bounded in any $k$-norm on
elements of $\mathcal{E}$ with support of diameter $\leq r.$
\end{lemma}

%

%TCIMACRO{\TeXButton{Proof}{\proof}}%
%BeginExpansion
\proof
%EndExpansion
Let $u$ have support of diameter $\leq r.$ Write $V=Supp\left(  u\right)  .$
Let $g_{V}=\max_{x\in V}\tau\left(  x\right)  ,$ $\ell_{V}=\min_{x\in V}%
\tau\left(  x\right)  .$ It is clear that $\ell_{V}\left\Vert u\right\Vert
\leq\left\Vert u\right\Vert _{k}\leq g_{V}\left\Vert u\right\Vert .$ Let
$\ell=\ell_{V}=\tau\left(  b\right)  $, $g=g_{N_{R}\left(  V\right)  }%
=\tau\left(  a\right)  .$ Then $d\left(  a,b\right)  \leq r+R.$ If $C$ is a
Lipschitz constant for $\rho,$ $\rho\left(  a\right)  -\rho\left(  b\right)
<C\left(  r+R\right)  .$ It follows that $\dfrac{g}{\ell}$ is uniformly
bounded for all such $u.$ Since $Supp\left(  Tu\right)  \subset N_{R}\left(
Supp\left(  u\right)  \right)  ,$ $\left\Vert Tu\right\Vert _{k}\leq
g\left\Vert Tu\right\Vert \leq g\left\Vert T\right\Vert \left\Vert
u\right\Vert \leq\dfrac{g}{\ell}\left\Vert T\right\Vert \left\Vert
u\right\Vert _{k}.$%
%TCIMACRO{\TeXButton{End Proof}{\endproof}}%
%BeginExpansion
\endproof
%EndExpansion

The next result is a variant of Proposition \ref{pr2.9}.

\begin{proposition}
\label{pr2.12}Assume the hypotheses of Proposition \ref{pr2.9} except for
(III). In addition, suppose that $\sum S_{i}=I$. Let $T\in\mathcal{B}\left(
\mathcal{E},\mathcal{F}\right)  $ have spatial bound $R.$ Then $T$ has an
extension in $\mathcal{B}\left(  \overline{\mathcal{E}}_{k},\mathcal{\bar{F}%
}_{k}\right)  $ which has spatial bound $\leq R.$
\end{proposition}

%

%TCIMACRO{\TeXButton{Proof}{\proof}}%
%BeginExpansion
\proof
%EndExpansion
Let $T_{i}=T_{ii}=T|\operatorname{Im}S_{i}.$ Point (III) is replaced by the
above lemma, and by hypothesis. Thus $\sum T_{i}S_{i}$ is bounded in the
$k$-norms. But%

\[
\left(  \sum T_{i}S_{i}\right)  u=\sum T_{i}S_{i}u=\sum TS_{i}u=T\sum
S_{i}u=Tu.
\]

\noindent Therefore $T$ extends to $\overline{\mathcal{E}}_{k}$. Spatial
boundedness follows from Lemma \ref{le2.6}.%
%TCIMACRO{\TeXButton{End Proof}{\endproof}}%
%BeginExpansion
\endproof
%EndExpansion

Using this Proposition and Theorem \ref{th2.2},

\begin{theorem}
\label{th2.13}The de Rham equivalence extends to a bounded and spatially
bounded equivalence between $\bar{\Omega}_{d,k}$ and $\check{C}_{k}^{\ast
}\left(  \mathcal{U};V\right)  ,$ for $\mathcal{U}$ a uniform cover by small
balls.%
%TCIMACRO{\TeXButton{End Proof}{\endproof}}%
%BeginExpansion
\endproof
%EndExpansion

\end{theorem}

\subsection{\textbf{\label{ss2.5}}}

For our purposes it is convenient to work with simplicial rather than \v{C}ech
cochains. Let $K\rightarrow M$ be a smooth triangulation. Let $C_{c}^{\ast
}\left(  K;V\right)  $ be the compactly supported cochains of $K$ with local
coefficients in $V$ \cite[Sections 30, 31]{ste}. It is a right $A$-module. Let
the $j$-simplexes of $K$ be $\left\{  \sigma_{i}\right\}  .$ We view the
$j$-cochain associated to $\sigma_{i}$ as being localized at the barycenter
$x_{i}\in\sigma_{i}.$ Then $C_{c}^{j}\left(  K;V\right)  \cong\bigoplus
_{i}V_{x_{i}}.$ For $e,f\in C_{c}^{j}\left(  K;V\right)  ,$ $\left(
e,f\right)  =\sum_{i}\left\langle e\left(  x_{i}\right)  ,f\left(
x_{i}\right)  \right\rangle .$ $\left\langle \cdot,\cdot\right\rangle $
denotes the fiber inner products. More generally, $\left(  e,f\right)
_{k}=\sum_{i}\left\langle e\left(  x_{i}\right)  ,f\left(  x_{i}\right)
\right\rangle k^{2\rho\left(  x_{i}\right)  }$. The weighted $\mathcal{L}^{2}$
simplicial cochains $C_{k}^{\ast}\left(  K;V\right)  $ are the completions of
$C_{c}^{\ast}\left(  K;V\right)  $ with respect to these inner products.

$C_{c}^{j}\left(  K;V\right)  $ gives rise to a Hilbert presheaf. The group of
sections over $U$ is defined to be $\left\{  \bigoplus_{i}V_{x_{i}}|x_{i}\in
U\right\}  ,$ with $r_{MU}$ the corresponding projection. The pointwise inner
product $\left\langle e,f\right\rangle _{U}\left(  x_{i}\right)  =\left\langle
e\left(  x_{i}\right)  ,f\left(  x_{i}\right)  \right\rangle $ if $x_{i}\in
U,$ $0$ otherwise. The measure $\mu$ is the counting measure on $\left\{
x_{i}\right\}  .$ Condition $\mathcal{A}$ holds. The proof is similar to that
for $\check{C}_{c}^{\ast}\left(  \mathcal{U};V\right)  $ in \ref{ss2.2}.

A homeomorphism $h:X\rightarrow Y$ of metric spaces is a quasi-isometry if
there exists $C>1$ such that for all $x\in X,$ $\dfrac{1}{C}d\left(
x,y\right)  <d\left(  h\left(  x\right)  ,h\left(  y\right)  \right)
<Cd\left(  x,y\right)  .$

\begin{definition}
\begin{enumerate}
\item \label{de2.14} A \textit{bounded geometry }(BG)\textit{ simplicial
complex }is one in which each vertex is a face of a uniformly bounded number
of simplexes.\textit{ }

\item A BG \textit{triangulation} of $M$ is a smooth triangulation
$K\rightarrow M$ by a BG simplicial complex which is a quasi-isometry when $K$
is equipped with the path metric for which each simplex has the standard metric.
\end{enumerate}
\end{definition}

\noindent The idea is that all images of simplexes of $K$ of the same
dimension have approximately the same size and shape. BG triangulations
clearly admit BG subdivisions of arbitrarily small mesh. The existence of BG
triangulations of BG manifolds is sometimes referred to as an unpublished
result of Calabi. However no detailed proof has ever been published. It must
be considered to be an open question. We will make use of BG triangulations
only in cases where they may be constructed \textquotedblleft by
hand\textquotedblright.

The condition (1) implies that the differentials of $C_{1}^{\ast}\left(
K;V\right)  $ are bounded. Those of $C_{k}^{\ast}\left(  K;V\right)  $ are
then bounded by Proposition \ref{pr2.12}.

Let $K\rightarrow M$ be a BG triangulation and $\mathcal{V}$ the cover of $M$
by the open vertex stars of $K.$ It is uniform.

\begin{lemma}
\label{le2.15}There are bounded and spatially bounded isomorphisms\newline%
$C_{k}^{\ast}\left(  K;V\right)  \rightarrow\check{C}_{k}^{\ast}\left(
\mathcal{V};V\right)  .$
\end{lemma}

%

%TCIMACRO{\TeXButton{Proof}{\proof}}%
%BeginExpansion
\proof
%EndExpansion
The map is induced by a bijection between the $j$-simplexes of $K$ and the
$\left(  j+1\right)  $-fold intersections of the vertex stars. For a vertex
$y_{\alpha}$ let $U_{\alpha}$ be its star. A simplex $\sigma_{\beta}=\left\{
y_{\alpha_{0}},\cdots,y_{\alpha_{j}}\right\}  $ then corresponds to $U_{\beta
}.$ The value of a cochain in $V_{x_{\beta}}$ determines a locally constant
section over $U_{\beta}$ by parallel transport. This gives an isomorphism
$C_{c}^{j}\left(  K;V\right)  \rightarrow\check{C}_{c}^{j}\left(
\mathcal{V};V\right)  .$ It is clearly spatially bounded. The bounded geometry
condition implies that there are only a finite number of combinatorial types
of vertex stars and of their $\left(  j+1\right)  $-fold intersections. Since
the triangulation is a quasi-isometry, the volumes in $M$ of the $U_{\beta}$
are uniformly bounded above and below. Let $c\in\check{C}_{c}^{j}\left(
\mathcal{V};V\right)  .$ We noted previously that $\left(  c,c\right)
_{d}=\left(  c,c\right)  .$ For any $\beta,$ by compatibility of the
connection, $d\left\langle c_{\beta},c_{\beta}\right\rangle =0.$ Since
$U_{\beta}$ is connected, $\left\langle c_{\beta},c_{\beta}\right\rangle $ is
constant, so $\left(  c_{\beta},c_{\beta}\right)  =\left\langle c\left(
x\right)  ,c\left(  x\right)  \right\rangle Vol\left(  U_{\beta}\right)  $ for
any $x\in U_{\beta}.$ Therefore for some $C>0$ and all $\beta,$ $\dfrac{1}%
{C}\left(  c_{\beta},c_{\beta}\right)  <\left\langle c\left(  x_{\beta
}\right)  ,c\left(  x_{\beta}\right)  \right\rangle <C\left(  c_{\beta
},c_{\beta}\right)  ,$ and the groups are boundedly isomorphic. The
equivalence in the $k$-norms is an application of Proposition \ref{pr2.12}. In
the simplicial groups we take the $S_{i}$ to be the projections of
$\bigoplus_{i}V_{x_{i}}$ onto its summands.%
%TCIMACRO{\TeXButton{End Proof}{\endproof}}%
%BeginExpansion
\endproof
%EndExpansion

\begin{remark}
\label{re2.16}This proof illustrates a general principle. Because of the
finiteness of the combinatorial types of vertex stars in a BG simplicial
complex, any construction on vertex stars depending only on the combinatorial
structure involves a bounded number of choices. Since a BG triangulation is a
quasi-isometry, local operators on $M$ produced by such a construction will be
uniformly bounded and uniformly spatially bounded.
\end{remark}

\begin{theorem}
\label{th2.17}If $K$ is a BG triangulation of $M,$ then for every $k,$
$C_{k}^{\ast}\left(  K;V\right)  $ is boundedly equivalent to $\bar{\Omega
}_{d,k}$ by a spatially bounded equivalence.
\end{theorem}

%

%TCIMACRO{\TeXButton{Proof}{\proof}}%
%BeginExpansion
\proof
%EndExpansion

Let $\mathcal{V}$ be as above. Any uniform cover has a uniform refinement by
small balls. Let $\mathcal{U}$ be such a refinement of $\mathcal{V}.$ We will
show that any function $\alpha\rightarrow s\left(  \alpha\right)  $ with
$U_{\alpha}\subset V_{s\left(  \alpha\right)  }$ induces a bounded and
spatially bounded equivalence $\check{C}_{c}^{\ast}\left(  \mathcal{V}%
;V\right)  \rightarrow\check{C}_{c}^{\ast}\left(  \mathcal{U};V\right)  .$ In
light of Theorem \ref{th2.13} and Proposition \ref{pr2.12}, this will complete
the proof. Any refining map $\mathcal{U}^{\prime}\rightarrow\mathcal{U}$ of
uniform covers induces a bounded and spatially bounded map of double complexes
$\check{C}_{c}^{\ast}\left(  \mathcal{U};\Omega_{d}\right)  \rightarrow
\check{C}_{c}^{\ast}\left(  \mathcal{U}^{\prime};\Omega_{d}\right)  .$ This is
an application of Proposition \ref{pr2.9}. The $T_{\gamma\beta}$ are the
restrictions induced by the $U_{\beta}^{\prime}\rightarrow U_{\gamma}.$ The
conditions are evident.

We choose covers as follows: Let $K^{\prime}$ be a BG subdivision of $K$ so
that the associated $\mathcal{V}^{\prime}$ refines $\mathcal{U}.$ Let
$\mathcal{U}^{\prime}$ be a uniform refinement of $\mathcal{V}^{\prime}$ by
small balls. We thus have refinements
\[
\mathcal{U}^{\prime}\rightarrow\mathcal{V}^{\prime}\rightarrow
\mathcal{U\rightarrow V}.
\]

The maps of $\Omega_{d}$ and $\check{C}_{c}^{\ast}\left(  \mathcal{\cdot
};V\right)  $ into $\check{C}_{c}^{\ast}\left(  \mathcal{\cdot};\Omega
_{d}\right)  $ are natural under refinement. Using Theorem \ref{th2.10}, they
are bounded and spatially bounded equivalences for $\mathcal{U}$ and
$\mathcal{U}^{\prime}.$ The same is then true of $\check{C}_{c}^{\ast}\left(
\mathcal{U};V\right)  \rightarrow\check{C}_{c}^{\ast}\left(  \mathcal{U}%
^{\prime};V\right)  $. Refinement induces $\check{C}_{c}^{\ast}\left(
\mathcal{V};V\right)  \rightarrow\check{C}_{c}^{\ast}\left(  \mathcal{V}%
^{\prime};V\right)  .$ A homotopy inverse is induced from any standard
subdivision map on simplicial cochains \cite[Ch. IV]{lef}. The $T_{\gamma
\beta}$ for Proposition \ref{pr2.8} are the matrix coefficients of the maps
and homotopies. This uses Remark \ref{re2.16}.

The equivalence of $\check{C}_{c}^{*}\left(  \mathcal{V};V\right)  $ and
$\check{C}_{c}^{*}\left(  \mathcal{U};V\right)  $ now follows from a general
fact: in any category, if there are morphisms
\[
C\overset{f}{\rightarrow}D\overset{g}{\rightarrow}E\overset{h}{\rightarrow}F
\]

\noindent with $gf$ and $hg$ equivalences, then $f$ is an equivalence.%
%TCIMACRO{\TeXButton{End Proof}{\endproof}}%
%BeginExpansion
\endproof
%EndExpansion

In the next section it will be clearer to work with chains than cochains. Let
$C_{j}\left(  K;V\right)  $ be the local coefficient chains. These are finite
sums $\sum_{i}c_{i}\sigma_{i},$ with $c_{i}\in V_{x_{i}.}$ The $k$-inner
product is $\left(  c,d\right)  _{k}=\sum_{i}\left\langle c_{i},d_{i}%
\right\rangle k^{2\rho\left(  x_{i}\right)  }.$ The completions are $C_{j}%
^{k}\left(  K;V\right)  .$ For a BG triangulation, there is a bounded and
spatially bounded equivalence (up to sign) $C_{k}^{\ast}\left(  K;V\right)
\rightarrow C_{n-\ast}^{k}\left(  K;V\right)  .$ This follows our standard
pattern and uses Remark \ref{re2.16}: the maps occurring in Poincar\'{e}
duality are locally defined with a bounded number of choices in each vertex star.

We will also use the ordinary de Rham theorems for simplicial cochains and
compactly supported simplicial cochains, with coefficients in $V.$ The proof
in \cite[Ch. IV]{whi} adapts readily.

\section{\label{s3}Homology of mapping telescopes}

In this section we establish the equivalences between weighted forms and
ordinary cochain complexes on certain manifolds of bounded geometry, as stated
in Theorem \ref{th0.3}.

\subsection{\textbf{\label{ss3.1}}}

We construct an infinite cyclic covering associated to an end.\textbf{\ }Let
$M$ be a complete connected Riemannian manifold with finitely many ends.
Suppose that there exists a cocompact open neighborhood $U$ of one of the ends
and a proper smooth embedding $h:U\rightarrow U$ such that $\bigcap_{n}%
h^{n}U=\emptyset.$ Let $\bigcup_{n=1}^{\infty}U_{n}$ be the disjoint union of
copies of $U.$ Let $\bar{N}=\bigcup_{n=1}^{\infty}U_{n}/\left\{  x_{n}%
\sim\left(  hx\right)  _{n+1}\right\}  .$ This is a smooth manifold with two
ends. The map $z$ defined by $z\left[  x_{n}\right]  =\left[  \left(
hx\right)  _{n}\right]  $ is a diffeomorphism, and extends to a properly
discontinuous action of $\mathbb{Z}$ by letting $z^{-1}\left[  x_{n}\right]
=\left[  x_{n+1}\right]  .$ Let $N$ be the quotient. By \cite[Theorem
13.11]{hura} there exist closed cocompact connected neighborhoods $\bar{N}%
^{+}$ and $\bar{N}^{-}$ of the ends of $\bar{N}$ with the following
properties: $\bar{N}=\bar{N}^{+}\cup\bar{N}^{-}$, $\bar{N}^{+}\cap\bar{N}%
^{-}=V_{0}$ is a closed codimension one submanifold, and $z\bar{N}^{+}%
\subset\bar{N}^{+}.$ Then $\bar{N}^{+}$ can be identified with a neighborhood
of the end of $M.$

We introduce weights on $\bar{N}$ of the type described in Section
\ref{ss1.1}. Let $V_{n}=z^{n}V_{0},$ and $W_{n}$ be the closure of
$z^{n+1}\bar{N}^{-}-z^{n}\bar{N}^{-}.$ Each $W_{n}$ is a fundamental domain
for $\mathbb{Z}.$ Let $\rho\left(  x\right)  $ be any $C^{\infty}$ real-valued
function on $\bar{N}$ with bounded gradient such that $\rho|V_{n}=n$ and
$\rho|W_{n}$ has values in $\left[  n,n+1\right]  $. Then the weight functions
are $k^{2\rho\left(  x\right)  }$. We index the ends of $M$ by subscripts. For
weights on $M$, extend the $\rho_{i}|\bar{N}^{+}$ to a function with values in
$\left[  -1,0\right]  $ outside the union of the $\bar{N}_{i}^{+}.$

An end is said to be quasi-periodic if the restriction of the metric
$\left\langle \cdot,\cdot\right\rangle $ on $M$ to $U$ is quasi-isometric to
the restriction of the lift of some (and thus any) metric on $N.$ Suppose now
that the ends of $M$ are quasi-periodic with disjoint neighborhoods $U_{i}.$
We extend the restrictions of the lifted metrics in any way to a metric
$\left\langle \cdot,\cdot\right\rangle ^{\prime}$ on $M.$ Then $\left\langle
\cdot,\cdot\right\rangle $ and $\left\langle \cdot,\cdot\right\rangle
^{\prime}$ are quasi-isometric. By \ref{ss1.2} the de Rham complexes
$\bar{\Omega}_{d,k}$ for the two metrics are boundedly isomorphic. We can
therefore replace $\left\langle \cdot,\cdot\right\rangle $ by $\left\langle
\cdot,\cdot\right\rangle ^{\prime}.$

We apply Theorem \ref{th2.17}. Choose any smooth triangulations of the $N_{i}$
with the images of $V_{0i}$ subcomplexes. These lift to BG triangulations of
the $\bar{N}_{i}.$ Extending their restrictions to the $\bar{N}_{i}$ in any
way gives a BG triangulation of $M.$ Let $\pi$ be the group of covering
transformations of a regular covering $\tilde{M}$ of $M.$ Let $\psi$ be the
canonical $C^{\ast}\left(  \pi\right)  $-bundle over $M$. Then $\bar{\Omega
}_{d,k}$ is boundedly and spatially boundedly equivalent to $C_{k}^{\ast
}\left(  M;\psi\right)  .$ (We have removed the triangulating complex from the
notation.) In light of the remarks on duality at the end of the last section,
the proof of Theorem \ref{th0.3} is reduced to showing that the inclusions
$C_{\ast}\left(  M;\psi\right)  \rightarrow C_{\ast}^{k}\left(  M;\psi\right)
$ and $C_{\ast}^{k}\left(  M;\psi\right)  \rightarrow C_{\ast}^{\ell f}\left(
M;\psi\right)  $ are equivalences for the stated values of $k.$ In this
section we will identify $\pi$ with a quotient of $\pi_{1}\left(  M\right)  $
by choosing a lift of the basepoint to $\tilde{M}.$

\subsection{\textbf{\label{ss3.2}}}

Let $\kappa_{i}=\pi_{1}\left(  \bar{N}_{i}\right)  .$ $V_{0i}$ may be chosen
so that the inclusions induce isomorphisms $\kappa_{i}\cong\pi_{1}\left(
V_{0i}\right)  \cong\pi_{1}\left(  \bar{N}_{i}^{+}\right)  \cong\pi_{1}\left(
\bar{N}_{i}^{-}\right)  $ \cite[Theorem 13.11]{hura}. Let $r_{i}:\kappa
_{i}\rightarrow\pi_{1}\left(  M\right)  \rightarrow\pi$ be induced by $\bar
{N}_{i}^{+}\rightarrow M.$ Composing $r_{i}$ with the inclusion $\pi
\rightarrow C^{\ast}\left(  \pi\right)  $ gives a homomorphism $\kappa
_{i}\rightarrow C^{\ast}\left(  \pi\right)  .$ $\kappa_{i}$ acts on $C^{\ast
}\left(  \pi\right)  $ via this map. Let $\tilde{N}_{i}$ be the universal
cover of $\bar{N}_{i}$ and $\phi_{i}=\tilde{N}_{i}\times_{r_{i}}C^{\ast
}\left(  \pi\right)  $. The restrictions of $\phi_{i}$ and $\psi$ to $\bar
{N}_{i}^{+}$ may be identified, since they have the same holonomy. Thus
$C_{\ast}\left(  \bar{N}_{i}^{+};\phi_{i}\right)  $ may be identified with the
subcomplex $C_{\ast}\left(  \bar{N}_{i}^{+};\psi\right)  \subset$ $C_{\ast
}\left(  M;\psi\right)  .$

Let $C$ be a complex of $A$-modules. It is $A$-finitely dominated if it is
equivalent to a complex of finitely generated $A$-modules. According to
\cite[Proposition 6.1]{hura}, this is equivalent to the following: there is a
complex $E$ of finitely generated free $A$-modules and maps $i:C\rightarrow E
$ and $j:E\rightarrow C$ such that $ji$ is homotopic to the identity. A
subcomplex of an $A$-module complex is cofinite if the quotient complex is
finitely generated.

\begin{lemma}
\label{le3.1}If $C_{\ast}\left(  M;\psi\right)  $ is $C^{\ast}\left(
\pi\right)  $-finitely dominated, each $C_{\ast}\left(  \bar{N}_{i};\phi
_{i}\right)  $ is $C^{\ast}\pi$-finitely dominated.
\end{lemma}

%

%TCIMACRO{\TeXButton{Proof}{\proof}}%
%BeginExpansion
\proof
%EndExpansion
Since $\bigoplus_{i}C_{\ast}\left(  \bar{N}_{i}^{+};\phi_{i}\right)  $ is a
cofinite subcomplex of $C_{\ast}\left(  M;\psi\right)  ,$ it is finitely
dominated \cite[Proposition 6.9(iii)]{hura}. This plus an additional condition
is sufficient for the finite domination of $\bigoplus_{i}C_{\ast}\left(
\bar{N}_{i};\phi_{i}\right)  :$ there is a cofinite subcomplex $Y\subset
\bar{N}^{+}=\bigcup_{i}\bar{N}_{i}^{+}$ such that the inclusion $C_{\ast}%
^{lf}\left(  Y;\psi\right)  \rightarrow C_{\ast}^{lf}\left(  \bar{N}^{+}%
;\psi\right)  $ is nullhomotopic \cite[Propositions 23.15-23.17]{hura}.
Henceforth we omit the coefficients. Since $C_{\ast}\left(  M\right)  $ is
finitely dominated, there exists a chain homotopy $H$ of the identity of
$C_{\ast}\left(  M\right)  $ to a chain map whose image is a finitely
generated subcomplex $F.$ There are cofinite subcomplexes $Y_{i}\subset\bar
{N}_{i}^{+}$ with union $Y$ which is a manifold with boundary such that
$F\subset C_{\ast}\left(  \overline{M-Y}\right)  $ and $\operatorname{Im}%
H|C_{\ast}\left(  \partial N^{+}\right)  \subset C_{\ast}\left(
\overline{M-Y}\right)  .$ This is possible since $F$ and $C_{\ast}\left(
\partial N^{+}\right)  $ are finitely generated. Then $H$ gives a nullhomotopy
homotopy of pairs of $C_{\ast}\left(  \bar{N}^{+},\partial\bar{N}^{+}\right)
\rightarrow C_{\ast}\left(  M,\overline{M-Y}\right)  $. By Alexander-Lefschetz
duality, $C_{c}^{\ast}\left(  \bar{N}^{+}\right)  \rightarrow C_{c}^{\ast
}\left(  Y\right)  $ is nullhomotopic. Transposing, $C_{\ast}^{lf}\left(
Y\right)  \rightarrow C_{\ast}^{lf}\left(  \bar{N}^{+}\right)  $ is
nullhomotopic. Therefore $\bigoplus_{i}C_{\ast}\left(  \bar{N}_{i};\phi
_{i}\right)  $ is finitely dominated

If a sum of complexes is finitely dominated, then each summand is. For let
$\bigoplus_{i}C_{i}\rightarrow E\rightarrow$ $\bigoplus_{i}C_{i}$ be a
domination. Restriction and projection induce dominations $C_{i}\rightarrow
E\rightarrow C_{i}.$%
%TCIMACRO{\TeXButton{End Proof}{\endproof}}%
%BeginExpansion
\endproof
%EndExpansion

\noindent The converse of this Lemma is also true by \cite[23.17, 6.2ii]{hura}.

Let $K$ be a subcomplex of $M$. Consider the algebraic mapping cones of the
inclusions
\begin{align*}
\hat{C}_{\ast}^{k}\left(  K\right)   &  =\mathcal{C}\left(  C_{\ast}\left(
K\right)  \rightarrow C_{\ast}^{k}\left(  K\right)  \right)  ,\\
\check{C}_{\ast}^{k}\left(  K\right)   &  =\mathcal{C}\left(  C_{\ast}%
^{k}\left(  K\right)  \rightarrow C_{\ast}^{lf}\left(  K\right)  \right)  .
\end{align*}
We will show that if $C_{\ast}\left(  M\right)  $ is finitely dominated,
$\hat{C}_{\ast}^{k}\left(  M\right)  $ is contractible for $k$ sufficiently
large, and $\check{C}_{\ast}^{k}\left(  M\right)  $ is contractible for $k>0$
sufficiently close to 0. This will give the claimed equivalences.

\begin{lemma}
\label{le3.2}Let $L\subset K$ be a cofinite subcomplex. Then the inclusion
induces equivalences on $\hat{C}_{\ast}^{k}$ and $\check{C}_{\ast}^{k}$ for
all $k.$
\end{lemma}

%

%TCIMACRO{\TeXButton{Proof}{\proof}}%
%BeginExpansion
\proof
%EndExpansion
This is a small adaptation of an argument in \cite[Prop. 3.13]{hura}. We
sketch the first, the second being similar. The map
\[
q:C_{\ast}\left(  K\right)  /C_{\ast}\left(  L\right)  \rightarrow C_{\ast
}^{k}\left(  K\right)  /C_{\ast}^{k}\left(  L\right)  .
\]
is an isomorphism. For let $c\in C_{\ast}^{k}\left(  K\right)  ,$ and
$\tilde{c}$ be gotten by setting $c$ to zero outside of $\overline{K-L}.$ Then
$\tilde{c}\in C_{\ast}\left(  K\right)  $ and $c-\tilde{c}\in C_{\ast}%
^{k}\left(  L\right)  ,$ so $q$ is surjective. Let $e\in C_{\ast}\left(
K\right)  \cap C_{\ast}^{k}\left(  L\right)  .$ Then there are $e_{i}\in
C_{\ast}\left(  L\right)  $ which converge to $e$ in the $k$-norm. Since each
$e_{i}$ has support in $L,$ so does $e,$ so $e\in C_{\ast}\left(  L\right)  $
and $q$ is injective. There is an exact sequence
\[
0\rightarrow\hat{C}_{\ast}^{k}\left(  L\right)  \rightarrow\hat{C}_{\ast}%
^{k}\left(  K\right)  \rightarrow\mathcal{C}\left(  q\right)  \rightarrow0.
\]

\noindent$C\left(  q\right)  $ is a free $A$-module and contractible since $q$
is an isomorphism. Therefore the first map is an equivalence.%
%TCIMACRO{\TeXButton{End Proof}{\endproof}}%
%BeginExpansion
\endproof
%EndExpansion

\subsection{\textbf{\label{ss3.3}}}

We apply this to replace $M$ by the union of the $\bar{N}_{i}^{+}.$ This
reduces the problem to working on the $\bar{N}_{i}.$ From this point on the
ends may be treated separately. The subscripts will therefore be omitted. We
put things into a more algebraic context. For a unital $C^{\ast}$-algebra $A$,
we consider the category of extended $A\left[  z,z^{-1}\right]  $-modules.
$A\left[  z,z^{-1}\right]  $ is the ring of Laurent polynomials. Such a module
$P $ is of the form $P^{0}\otimes_{A}A\left[  z,z^{-1}\right]  =P^{0}\left[
z,z^{-1}\right]  $ for some finitely generated Hilbert $A$-module $P^{0}.$
Thus we can write $P=\bigoplus_{n}P^{n}=\bigoplus_{n}P^{0}z^{n}.$ Finitely
generated free $A\left[  z,z^{-1}\right]  $-modules are included, since
$\left(  A\left[  z,z^{-1}\right]  \right)  ^{N}\cong A^{N}\left[
z,z^{-1}\right]  .$ If $\left\langle \cdot,\cdot\right\rangle $ is the inner
product on $P^{0},$ one is defined on $P$ by $\left(  \sum_{n}c_{n}z^{n}%
,\sum_{n}d_{n}z^{n}\right)  =\sum_{n}\left\langle c_{n},d_{n}\right\rangle .$
More generally, there are $k$-inner products $\left(  \cdot,\cdot\right)
_{k}$ where the right hand side is replaced by $\sum_{n}\left\langle
c_{n},d_{n}\right\rangle k^{2n}.$ Note that the $P^{n}$ are orthogonal for any
$k.$

We denote the completions of $P$ by $P_{\left(  k\right)  }.$ $P_{\left(
k\right)  }$ is the Hilbert module exterior tensor product $P^{0}%
\otimes\mathbb{C}\left[  z,z^{-1}\right]  _{\left(  k\right)  }.$ The set
$\left\{  e_{n}\right\}  =\left\{  k^{-n}z^{n}\right\}  $ is an orthonormal
basis for $\mathbb{C}\left[  z,z^{-1}\right]  _{\left(  k\right)  }$. Any
element $c$ of $P_{\left(  k\right)  }$ may therefore be written as $\sum
_{n}a_{n}e_{n}$ with $a_{n}\in P^{0}$ and $\sum_{n}\left\langle a_{n}%
,a_{n}\right\rangle $ norm convergent, or as $\sum_{n}c_{n}z^{n}$ with
$c_{n}=k^{-n}a_{n}.$ Since $e_{n}z=ke_{n+1},$ multiplication by $z$ has
operator norm $k.$ From this it follows that $cz=\sum c_{n}z^{n+1}.$ We will
also use $A\left[  z\right]  $- and $A\left[  z^{-1}\right]  $- extended
modules. There are similar discussions for them.

An homomorphism $T:$ $P\rightarrow Q$ of extended $A\left[  z,z^{-1}\right]
$-modules may be described by a finite sum $\sum_{n}z^{n}T_{n}$, where each
$T_{n}:P^{0}\rightarrow P^{0}.$ The analog of spatial boundedness is
finiteness of the sum. An $A\left[  z\right]  $- $\left(  A\left[
z^{-1}\right]  \text{-}\right)  $ module homomorphism may be described by a
similar sum with $n\ge0$ $\left(  n\le0\right)  .$

\begin{lemma}
\label{le3.3}$T$ is bounded in any $k$-norm.
\end{lemma}

%

%TCIMACRO{\TeXButton{Proof}{\proof}}%
%BeginExpansion
\proof
%EndExpansion
This is the matrix case of Proposition \ref{pr2.8}. The matrix entries are
$T_{nm}=T_{n-m}z^{n-m}:P^{m}\rightarrow P^{n}.$ Since there are finitely many
$T_{n}$ and $P^{0}$ is finitely generated the $T_{n-m}$ are uniformly bounded.
$\left\Vert cz\right\Vert _{k}=k\left\Vert c\right\Vert _{k},$ so the $T_{nm}$
are uniformly bounded in the $k$-norm.%
%TCIMACRO{\TeXButton{End Proof}{\endproof}}%
%BeginExpansion
\endproof
%EndExpansion

By continuity, the extension of $T$ to $P_{\left(  k\right)  }$ is an
$A\left[  z,z^{-1}\right]  $- $\left(  A\left[  z\right]  \text{-},A\left[
z^{-1}\right]  \text{-}\right)  $ module homomorphism. Since $cz=\sum
c_{n}z^{n+1},$ it is again given by $\sum_{n}z^{n}T_{n}.$

In general, $z$ induces an automorphism $\alpha$ of $\kappa=\pi_{1}\left(
\bar{N}\right)  ,$ which is well-defined up to inner automorphism. We will
assume the following.

\begin{enumerate}
\item[ ] $\mathcal{G}$: For each $i,$ $\pi_{1}\left(  N_{i}\right)  =\pi
_{1}\left(  \bar{N}_{i}\right)  \times\mathbb{Z}=\kappa_{i}\times\mathbb{Z}.$
\end{enumerate}

$\phi$ was defined as $\tilde{N}\times_{r}C^{\ast}\left(  \pi\right)  .$ We
define a flat bundle $\phi^{\prime}$ over $N.$ Let $r^{\prime}=rp_{1}%
:\kappa_{i}\times\mathbb{Z\rightarrow}C^{\ast}\left(  \pi\right)  .$ Then
$\phi^{\prime}=\tilde{N}\times_{r^{\prime}}C^{\ast}\left(  \pi\right)  .$
$\left\{  e\right\}  \times\mathbb{Z\subset\pi}_{1}\left(  N\right)  $ acts
freely on $\phi$ preserving fibers, with quotient $\phi^{\prime}.$ Let
$\bar{N}$ be triangulated as described in Section \ref{ss3.1}. It follows that
$C=C_{\ast}\left(  \bar{N};\phi\right)  $ is a complex of finitely generated
free $C^{\ast}\left(  \pi\right)  \left[  z,z^{-1}\right]  $-modules.

To fix a generating module, let $C^{0}$ be the $A$-module generated by
simplexes in $W_{0}-V_{1}.$ Then $C=C^{0}\left[  z,z^{-1}\right]  .$ By
construction, $C^{+}=C^{0}\left[  z\right]  $ is a subcomplex, corresponding
to $\bar{N}^{+}.$ Two slightly different $k$-inner products have been
described for $C:$ one using $\rho\left(  x\right)  ,$ the other in this
subsection. If $\sigma$ is a simplex in $W_{n}-V_{n+1},$ and $x\in\sigma,$
then $n\le\rho\left(  x\right)  \le n+1.$ It follows that the two $k$-norms
are equivalent.

We discuss a general notion of locally finite chains. Let $P$ be a module over
any ring with a decomposition $P=\bigoplus_{i}P^{i}.$ The locally finite
module is $P^{\ell f}=\prod_{i}P^{i}.$ Given a complex $D$ with a
decomposition of each $D_{j},$ $D^{\ell f}$ is also a complex with the
extended differentials. For simplicial chains, decomposed by the simplexes,
this gives the locally finite chains. We can therefore identify $C_{\ast
}^{\ell f}\left(  \bar{N};\phi\right)  $ and $C^{\ell f}$ in the present
sense. For a complex $C$ of extended $A{}\!{}\left[  z\right]  $-modules, we
use the decompositions $C_{j}=\bigoplus C_{j}^{n}.$ In the simplicial case
this is the same as that given by the simplexes, since the $C_{j}^{n}$ are
finitely generated. We can identify $C^{\ell f}$ with $C\otimes_{A\left[
z\right]  }A\left[  \left[  z\right]  \right]  $ (the formal power series
ring). An $A\left[  z\right]  $-module chain map $T:C\rightarrow D$ induces
one $C^{\ell f}\rightarrow D^{\ell f}$ using the expression $T=\sum_{n\geq
0}z^{n}T_{n}.$ It follows that the action of $T$ on $C^{\ell f}$ is an
extension of its action on any $C_{\left(  k\right)  }.$

Let $C$ be an $A\left[  z,z^{-1}\right]  $-complex and $C^{+}=C^{0}\left[
z\right]  .$ We assume that $C^{+}$ is a subcomplex. Then $C^{+}z^{n}$ is a
subcomplex of $C,$ and $C^{+,\ell f}z^{n}$ is a subcomplex of $C^{\ell f}.$ In
analogy with Lemma \ref{le3.2}, for any $k>0$ we define
\begin{align*}
\hat{C}_{\left(  k\right)  }^{+}  &  =\mathcal{C}\left(  C^{+}\rightarrow
C_{\left(  k\right)  }^{+}\right)  ,\\
\check{C}_{\left(  k\right)  }^{+}  &  =\mathcal{C}\left(  C_{\left(
k\right)  }^{+}\rightarrow C^{+,lf}\right)  .
\end{align*}
More generally for $n\in\mathbb{Z}$ there are
\begin{align}
\hat{C}_{\left(  k\right)  }^{+,n}  &  =\mathcal{C}\left(  C^{+}%
z^{n}\rightarrow\left(  C^{+}z^{n}\right)  _{\left(  k\right)  }\right)
,\nonumber\\
\check{C}_{\left(  k\right)  }^{+,n}  &  =\mathcal{C}\left(  \left(
C^{+}z^{n}\right)  _{\left(  k\right)  }\rightarrow C^{+,\ell f}z^{n}\right)
.
\end{align}
We sometimes omit the $k$ for simplicity. These constructions are natural. For
example, consider an $A\left[  z\right]  $-module chain map or homotopy
$s:C^{+}z^{n}\rightarrow D^{+}z^{m}.$ Since the extensions to $\left(
C^{+}z^{n}\right)  _{\left(  k\right)  }$ and $C^{+,\ell f}z^{n}$ are
compatible, there is an induced $A\left[  z\right]  $-module map or homotopy
$\check{s}:\check{C}_{\left(  k\right)  }^{+,n}\rightarrow\check{D}_{\left(
k\right)  }^{+,m}$. All these definitions may be repeated under the assumption
that $C^{-}=C^{0}\left[  z^{-1}\right]  $ is a subcomplex.

\begin{lemma}
\label{le3.4}An equivalence $C\rightarrow D$ of $A\left[  z,z^{-1}\right]
$-module complexes such that $C^{+}$ and $D^{+}$ are subcomplexes induces
$A$-module equivalences $\hat{C}_{\left(  k\right)  }^{+}\rightarrow\hat
{D}_{\left(  k\right)  }^{+}$ and $\check{C}_{\left(  k\right)  }%
^{+}\rightarrow\check{D}_{\left(  k\right)  }^{+}$ for any $k>0.$ There is a
similar statement for $C^{-}$ and $D^{-}.$
\end{lemma}

%

%TCIMACRO{\TeXButton{Proof}{\proof}}%
%BeginExpansion
\proof
%EndExpansion
We take the first case, the others differing only in notation. The proof
consists of constructing a functor $\mathcal{F}$ from the homotopy category of
$A\left[  z,z^{-1}\right]  $-module chain maps $C\rightarrow D$ to that of
$A$-module chain maps $\hat{C}^{+}\rightarrow\hat{D}^{+}.$ With a proof like
that of Lemma \ref{le3.2}, inclusions induce $A$-module equivalences
$h_{n}:\hat{C}^{+}\rightarrow\hat{C}^{+,-n}$ for $n>0.$ Let $r_{n}$ be
homotopy inverses.

Suppose given a map $f:C\rightarrow D.$ Since $C^{0}$ is finitely generated,
for any $m$ $f\left(  C^{+}z^{-m}\right)  \subset D^{+}z^{-n}$ for all
sufficiently large $n.$ Denote the induced map $\hat{C}^{+,-m}\rightarrow$
$\hat{D}^{+,-n}$ by $\hat{f}_{mn}.$ $\mathcal{F}\left(  f\right)  $ is
represented by $r_{n}\hat{f}_{0n}:\hat{C}^{+}\rightarrow\hat{D}^{+}$ for any
$n$ such that $\hat{f}_{0n}$ is defined. We show that different choices of $n$
give homotopic maps. Suppose that $m>n$ and let $j:\hat{D}^{+,-n}%
\rightarrow\hat{D}^{+,-m}$ be the inclusion. $\left(  r_{m}j\right)
h_{n}=r_{m}\left(  jh_{n}\right)  =r_{m}h_{m}\sim I.$ Since $h_{m}$ is an
equivalence, $r_{m}j$ is a homotopy inverse of $h_{n},$ so is homotopic to
$r_{n}.$ Then $r_{m}\hat{f}_{0m}=r_{m}j\hat{f}_{0n}\sim r_{n}\hat{f}_{0n}.$

If $H:C\rightarrow D$ is a homotopy between $f$ and $g$, $\mathcal{F}\left(
H\right)  $ is represented by $r_{n}\hat{H}_{0n}$ for any $n$ such that
$\hat{f}_{0n},$ $\hat{g}_{0n},$ and $\hat{H}_{0n}$ are defined. $\partial
r_{n}\hat{H}_{0n}+r_{n}\hat{H}_{0n}\partial=r_{n}\hat{f}_{0n}-r_{n}\hat
{g}_{0n}. $

Given $f:C\rightarrow D$ and $g:D\rightarrow E,$ choose $n$ so that $\hat
{f}_{0n}$ is defined, then $m$ so that $\hat{g}_{nm}$ is defined. Then
$\mathcal{F}\left(  gf\right)  $ is represented by $r_{m}\widehat{\left(
gf\right)  }_{0m}$, and $\mathcal{F}\left(  g\right)  \mathcal{F}\left(
f\right)  $ by $r_{m}\hat{g}_{0m}r_{n}\hat{f}_{0n}.$ $\hat{g}_{nm}h_{n}%
=\hat{g}_{0m},$ so $\hat{g}_{nm}\sim\hat{g}_{0m}r_{n}.$ Therefore $r_{m}%
\hat{g}_{0m}r_{n}\hat{f}_{0n}\sim r_{m}\hat{g}_{nm}\hat{f}_{0n}=r_{m}\left(
\widehat{gf}\right)  _{0m}.$ Therefore $\mathcal{F}$ preserves composition up
to homotopy.%
%TCIMACRO{\TeXButton{End Proof}{\endproof}}%
%BeginExpansion
\endproof
%EndExpansion

\subsection{\textbf{\label{ss3.4}}}

Let $X$ be a space and $h$ a self-map. The mapping torus $\mathcal{T}\left(
h\right)  $ is the quotient $X\times I/\left\{  \left(  x,1\right)  =\left(
h\left(  x\right)  ,0\right)  \right\}  .$ It has an infinite cyclic cover
\[
\mathcal{\bar{T}}\left(  h\right)  =\bigcup_{j=-\infty}^{\infty}X\times
I\times\left\{  j\right\}  /\left\{  \left(  x,1,j\right)  =\left(  h\left(
x\right)  ,0,j+1\right)  \right\}  ,
\]
the doubly infinite mapping telescope. $\mathbb{Z}$ acts on$\mathcal{\ \bar
{T}}\left(  h\right)  $ by $\left(  n,\left(  x,t,j\right)  \right)
\rightarrow\left(  x,t,j+n\right)  .$ Suppose that $X$ is a CW complex and $h$
is a cellular map. Ranicki observed that the cellular chain complex of
$\mathcal{\bar{T}}\left(  h\right)  $ is the algebraic mapping torus
\[
T_{\ast}\left(  h_{\ast}\right)  =\mathcal{C}\left(  I-zh_{\ast}:C_{\ast
}\left(  X\right)  \left[  z,z^{-1}\right]  \rightarrow C_{\ast}\left(
X\right)  \left[  z,z^{-1}\right]  \right)  .
\]

Now let $C$ be a complex of extended $A\left[  z,z^{-1}\right]  $-modules. Let
the $A$-module homomorphism of $C$ given by $z$ be $\zeta.$ By \cite[p.263]%
{hura} there is an $A\left[  z,z^{-1}\right]  $-module chain equivalence
$s:C\rightarrow T\left(  \zeta^{-1}\right)  .$ If $C$ is finitely dominated,
$C$ is equivalent to a complex of finitely generated $A$-modules $P.$ There is
then an induced $A\left[  z,z^{-1}\right]  $-module equivalence
$t:C\rightarrow T\left(  h\right)  ,$ where $h$ is a self-equivalence of $P$
induced from $\zeta^{-1}.$ We equip $P$ with any $A$-valued inner product, and
$T\left(  h\right)  $ with a $k$-inner product as described in \ref{ss3.3}$.$
From now on we will write $T$ for $T\left(  h\right)  .$ The composition
$ts:C\rightarrow T$ is an $A\left[  z,z^{-1}\right]  $-module chain
equivalence. By Lemma \ref{le3.3} it extends to an equivalence of the
completions $C_{\left(  k\right)  }$ and $T_{\left(  k\right)  }.$ According
to Theorem \ref{th2.17}, Lemma \ref{le3.2}, and Lemma \ref{le3.4}, the
equivalence of $\bar{\Omega}_{d,k}$ and $C^{\ast}\left(  M;\psi\right)  $ for
all $k>0$ which are sufficiently small will follow if we show that $\check
{T}_{\left(  k\right)  }^{+}$ is contractible. By Lemma \ref{le3.4}, this
doesn't depend on the choice of $T^{+}.$

For the equivalence of $\bar{\Omega}_{d,k}$ and $C_{c}^{\ast}\left(
M;\psi\right)  $ for $k$ large, it is notationally convenient to use the
reversed complex of $C.$ There are two choices for the generator $z$ of the
action of $\mathbb{Z}$ on $C.$ The reversed complex $_{r}C$ is $C$ with the
actions of $z$ and $z^{-1}$ interchanged. This change has no topological
significance. The $\pm$ labels of the ends are switched. Replace $C$ by
$_{r}C.$ According to our notational conventions, $\left(  _{r}C\right)
_{\left(  k\right)  }=_{r}\!\left(  C_{\left(  1/k\right)  }\right)  .$ We
then wish to show that $\widehat{\left(  _{r}C\right)  }_{\left(  k\right)
}^{-}$ is contractible for all small $k>0$ . By Lemma \ref{le3.4}, it is
sufficient do the same for $\hat{T}_{\left(  k\right)  }^{-}$.

Let $\mathcal{\bar{T}}^{\infty}\left(  h\right)  $ be $\mathcal{\bar{T}%
}\left(  h\right)  $ with the positive end compactified by a point $\infty.$
There is an evident homotopy contracting $\mathcal{\bar{T}}^{\infty}\left(
h\right)  $ to $\infty.$ We consider the corresponding homotopy of $T.$ The
first part of the following proof is the analytic counterpart of Ranicki's
result on the vanishing of homology with Novikov ring coefficients. (The
reader may wish to consider the simplest example first: $P=\mathbb{C}$ in
degree 0, $h=I.$ This gives the standard chain complex of $\mathbb{R}.)$

\begin{proposition}
\label{pr3.5}$\check{T}_{\left(  k\right)  }^{+}$ and $\hat{T}_{\left(
k\right)  }^{-}$ are contractible for all $k>0$ which are sufficiently small.
\end{proposition}

%

%TCIMACRO{\TeXButton{Proof}{\proof}}%
%BeginExpansion
\proof
%EndExpansion
$T$ is described by%

\begin{align*}
T_{j}  &  =P_{j}\left[  z,z^{-1}\right]  \oplus P_{j-1}\left[  z,z^{-1}%
\right]  ,\\
\partial_{j}  &  =\left(
\begin{array}
[c]{cc}%
\partial_{j} & \left(  -\right)  ^{j}\left(  I-zh\right) \\
0 & \partial_{j-1}%
\end{array}
\right)  :\text{ }P_{j}\left[  z,z^{-1}\right]  \oplus P_{j-1}\left[
z,z^{-1}\right]  \rightarrow\\
&  \rightarrow P_{j-1}\left[  z,z^{-1}\right]  \oplus P_{j-2}\left[
z,z^{-1}\right]  .
\end{align*}

\noindent It is generated by $T^{0}=P\oplus P_{\ast-1.}$ Since the norm of
multiplication by $z$ is $k,$ $\left\Vert zh\right\Vert _{k}\leq k\left\Vert
h\right\Vert .$ Thus for all $k<\left\Vert h\right\Vert ^{-1},$ $I-zh$ is
invertible in the $k$-norm with inverse $r=\sum_{n=0}^{\infty}\left(
zh\right)  ^{n}.$ Then
\[
H_{j}=\left(
\begin{array}
[c]{cc}%
0 & 0\\
\left(  -\right)  ^{j}r & 0
\end{array}
\right)
\]

\noindent is a bounded $A\left[  z,z^{-1}\right]  $-module contraction of
$T_{\left(  k\right)  }$. There are now two cases.

Let $T^{+}=\mathcal{C}\left(  I-zh:P\left[  z\right]  \rightarrow P\left[
z\right]  \right)  .$ This is generated by $T^{+,0}=P\oplus P_{*-1.}$ Since
$r$ preserves $T_{\left(  k\right)  }^{+},$ $H$ restricts to a contraction
$H^{+}$ of $T_{\left(  k\right)  }^{+}.$ Since any $T^{+,n}$ is in the image
by $H^{+}$ of only finitely many others, $H^{+}$ extends to an $A\left[
z\right]  $-module contraction of $T^{+,lf}.$ Thus $\check{T}_{\left(
k\right)  }^{+}$ is contractible.

Let $T^{-}=\mathcal{C}\left(  I-zh:P\left[  z^{-1}\right]  z^{-1}\rightarrow
P\left[  z^{-1}\right]  \right)  .$ This is generated by $P\oplus
P_{*-1}z^{-1}. $ However, there seems to be no advantage in using the
associated decomposition, and we will continue to use the one above.
\[
T^{-,n}=\left\{
\begin{array}
[c]{ll}%
T^{n}, & n<0\\
P\oplus0, & n=0\\
0, & n>0.
\end{array}
\right.
\]
$T^{+}\cap T^{-}=P\oplus0$ will be identified with $P.$ Let $i^{-}$ and
$q^{-}$ be the injection of and projection onto $T^{-}.$ The latter isn't a
chain map. If $H^{-}=q^{-}Hj^{-}$ is expanded in a series using the series for
$r,$ only finitely many terms are nonzero on any element of $T^{-}.$ Therefore
$H^{-}$ induces an $A$-module homomorphism $T^{-}\rightarrow T^{-}. $%
\[
\partial H^{-}+H^{-}\partial=\partial q^{-}Hi^{-}+q^{-}\left(  I-\partial
H\right)  i^{-}=I_{T^{-}}+\left(  \partial q^{-}-q^{-}\partial\right)
Hi^{-}.
\]

We compute $\partial q^{-}-q^{-}\partial.$

\begin{enumerate}
\item[ ] On $T^{-},$ since it is a subcomplex, $\partial q^{-}-q^{-}%
\partial=\partial-\partial=0.$

\item[ ] On $T^{n}$ for $\dot{n}>0,$ $q^{-}=0$ and $q^{-}\partial=0,$ so
$\partial q^{-}-q^{-}\partial=0.$

\item[ ] On $0\oplus P_{*-1}\subset T^{0},$ $\partial q^{-}-q^{-}%
\partial=-q^{-}\partial=\left(  -\right)  ^{j+1}:0\oplus P_{j-1}\rightarrow
P_{j-1}\oplus0.$
\end{enumerate}

Thus $\partial H^{-}+H^{-}\partial=I$ $-\ell,$ where
\[
\ell=\left(  I+0\right)  +\left(  \Sigma_{n=1}^{\infty}\left(  \left(
zh\right)  ^{n}+0\right)  \right)  :\left(  P\oplus0\right)  \oplus\left(
\bigoplus_{n=1}^{\infty}\left(  P\oplus P_{\ast-1}\right)  z^{-n}\right)
\rightarrow P.
\]
The same relation holds on $T_{\left(  k\right)  }^{-}$ with all operators
replaced by their bounded extensions. By the definitions of the inner
products, $P^{k}=P.$ The extension of $H^{-}$ is therefore a homotopy from the
identity of $T_{\left(  k\right)  }^{-}$ to a map to $P\subset T^{-},$ which
takes $T^{-}$ to itself. Therefore the inclusion of $T^{-}$ in $T_{\left(
k\right)  }^{-}$ is an equivalence, and $\hat{T}_{\left(  k\right)  }^{-}$ is
contractible.%
%TCIMACRO{\TeXButton{End Proof}{\endproof}}%
%BeginExpansion
\endproof
%EndExpansion

\section{\label{s4}Examples}

In this section we give examples for Theorems \ref{th0.1} and \ref{th0.2}.

\subsection{\textbf{\label{ss4.1}}}

The Euler characteristic takes all integer values in all dimensions $\geq4,$
even for manifolds with cylindrical ends. There exists a closed surface with
any given value of $\chi.$ It may be embedded in $\mathbb{R}^{n}$ for any
$n\geq4.$ The normal disk bundle is an orientable manifold with boundary with
the same $\chi.$ Then attach a cylinder over the boundary.

We showed that the complex chains on an end satisfying the hypotheses of
Theorem \ref{th0.1} are equivalent near infinity to the algebraic mapping
torus of a homotopy equivalence. This means that rationally, the end looks
like a cylinder. However, if torsion is taken into account, this need not be
the case. Let $\bar{N}$ be the connected sum of $S^{n-1}\times\lbrack
0,\infty)$ with countably many copies of $\mathbb{RP}^{n},$ attached
periodically. Attach $D^{n}$ to $\bar{N}$ along $S^{n-1}\times\left\{
0\right\}  $ to obtain $M.$ Then $M$ is rationally acyclic and is orientable
for $n$ odd, but has infinitely generated 2-torsion.

\subsection{\textbf{\label{ss4.2} }}

We will first relate the $\tilde{K}_{0}\left(  C^{\ast}\left(  \pi\right)
\right)  $-valued Euler characteristic $\tilde{\chi}_{C^{\ast}\left(
\pi\right)  }$ to Wall's finiteness obstruction \cite{wal1}, \cite{wal2}. We
will then give examples of manifolds satisfying the hypotheses of Theorem
\ref{th0.2} (for the universal cover) for which $\tilde{\chi}_{C^{\ast}\left(
\pi\right)  }\neq0.$ It follows that the index of $\bar{D}_{k}^{even}$ in
$\tilde{K}_{0}\left(  C^{\ast}\left(  \pi\right)  \right)  $ is nonzero for
$k>0$ large or small. In the basic examples, $\pi$ is a finite group, and the
invariant is an equivariant Euler characteristic taking values in the reduced
representation ring $\tilde{R}\left(  \pi\right)  =\tilde{K}_{0}\left(
\mathbb{C}\left[  \pi\right]  \right)  =\tilde{K}_{0}\left(  C^{\ast}\left(
\pi\right)  \right)  $ . Examples with infinite groups are constructed using
free products and semidirect products. Examples with torsion-free $\pi$ are
not known and are unlikely.

C. T. C. Wall introduced an obstruction to finiteness up to homotopy for
certain CW complexes $X.$ Let $C_{\ast}\left(  \tilde{X}\right)  $ be the
cellular chain complex of the universal cover of $X.$ Let $\pi=\pi\left(
X\right)  $ be the group of covering transformations of $\tilde{X}.$ Suppose
that $C_{\ast}\left(  \tilde{X}\right)  $ is $\mathbb{Z}\left[  \pi\right]
$-finitely dominated, i.e. chain homotopy equivalent to a finite-dimensional
complex of finitely generated projective $\mathbb{Z}\left[  \pi\right]
$-modules $F$. Define $o_{X}=\Sigma\left(  -\right)  ^{i}\left[  F_{i}\right]
\in\tilde{K}_{0}\left(  \mathbb{Z}\left[  \mathbb{\pi}\right]  \right)  .$
This is independent of the choice of $F.$ If $\pi$ is finitely presented, $X $
is homotopy equivalent to a finite CW complex if and only if $o_{X}=0.$ Wall
\cite{wal2} considered the effect of a change of rings. Let $R$ be any ring,
and $v:$ $\mathbb{Z}\left[  \pi\right]  \rightarrow R$ a homomorphism,
inducing $v_{\ast}:\tilde{K}_{0}\left(  \mathbb{Z}\left[  \mathbb{\pi}\right]
\right)  \rightarrow\tilde{K}_{0}\left(  R\right)  .$ $\tilde{\chi}%
_{R}=v_{\ast}\left(  o_{X}\right)  $ is the Euler characteristic of $C_{\ast
}\left(  \tilde{X}\right)  \otimes_{v}R.$ The point is that $\tilde{\chi}_{R}$
may be defined in cases where $o_{X}$ is not. We will consider the inclusion
$v:\mathbb{Z\pi}\rightarrow C^{\ast}\left(  \pi\right)  .$ $C_{\ast}\left(
\tilde{X}\right)  \otimes_{v}C^{\ast}\left(  \pi\right)  $ may be identified
with the local coefficient chains of $X$ with coefficients in the bundle
$\psi=\tilde{X}\times_{v}C^{\ast}\left(  \pi\right)  .$ (See Lemma
\ref{le4.1}.) Unfortunately, there seem to be no known cases where $v_{\ast
}:\tilde{K}_{0}\left(  \mathbb{Z}\left[  \mathbb{\pi}\right]  \right)
\rightarrow\tilde{K}_{0}\left(  C^{\ast}\left(  \pi\right)  \right)  $ is
nonzero. However, we give examples where $\tilde{\chi}_{C^{\ast}\left(
\pi\right)  }$ is nonzero. The basic ingredients are idempotents in
$\mathbb{Q}\left[  \pi\right]  $ which represent nonzero elements of
$\tilde{K}_{0}\left(  C^{\ast}\left(  \pi\right)  \right)  .$

Let $\pi$ be a finite group. Then $C^{\ast}\left(  \pi\right)  =\mathbb{C}%
\left[  \mathbb{\pi}\right]  .$ Let $p:\mathbb{C}\left[  \mathbb{\pi}\right]
\mathbb{\rightarrow C}\left[  \mathbb{\pi}\right]  $ be the idempotent given
by multiplication by a central idempotent. If $p$ is not $0$ or the identity,
its image $P$ represents a nonzero element of $\tilde{K}_{0}\left(  C^{\ast
}\left(  \mathbb{\pi}\right)  \right)  .$ Suppose that the idempotent has
rational coefficients. For example, this is always the case if $\pi$ is a
symmetric group \cite[Section II.3]{sim}. Then $\tilde{K}_{0}\left(
\mathbb{Q}\left[  \pi\right]  \right)  \rightarrow\tilde{K}_{0}\left(
C^{\ast}\left(  \pi\right)  \right)  $ is an isomorphism. The simplest example
is $\pi=\mathbb{Z}_{2}=\left\{  e,g\right\}  $ with the idempotent $\frac
{1}{2}\left(  e+g\right)  $ corresponding to the trivial 1-dimensional representation.

Let $\pi$ and $\rho$ be any groups, and $\pi*\rho$ their free product. By
\cite[Theorem 5.4]{lan1}, the bottom row of
\[%
\begin{array}
[c]{ccc}%
\tilde{K}_{0}\left(  \mathbb{Q}\left[  \pi\right]  \right)  \oplus\tilde
{K}_{0}\left(  \mathbb{Q}\left[  \rho\right]  \right)  & \longrightarrow &
\tilde{K}_{0}\left(  \mathbb{Q}\left[  \pi*\rho\right]  \right) \\
\downarrow &  & \downarrow\\
\tilde{K}_{0}\left(  C^{*}\left(  \pi\right)  \right)  \oplus\tilde{K}%
_{0}\left(  C^{*}\left(  \rho\right)  \right)  & \longrightarrow & \tilde
{K}_{0}\left(  C^{*}\left(  \pi*\rho\right)  \right)
\end{array}
\]

\noindent is an isomorphism. Therefore if either $\tilde{K}_{0}\left(
\mathbb{Q}\left[  \pi\right]  \right)  \rightarrow\tilde{K}_{0}\left(
C^{\ast}\left(  \pi\right)  \right)  $ or $\tilde{K}_{0}\left(  \mathbb{Q}%
\left[  \rho\right]  \right)  \allowbreak\rightarrow$ $\tilde{K}_{0}\left(
C^{\ast}\left(  \rho\right)  \right)  $ is nonzero, the map on the right is as well.

Let $\pi$ be any group, and $\alpha:\mathbb{Z\rightarrow}Aut\left(
\pi\right)  $ a homomorphism. Let $\pi\rtimes_{\alpha}\mathbb{Z}$ be the
semidirect product. $C^{\ast}\left(  \pi\rtimes_{\alpha}\mathbb{Z}\right)
=C^{\ast}\left(  \pi\right)  \rtimes_{\alpha}\mathbb{Z},$ the crossed product
algebra. Suppose that the composition $\tilde{K}_{0}\left(  \mathbb{Q}\left[
\pi\right]  \right)  \rightarrow\tilde{K}_{0}\left(  C^{\ast}\left(
\pi\right)  \right)  \rightarrow\tilde{K}_{0}\left(  C^{\ast}\left(
\pi\rtimes_{\alpha}\mathbb{Z}\right)  \right)  $ is nonzero. Then by a
naturality argument like the preceding, $\tilde{K}_{0}\left(  \mathbb{Q}%
\left[  \pi\rtimes_{\alpha}\mathbb{Z}\right]  \right)  \rightarrow\tilde
{K}_{0}\left(  C^{\ast}\left(  \pi\rtimes_{\alpha}\mathbb{Z}\right)  \right)
$ is nonzero. For example, let $\pi=\mathbb{Z}_{2}\times\mathbb{Z}_{2}$ with
generators $g_{0}$ and $g_{1}$ and $\alpha\left(  1\right)  \left(
g_{i}\right)  =g_{1-i}.$ By the Pimsner-Voiculescu sequence \cite{pivo},
$\tilde{K}_{0}\left(  \mathbb{C}^{\ast}\left(  \pi\right)  \right)
\rightarrow\tilde{K}_{0}\left(  C^{\ast}\left(  \pi\right)  \rtimes_{\alpha
}\mathbb{Z}\right)  \cong\mathbb{Z}^{2}$ is surjective.

In these situations, if we start with an idempotent in $\mathbb{Q}\left[
\pi\right]  ,$ we obtain an idempotent in $\mathbb{Q}\left[  \pi\ast
\rho\right]  $ or $\mathbb{Q}\left[  \pi\rtimes_{\alpha}\mathbb{Z}\right]  .$

Let $\pi$ be any group and $p$ an idempotent in $\mathbb{Q}\left[  \pi\right]
$ representing a nonzero element of $\tilde{K}_{0}\left(  C^{\ast}\left(
\pi\right)  \right)  .$ We also denote by $p$ the corresponding multiplication
operator with image $P$. We construct a chain complex $C$ of $\mathbb{Z}%
\left[  \mathbb{\pi}\right]  \left[  z\right]  $-modules. For a suitable
integer $\ell,$ $\ell p$ is a module homomorphism which is defined
$\mathbb{Z}\left[  \mathbb{\pi}\right]  \mathbb{\rightarrow Z}\left[
\mathbb{\pi}\right]  .$ Let
\begin{align*}
C_{j}  &  =\left\{
\begin{array}
[c]{ll}%
\mathbb{Z}\left[  \mathbb{\pi}\right]  \left[  z\right]  , & j=0,1,\\
0 & \text{otherwise,}%
\end{array}
\right. \\
\partial &  =\ell\left(  I-zp\right)  .
\end{align*}

\noindent$\partial$ will in general have an infinitely generated cokernel of
exponent $\ell,$ so $C$ will not be finitely dominated. However,
$\partial\otimes I:C_{1}\otimes\mathbb{Q\rightarrow}C_{0}\otimes\mathbb{Q}$ is
injective with cokernel $P.$ First, $C\otimes\mathbb{Q}$ is chain equivalent
to the complex $I-zp:\mathbb{Q}\left[  \mathbb{\pi}\right]  \left[  z\right]
\mathbb{\rightarrow Q}\left[  \mathbb{\pi}\right]  \left[  z\right]  $ by
\[%
\begin{array}
[c]{ccc}%
\left(  \mathbb{Z}\left[  \mathbb{\pi}\right]  \mathbb{\otimes Q}\right)
\left[  z\right]  & \overset{\ell\left(  I-zp\right)  }{\longrightarrow} &
\left(  \mathbb{Z}\left[  \mathbb{\pi}\right]  \mathbb{\otimes Q}\right)
\left[  z\right] \\
I\downarrow &  & \downarrow1/\ell\\
\mathbb{Q}\left[  \mathbb{\pi}\right]  \left[  z\right]  & \overset
{I-zp}{\longrightarrow} & \mathbb{Q}\left[  \mathbb{\pi}\right]  \left[
z\right]  .
\end{array}
\]

\noindent We use the convention that $z^{-n}$ acts as $0$ on $\mathbb{Q}%
\left[  \mathbb{\pi}\right]  z^{j}$ if $n>j.$ Then $H=I-p\sum_{n=0}^{\infty
}z^{-n}$ satisfies $H\partial=I$ and $\partial H=I-\left[  p\;p\;p\cdots
\right]  ,$ where the vector goes in the first row. Therefore $C\otimes
\mathbb{Q}$ is equivalent to $P$ in degree 0. We also consider $C^{t},$ which
is the same except that $\partial=\ell\left(  I-z^{-1}\bar{p}\right)  .$ The
bar denotes conjugation in the group ring. $\partial\otimes I_{\mathbb{Q}}$ is
invertible with inverse $\ell^{-1}\left(  I+\bar{p}\sum_{n=1}^{\infty}%
z^{-n}\right)  .$

We will realize $C$ and $C^{t}$ geometrically. The following construction is
mostly due to Hughes and Ranicki \cite[Remark 10.3 (iii)]{hura}. Let $\pi$ be
any finitely presented group. For any $n\ge5$ there exists a paralellizable
manifold $L$ of dimension $n$ with boundary $V$ such that $\pi\left(
V\right)  =\pi\left(  N\right)  =\pi.$ We can embed a 2-complex with
fundamental group $\pi$ in $\mathbb{R}^{n}$ for $n\ge5$ and let $L$ be a
smooth regular neighborhood. Let $n\ge6.$

Let $N=S^{1}\times V.$ $N_{0}$ is the boundary component $N\times\left\{
0\right\}  $ of $N\times I.$ Attach a trivial 2-handle to $N\times\left\{
1\right\}  .$ The corresponding boundary component is the connected sum
$N^{\prime}=\left(  N\times\left\{  1\right\}  \right)  \#\left(  S^{2}\times
S^{n-3}\right)  .$ $\pi\left(  N^{\prime}\right)  \cong\pi\times\mathbb{Z}.$
Identify $\pi\left(  N^{\prime}\right)  $ with $\pi_{1}\left(  N^{\prime
}\right)  $ by choosing a basepoint and a lift of it to $\tilde{N}^{\prime}.$
Choose $h\in\pi_{2}\left(  N^{\prime}\right)  $ representing the cycle
$S^{2}\times\ast.$ Let $z$ be the generator of $\pi_{1}\left(  S^{1}\right)
.$ Attach a 3-handle using $\ell\left(  1-zp\right)  h.$ Let $\left(
W,N_{0},N_{1}\right)  $ be the resulting cobordism. $\pi\left(  N_{1}\right)
\cong\pi\left(  W\right)  \cong\pi\times\mathbb{Z}.$ We describe the complex
of the universal covers $C_{\ast}\left(  \tilde{W},\tilde{N}_{0}\right)  $
defined by the handle structure. Let $\tilde{h}$ correspond to $h$ under
$\pi_{2}\left(  N^{\prime}\right)  \cong\pi_{2}\left(  \tilde{N}^{\prime
}\right)  .$ $\tilde{h}$ represents $S^{2}\times\ast$ for some 2-handle
$e_{2}$ in $\tilde{W}.$ This handle generates $C_{2}\left(  \tilde{W}%
,\tilde{N}_{0}\right)  $ as a free left $\mathbb{Z}\left[  \pi\right]  \left[
z,z^{-1}\right]  $-module. $C_{3}\left(  \tilde{W},\tilde{N}_{0}\right)  $ is
freely generated by the handle $e_{3}$ attached by $\ell\left(  I-zp\right)
\tilde{h}.$ Therefore $\partial_{3}$ is given by $\partial e_{3}=\ell\left(
1-zp\right)  e_{2}.$ $\ell\left(  1-zp\right)  $ can also be described as the
$\mathbb{Z}\left[  \pi\right]  \left[  z,z^{-1}\right]  $-valued intersection
number $\mu\cdot\nu$ of the attaching sphere of $e_{3}$ with the transverse
sphere $\ast\times S^{n-3}$ of $e_{2}$ \cite[Sections II.6-II.8]{hud}. Now
consider the dual handle decomposition of $\left(  \tilde{W},\tilde{N}%
_{1}\right)  .$ This consists of handles of dimensions $n-2$ and $n-3.$ As
cells, these are the same as the original handles, but the attaching and
transverse spheres are interchanged. Therefore $\partial_{n-2}$ on $C_{\ast
}\left(  \tilde{W}^{+},\tilde{N}_{0}^{+}\right)  $ is given by $\nu\cdot\mu.$
In the present dimensions, $\nu\cdot\mu=\overline{\mu\cdot\nu}.$ It follows
that $\partial_{n-2}$ is given by $\ell\left(  1-z^{-1}\bar{p}\right)  .$

Let $\bar{W}$ be the infinite cyclic covering of $W$ classified by a map
$W\rightarrow S^{1}$ corresponding to $\pi\times\mathbb{Z\rightarrow Z}.$
$\bar{W}$ has the form $\left(  V\times\left[  0,1\right]  \times
\mathbb{R}\right)  \cup$\{handles indexed by $z^{n},$ $n\in\mathbb{Z}$\}.
$\bar{W}$ contains a subspace $\bar{W}^{+}$ diffeomorphic to $\left(
V\times\left[  0,1\right]  \times\lbrack0,\infty\right)  )\cup$\{handles
indexed by $z^{n},$ $n\geq0$\}. Let $\bar{N}_{0}^{+}$ and $\bar{N}_{1}^{+}$ be
the boundary components of $\partial\bar{W}^{+}-V\times\left(  0,1\right)
\times\left\{  0\right\}  .$ $\partial\bar{N}_{1}$ is diffeomorphic to $V.$
Let $M=\bar{N}_{1}\cup_{V}L,$ a manifold without boundary with $\pi\left(
M\right)  =\pi.$ We will show that $\tilde{\chi}_{\mathbb{Q}\left[
\pi\right]  }(M)=\left[  P\right]  .$

Note that $\left(  \bar{W}^{+}\right)  ^{\sim}=\tilde{W}^{+}$ and so on. From
the above, $C_{\ast}\left(  \tilde{W}^{+},\tilde{N}_{0}^{+}\right)  $ is the
complex $C$ with a dimension shift of 2, and the $\tilde{K}_{0}\left(
\mathbb{Q}\left[  \pi\right]  \right)  $-valued Euler characteristic of
$C_{\ast}\left(  \tilde{W}^{+},\tilde{N}_{0}^{+}\right)  \otimes\mathbb{Q}$ is
$\left[  P\right]  .$ $C_{\ast}\left(  \tilde{W}^{+},\tilde{N}_{1}^{+}\right)
$ is $C^{t}$ with a dimension shift, so the Euler characteristic of $C_{\ast
}\left(  \tilde{W}^{+},\tilde{N}_{1}^{+}\right)  \otimes\mathbb{Q}$ is 0. For
the following computations we use the chains of a smooth triangulation of
$\tilde{W}^{+}$ lifted from one of $\bar{W}^{+}.$ $\bar{N}_{0}^{+}=$
$V\times\lbrack0,\infty)$ is homotopy equivalent to $V.$ Therefore $C_{\ast
}\left(  \tilde{N}_{0}^{+}\right)  \otimes\mathbb{Q}$ is $\mathbb{C}\left[
\pi\right]  $-module equivalent to the finitely generated free complex
$C_{\ast}\left(  \tilde{V}\right)  \otimes\mathbb{Q},$ so represents
$0\in\tilde{K}_{0}\left(  \mathbb{Q}\left[  \pi\right]  \right)  .$ The sum
theorem for Euler characteristics \cite[Lemma 7]{wal2} applied to
\[
0\rightarrow C_{\ast}\left(  \tilde{N}_{0}^{+}\right)  \otimes
\mathbb{Q\rightarrow}C_{\ast}\left(  \tilde{W}^{+}\right)  \otimes
\mathbb{Q\rightarrow}C_{\ast}\left(  \tilde{W}^{+},\tilde{N}_{0}^{+}\right)
\otimes\mathbb{Q}\rightarrow0\mathbb{\;}%
\]

\noindent implies that $C_{*}\left(  \tilde{W}^{+}\right)  \otimes\mathbb{Q}$
represents $\left[  P\right]  .$ Then from
\[
0\rightarrow C_{*}\left(  \tilde{N}_{1}^{+}\right)  \otimes
\mathbb{Q\rightarrow}C_{*}\left(  \tilde{W}^{+}\right)  \otimes
\mathbb{Q\rightarrow}C_{*}\left(  \tilde{W}^{+},\tilde{N}_{1}^{+}\right)
\otimes\mathbb{Q}\rightarrow0,
\]

\noindent$C_{\ast}\left(  \tilde{N}_{1}^{+}\right)  \otimes\mathbb{Q}$
represents $\left[  P\right]  .$ The Mayer-Vietoris sequence%
\[
0\rightarrow C_{\ast}\left(  \tilde{V}\right)  \otimes\mathbb{Q}%
\rightarrow\left(  C_{\ast}\left(  \tilde{N}_{1}^{+}\right)  \oplus C_{\ast
}\left(  \tilde{L}\right)  \right)  \otimes\mathbb{Q\rightarrow}C_{\ast
}\left(  \tilde{M}\right)  \otimes\mathbb{Q\rightarrow}0
\]
shows that the Euler characteristic of $C_{\ast}\left(  \tilde{M}\right)
\otimes\mathbb{Q}$ in $\tilde{K}_{0}\left(  \mathbb{Q}\left[  \pi\right]
\right)  $ is $\left[  P\right]  .$

We wish to deal with right modules. From now on the above chain groups will be
equipped with the right action of the group ring defined by $ca=\bar{a}c.$
This change induces an equivalence between the categories of left and right
modules, so has no effect on the above computations. $\tilde{\chi}%
_{\mathbb{Q}\left[  \pi\right]  }$ was defined in terms of local coefficient
chains. The following well-known fact identifies these with chains of the
universal cover. Let $K$ be a simplicial complex and $\pi=\pi\left(  K\right)
.$ Let $\psi$ be the canonical bundle with fiber $\mathbb{Z}\left[
\pi\right]  .$

\begin{lemma}
\label{le4.1}There is an isomorphism of right $\mathbb{Z}\left[  \pi\right]
$-modules $C_{\ast}\left(  \tilde{K}\right)  \cong C_{\ast}\left(
K;\psi\right)  .$
\end{lemma}

%

%TCIMACRO{\TeXButton{Proof}{\proof}}%
%BeginExpansion
\proof
%EndExpansion
This is a simpler version of Section \ref{ss5.1}. A local coefficient
$j$-chain is a finitely-supported function which assigns to each $j$-simplex
of $K$ an element of the fiber of $\psi$ above its barycenter. Equivalently,
it is determined by a function $v$ from $j$-simplexes of $\tilde{K}$ to
$\mathbb{Z}\left[  \pi\right]  $ such that $v\left(  g\sigma\right)
=gv\left(  \sigma\right)  ,$ whose support intersects finitely many orbits of
$\pi.$ Let $S_{j}$ be the set of such functions. We define $vg$ by $vg\left(
\sigma\right)  =g^{-1}v\left(  g\sigma\right)  .$ Then $vg=v.$ For $u\in
C_{j}\left(  \tilde{K}\right)  $ let $\tau u=\sum_{g}\left(  ug^{-1}\right)
g.$ $\tau$ is an isomorphism to $S_{j}$. The inverse takes $v$ to
$C_{j}\left(  \tilde{K}\right)  \overset{v}{\rightarrow}$ $\mathbb{Z}\left[
\pi\right]  \rightarrow\mathbb{Z},$ where the last map is the component of the
identity of $\pi.$ Right multiplication in the fibers of $\psi$ by
$\mathbb{Z}\left[  \pi\right]  $ corresponds to right multiplication of values
of elements of $S_{\ast}.$ This corresponds under $\tau$ to the usual action
$ua\left(  \sigma\right)  =u\left(  a\sigma\right)  .$ These isomorphisms
commute with $\partial.$ This is clear for $\tau^{-1}.$ Consider the
isomorphism between $C_{\ast}\left(  K;\psi\right)  $ and $S_{\ast}.$ The
boundary for the first contains operators of parallel translation in $\psi$
along curves in $K.$ If a curve is lifted to $\tilde{K},$ the lift of the
parallel translation to $\tilde{K}\times\mathbb{Z}\left[  \pi\right]  $
projects to the identity of $\mathbb{Z}\left[  \pi\right]  .$%
%TCIMACRO{\TeXButton{End Proof}{\endproof}}%
%BeginExpansion
\endproof
%EndExpansion

As a consequence, $\tilde{\chi}_{\mathbb{Q}\left[  \pi\right]  }\left(
M\right)  =\left[  P\right]  .$ By Theorem 0.2, this construction gives a
manifold for which the index of $\bar{D}_{k}^{even}$ is $\left[  P\right]  $
for $k$ large.

Hughes and Ranicki \cite{hura} have introduced the locally finite finiteness
obstruction. If $C_{*}^{\ell f}\left(  X;\mathbb{Z}\left[  \pi\right]
\right)  $ is equivalent to a complex of finitely generated projective
modules, then its Euler characteristic is $o^{\ell f}\in K_{0}\left(
\mathbb{Z}\left[  \pi\right]  \right)  .$ It doesn't appear to have a direct
geometrical interpretation. If $C_{*}^{\ell f}\left(  X;\psi\right)  $ is
$C^{*}\left(  \pi\right)  $-finitely dominated, we refer to its Euler
characteristic in $K_{0}\left(  C^{*}\left(  \pi\right)  \right)  $ as
$\chi_{C^{*}\left(  \pi\right)  }^{\ell f}.$

\begin{lemma}
\label{le4.2}If $M^{n}$ is orientable and either $\chi_{C^{\ast}\left(
\pi\right)  }$ or $\chi_{C^{\ast}\left(  \pi\right)  }^{\ell f}$ is defined,
then so is the other, and $\chi_{C^{\ast}\left(  \pi\right)  }^{\ell
f}=\left(  -\right)  ^{n}\chi_{C^{\ast}\left(  \pi\right)  }$
\end{lemma}

%

%TCIMACRO{\TeXButton{Proof}{\proof}}%
%BeginExpansion
\proof
%EndExpansion
Duality gives an equivalence (up to sign) $C_{c}^{\ast}\left(  M;\psi\right)
\rightarrow C_{n-\ast}\left(  M;\psi\right)  .$ Therefore $C_{c}^{\ast}$ is
finitely dominated if and only if $C_{\ast}$ is. If so, $\chi\left(
C_{c}^{\ast}\right)  =\left(  -\right)  ^{n}\chi_{C^{\ast}\left(  \pi\right)
},$ since if $n$ is odd, duality exchanges the parities of the degrees.
$C_{\ast}^{\ell f}=\left(  C_{c}^{\ast}\right)  ^{\prime},$ so $C_{\ast}^{\ell
f}$ is finitely dominated if and only if $C_{c}^{\ast}$ is. Suppose that
$C_{c}^{\ast}$ is equivalent to the complex $F$ of finitely generated modules.
Then $C_{\ast}^{\ell f}$ is equivalent to $F^{\prime}.$ Since finitely
generated Hilbert modules are self-dual, $\chi\left(  C_{c}^{\ast}\right)
=\chi_{C^{\ast}\left(  \pi\right)  }^{\ell f}.$%
%TCIMACRO{\TeXButton{End Proof}{\endproof}}%
%BeginExpansion
\endproof
%EndExpansion

\section{\label{s5}Differential operators}

This section contains the proof that certain differential operators over
$C^{\ast}$-algebras are symmetric with nonnegative spectrum. This is a
generalization to bounded geometry manifolds of a special case of a theorem of
Kasparov. A proof is briefly sketched in \cite{kas1}. The one given here is
another application of weighted spaces.

\subsection{\textbf{\label{ss5.1}}}

Let $M$ be a manifold of bounded geometry and $E$ be an Hermitian vector
bundle over $M.$ Let $\pi$ be the group of covering transformations of a
normal covering space $\tilde{M}$. There is an Hilbert $C^{\ast}\left(
\pi\right)  $-module $\mathcal{E}$ associated to $E$ and $\tilde{M}$
\cite[Theorem 9.1]{kas2}, \cite[Section 1]{fac}. This is a reinterpretation of
the $\mathcal{L}^{2}$-type space associated to $E\otimes\psi.$

$\psi=\tilde{M}\times_{\pi}C^{\ast}\left(  \pi\right)  ,$ where the
equivalence relation is $\left(  x,a\right)  \sim\left(  gx,ga\right)  .$ Let
$\tilde{\psi}=\tilde{M}\times C^{\ast}\left(  \pi\right)  .$ The projection of
$\psi$ is induced from that of $\tilde{\psi},$ $\tilde{M}\times C^{\ast
}\left(  \pi\right)  \rightarrow\tilde{M}.$ $E$ lifts to a $\pi$-bundle
$\tilde{E}$ on $\tilde{M}.$ There is a one-to-one correspondence between
$C^{\infty}$ sections $v$ of $\tilde{E}\otimes\tilde{\psi}$ satisfying
$v\left(  gx\right)  =gv\left(  x\right)  $ and $C^{\infty}\left(
E\otimes\psi\right)  .$ Given $v$ and $y\in M,$ define $\left(  \kappa
v\right)  \left(  y\right)  $ to be the class of $\left(  x,v\left(  x\right)
\right)  ,$ where $x$ is any lift of $y.$ $\kappa$ clearly preserves the
$C^{\ast}\left(  \pi\right)  $-module structures defined by right
multiplication on fibers. The inverse $\lambda$ is given as follows. Let
$\ell_{x}$ be the canonical isomorphism of $\left(  E\otimes\psi\right)  _{y}$
with $\tilde{E}_{x}\otimes\tilde{\psi}_{x}\simeq\tilde{E}_{x}\otimes C^{\ast
}\left(  \pi\right)  $ given by the identifications. Then $\ell_{gx}=g\ell
_{x}.$ If $w$ is a section of $E\otimes\psi,$ let $\left(  \lambda w\right)
\left(  x\right)  =\ell_{x}\left(  w\left(  y\right)  \right)  .$ Then
$\left(  \lambda w\right)  \left(  gx\right)  =\ell_{gx}w\left(  y\right)
=g\ell_{x}w\left(  y\right)  =g\left(  \lambda w\right)  \left(  x\right)  .$
The $C^{\ast}\left(  \pi\right)  $-valued inner product on $C_{c}^{\infty
}\left(  E\otimes\psi\right)  $ corresponds to $\left(  u_{1}\otimes
u_{2},v_{1}\otimes v_{2}\right)  _{C^{\ast}}=\int_{F}\left\langle u_{1}\left(
x\right)  ,v_{1}\left(  x\right)  \right\rangle _{\tilde{E}}u_{2}\left(
x\right)  ^{\ast}v_{2}\left(  x\right)  dx,$ where $F$ is a fundamental
domain. If we write $\left(  vg\right)  \left(  x\right)  =g^{-1}v\left(
gx\right)  ,$ the invariance condition becomes $vg=v.$ If $u\in C_{c}^{\infty
}\left(  \tilde{E}\right)  ,$ let $\left(  \tau u\right)  =\sum_{g}%
ug^{-1}\otimes g\in C^{\infty}\left(  \tilde{E}\otimes\tilde{\psi}\right)  .$
$g$ denotes the constant section. It satisfies the condition since if $k\in
\pi,$
\[
\left(  \tau u\right)  k=\sum_{g}ug^{-1}k\otimes k^{-1}g=\sum_{g}u\left(
k^{-1}g\right)  ^{-1}\otimes k^{-1}g=\tau u.
\]
The action of $\mathbb{C}\left[  \pi\right]  $ on $C_{c}^{\infty}\left(
\tilde{E}\right)  $ extending $ug\left(  x\right)  =u\left(  gx\right)  $
corresponds to the $C^{\ast}\left(  \pi\right)  $ action on $C^{\infty}\left(
\tilde{E}\otimes\tilde{\psi}\right)  .$ The composition $\kappa\tau$ takes
$C_{c}^{\infty}\left(  \tilde{E}\right)  $ to $C_{c}^{\infty}\left(
E\otimes\psi\right)  .$ The induced inner product on $C_{c}^{\infty}\left(
\tilde{E}\right)  $ is
\begin{align}
\left(  u,v\right)  _{C^{\ast}}  &  =\int_{F}\sum_{g,h}\left\langle \left(
ug^{-1}\right)  \left(  x\right)  ,\left(  vh^{-1}\right)  \left(  x\right)
\right\rangle g^{-1}hdx\label{5.1}\\
&  =\int_{F}\sum_{g,h}\left\langle \left(  ugh\right)  \left(  x\right)
,\left(  vh\right)  \left(  x\right)  \right\rangle gdx\nonumber\\
&  =\int_{\tilde{M}}\sum_{g}\left\langle \left(  ug\right)  \left(  x\right)
,v\left(  x\right)  \right\rangle gdx=\sum_{g}\left(  ug,v\right)  g.\nonumber
\end{align}

\noindent Let $\mathcal{E}$ be the completion of $C_{c}^{\infty}\left(
\tilde{E}\right)  $ in the norm $\left\Vert u\right\Vert _{C^{\ast}%
}=\left\Vert \left(  u,u\right)  _{C^{\ast}}\right\Vert _{C^{\ast}\left(
\pi\right)  }^{1/2}.$ We will show that $\kappa\tau:C_{c}^{\infty}\left(
\tilde{E}\right)  \rightarrow C_{c}^{\infty}\left(  E\otimes\psi\right)  $ has
dense range with respect to the usual topology on $C_{c}^{\infty}$. It follows
that $\mathcal{E}$ may be identified with the completion of $C_{c}^{\infty
}\left(  E\otimes\psi\right)  .$ In particular, it is a Hilbert $C^{\ast
}\left(  \pi\right)  $-module.

An invariant section $v\in C^{\infty}\left(  \tilde{E}\otimes\tilde{\psi
}\right)  $ is called locally finite if it is of the form $\sum_{g}%
v_{g}\otimes g,$ where the supports of the $v_{g}$ are a locally finite
collection. Thus, if $u\in C_{c}^{\infty}\left(  \tilde{E}\right)  ,$ $\tau u$
is locally finite. If $v$ is locally finite,
\[
vh=\sum_{g}v_{g}h\otimes h^{-1}g=\sum_{g}v_{hg}h\otimes g=\sum_{g}v_{g}\otimes
g.
\]

\noindent Thus for all $g,$ $h,$ $v_{gh}h=v_{g}.$ Taking $h=g^{-1},$
$v_{e}g^{-1}=v_{g}.$ Therefore $v=\sum_{g}v_{e}g^{-1}\otimes g,$ and $v$ is
locally finite exactly when the translates of the support of $v_{e}$ are a
locally finite collection. $\tau$ extends to such $u=v_{e}.$

\begin{lemma}
\label{le5.1}If $v$ is locally finite, $\kappa v$ has compact support if and
only if $v_{e}$ does.
\end{lemma}

%

%TCIMACRO{\TeXButton{Proof}{\proof}}%
%BeginExpansion
\proof
%EndExpansion
Let $p:\tilde{M}\rightarrow M$ be the projection. By invariance, $Supp\left(
\kappa v\right)  =pSupp\left(  v\right)  .$ $Supp\left(  v\right)
=\overline{\cup_{g}gSupp\left(  v_{e}\right)  }$ . Since the $gSupp\left(
v_{e}\right)  $ are locally finite, this is $\cup_{g}gSupp\left(
v_{e}\right)  .$ Thus $Supp\left(  \kappa v\right)  =pSupp\left(
v_{e}\right)  ,$ and if $Supp\left(  v_{e}\right)  $ is compact, so is
$Supp\left(  \kappa v\right)  .$

$p|Supp\left(  v_{e}\right)  $ is finite-to-one. For if $Supp\left(
v_{e}\right)  $ contained infinitely many translates of some point, its
translates wouldn't be point finite. We show that $p|Supp\left(  v_{e}\right)
$ is a closed map. Let $V\subset Supp\left(  v_{e}\right)  $ be closed. Then
$\cup_{g}gV$ is closed since the $gV$ are locally finite. $pV=p\left(
\cup_{g}gV\right)  $ is closed since $M$ has the quotient topology.
$Supp\left(  v_{e}\right)  $ is then compact by a standard result
\cite[Exercise 26.12]{mun}.%
%TCIMACRO{\TeXButton{End Proof}{\endproof}}%
%BeginExpansion
\endproof
%EndExpansion

\begin{proposition}
\label{pr5.2}$\kappa\tau:C_{c}^{\infty}\left(  \tilde{E}\right)  \rightarrow
C_{c}^{\infty}\left(  E\otimes\psi\right)  $ has dense range.
\end{proposition}

%

%TCIMACRO{\TeXButton{Proof}{\proof}}%
%BeginExpansion
\proof
%EndExpansion
Let the sections with support in a set $K$ be $C_{K}^{\infty}\left(
E\otimes\psi\right)  .$ Let $B\subset M$ be a closed ball. By Lemma
\ref{le5.1}, the elements of $C_{B}^{\infty}\left(  E\otimes\psi\right)  $
which are images by $\kappa$ of locally finite invariant sections of
$\tilde{E}\otimes\tilde{\psi}$ come from elements of $C_{c}^{\infty}\left(
\tilde{E}\right)  .$ A choice of a lift of $B$ to $\tilde{M}$ determines a
trivialization $\psi|_{B}\cong B\times C^{\ast}\left(  \pi\right)  .$ Also
choose a trivialization $E|_{B}\cong B\times\mathbb{C}^{k}.$ The images of the
locally finite invariant sections correspond to the algebraic tensor product
$C_{B}^{\infty}\odot\left(  \mathbb{C}^{k}\otimes\mathbb{C}\left(  \pi\right)
\right)  .$ This has a unique tensor product topology \cite[II.3]{gro}.
$C_{B}^{\infty}\odot\left(  \mathbb{C}^{k}\otimes C^{\ast}\left(  \pi\right)
\right)  $ also has a unique tensor product, with completion $C_{B}^{\infty
}\left(  \mathbb{C}^{k}\otimes C^{\ast}\left(  \pi\right)  \right)  $. Since
$\mathbb{C}^{k}\otimes\mathbb{C}\left(  \pi\right)  $ is dense in
$\mathbb{C}^{k}\otimes C^{\ast}\left(  \pi\right)  ,$ $C_{B}^{\infty}%
\odot\left(  \mathbb{C}^{k}\otimes\mathbb{C}\left(  \pi\right)  \right)  $ is
dense in $C_{B}^{\infty}\left(  \mathbb{C}^{k}\otimes C^{\ast}\left(
\pi\right)  \right)  $. Therefore, the images of elements of $C_{c}^{\infty
}\left(  \tilde{E}\right)  $ are dense in $C_{B}^{\infty}\left(  E\otimes
\psi\right)  .$

Let $\left\{  U_{i}\right\}  $ be a locally finite cover of $M$ by open balls
with closures $B_{i},$ and $\left\{  \phi_{i}\right\}  $ a subordinate
partition of unity. Let $w\in C_{c}^{\infty}\left(  E\otimes\psi\right)  .$
Then the sum $w=\sum\phi_{i}w=\sum w_{i}$ is finite. Let $w_{ij}\in C_{B_{i}%
}^{\infty}\left(  E\otimes\psi\right)  $ be images of locally finite sections
such that $w_{ij}$ converges to $w_{i}.$ Then the sections $\sum_{i}w_{ij}$
are images of elements of $C_{c}^{\infty}\left(  \tilde{E}\right)  ,$ and
converge to $w$ in $C_{c}^{\infty}\left(  E\otimes\psi\right)  . $%
%TCIMACRO{\TeXButton{End Proof}{\endproof}}%
%BeginExpansion
\endproof
%EndExpansion

Let $F$ be another bundle with associated module $\mathcal{F}$, and $D$ a
first order linear differential operator $C_{c}^{\infty}\left(  E\right)
\rightarrow C_{c}^{\infty}\left(  F\right)  .$ Then $D$ lifts to an invariant
operator $\tilde{D}$: $C_{c}^{\infty}\left(  \tilde{E}\right)  \rightarrow
C_{c}^{\infty}\left(  \tilde{F}\right)  ,$ in the sense that $\tilde{D}\left(
ug\right)  =\left(  \tilde{D}u\right)  g.$ We will relate $\tilde{D}$ to the
operator $D^{\wedge}:C_{c}^{\infty}\left(  E\otimes\psi\right)  \rightarrow
C_{c}^{\infty}\left(  F\otimes\psi\right)  $, $D$ with coefficients in $\psi.$
We recall the construction \cite[4.2]{mil1},\cite[IV.9]{pal}.

Let $\nabla^{E}$ be a unitary connection on $E.$ $D$ may be expressed as a
locally finite sum $D=B_{0}+\sum_{j>0}B_{j}\nabla_{X_{j}}^{E},$ where
$B_{j}\in C^{\infty}\left(  Hom\left(  E,F\right)  \right)  ,$ $X_{j}\in
C^{\infty}\left(  TM\right)  .$ Let $\nabla^{\psi}$ be the flat connection on
$\psi.$ Let $\nabla=\nabla^{E}\otimes I_{\psi}+I_{E}\otimes\nabla^{\psi}.$
Define $D^{\wedge}=B_{0}\otimes I_{\psi}+\sum_{j>0}\left(  B_{j}\otimes
I_{\psi}\right)  \nabla_{X_{j}}.$ This is independent of $\nabla^{E}.$ The
construction preserves formal adjoints. Using local sections of the covering
projection, all the elements of structure lift to $\tilde{M}$ to define
$\tilde{D}^{\wedge}.$ It is evident that $\widetilde{D^{\wedge}}=\tilde
{D}^{\wedge}$ and that for an invariant section $v,$ $\kappa\left(  \tilde
{D}^{\wedge}v\right)  =D^{\wedge}\left(  \kappa v\right)  .$ Since $\psi$ is
flat, $\tilde{\nabla}_{\tilde{X}_{j}}\left(  v_{g}\otimes g\right)  =\left(
\tilde{\nabla}_{\tilde{X}_{j}}^{\tilde{E}}v_{g}\right)  \otimes g,$ so
$\tilde{D}^{\wedge}\left(  v_{g}\otimes g\right)  =\left(  \tilde{D}%
v_{g}\right)  \otimes g.$ If $u\in C_{c}^{\infty}\left(  \tilde{E}\right)  ,$%
\[
\tau\left(  \tilde{D}u\right)  =\sum_{g}\left(  \tilde{D}u\right)
g^{-1}\otimes g=\sum_{g}\tilde{D}\left(  ug^{-1}\right)  \otimes g=\sum
_{g}\tilde{D}^{\wedge}\left(  ug^{-1}\otimes g\right)  =\tilde{D}^{\wedge
}\left(  \tau u\right)  .
\]

\noindent Therefore we may identify the operators $\tilde{D}$ and $D^{\wedge}
$ under the above identification of Hilbert modules.

\subsection{\textbf{\label{ss5.2}}}

We will assume that the principal symbol of $D$ is uniformly bounded in norm.
Let $\tilde{D}^{\#}$ be the formal adjoint of $\tilde{D}$ with respect to the
ordinary $\mathcal{L}^{2}$ inner products. Let
\[
T=\left(
\begin{array}
[c]{rr}%
0 & \tilde{D}^{\#}\\
\tilde{D} & 0
\end{array}
\right)  :C_{c}^{\infty}\left(  \tilde{E}\oplus\tilde{F}\right)  \rightarrow
C_{c}^{\infty}\left(  \tilde{E}\oplus\tilde{F}\right)  .
\]

\noindent The principal symbol of $T$ is also uniformly bounded. $T$ is
symmetric for the $C^{\ast}$-inner product. For
\[
\left(  Tu,v\right)  _{C^{\ast}}=\sum_{g}\left(  \left(  Tu\right)
g,v\right)  g=\sum_{g}\left(  T\left(  ug\right)  ,v\right)  g=\sum_{g}\left(
ug,Tv\right)  g=\left(  u,Tv\right)  _{C^{\ast}}.
\]

\noindent Thus $\bar{T},$ the closure of $T$ for the $C^{*}$-norm, is
symmetric. By an easy argument, the adjoint of a closable operator is equal to
the adjoint of its closure \cite[Vol. 1, Th. 4.1.3]{kari}.

\begin{theorem}
\label{th5.3}$\tilde{D}^{\ast}\overline{\tilde{D}}$ is symmetric with real
spectrum contained in $[0,\infty).$
\end{theorem}

We use this terminology rather than \textquotedblleft
self-adjoint\textquotedblright\ since self-adjoint operators over $C^{\ast}%
$-algebras need not have real spectrum \cite{hil}. The main point is to show
that $\bar{T}\pm\lambda i$ has dense range for some $\lambda>0.$ The proof
involves comparing $\bar{T}$ and the closures of $T$ on weighted spaces. For
the present, $\lambda$ is a free parameter which eventually will be chosen to
be sufficiently large. Until further notice we consider the closure of $T$ as
an operator on $\mathcal{L}^{2}\left(  \tilde{E}\oplus\tilde{F}\right)  ,$
still denoted $\bar{T}.$ According to Chernoff \cite{che}, $T$ is essentially
self-adjoint. Let $x_{0}\in\tilde{M}$ be a fixed point, and $d\left(
x,x_{0}\right)  $ be the distance function. Gaffney has shown that there
exists a $C^{\infty}$ function $\rho\left(  x\right)  $ such that $\left\vert
d\left(  x,x_{0}\right)  -\rho\left(  x\right)  \right\vert $ is bounded and
$\left\Vert d\rho\left(  x\right)  \right\Vert $ is bounded \cite[Lemma
A1.2.1]{shu}. Let $\sigma_{T}$ be the principal symbol of $T$ and
$\delta=\left(  \sup_{x\in\tilde{M}}\left\Vert \sigma_{T}\left(
x,d\rho\left(  x\right)  \right)  \right\Vert \right)  .$ Let $\mathcal{L}%
_{k}^{2}$ be the completion of $C_{c}^{\infty}\left(  \tilde{E}\oplus\tilde
{F}\right)  $ in the inner product with weight function $k^{\rho\left(
x\right)  }$. Let $\bar{T}_{k}$ be the closure of $T$ acting on $\mathcal{L}%
_{k}^{2}.$ The following argument is well known.

\begin{lemma}
\label{le5.4}$\bar{T}_{k}\pm i\lambda$ is boundedly invertible if $|\log
k|<\delta\lambda.$
\end{lemma}

%

%TCIMACRO{\TeXButton{Proof}{\proof} }%
%BeginExpansion
\proof
%EndExpansion
Multiplication by $k^{\rho\left(  x\right)  }$ induces a unitary
$\mathcal{L}_{k}^{2}\rightarrow\mathcal{L}^{2}.$ $\bar{T}_{k}\pm i\lambda$ is
unitarily equivalent to the closure of
\[
k^{\rho\left(  x\right)  }\left(  T\pm i\lambda\right)  k^{-\rho\left(
x\right)  }=T+\left(  \log k\right)  \sigma\left(  x,d\rho\left(  x\right)
\right)  \pm i\lambda
\]
acting on $\mathcal{L}^{2},$ which is $\bar{T}+\left(  \log k\right)
\sigma\left(  x,d\rho\left(  x\right)  \right)  \pm i\lambda.$ Since $\bar{T}$
is self-adjoint, $\bar{T}\pm i\lambda$ is boundedly invertible and
\[
\left\Vert \left(  \log k\right)  \sigma\left(  x,d\rho\left(  x\right)
\right)  \left(  \bar{T}\pm i\lambda\right)  ^{-1}\right\Vert \leq\left\vert
\log k\right\vert \delta\lambda^{-1}%
\]
by \cite[Theorem 5.18]{wei}. This is $<1$ provided that $|\log k|<\delta
\lambda$ and then
\[
\bar{T}+\left(  \log k\right)  \sigma\left(  x,d\rho\left(  x\right)  \right)
\pm i\lambda
\]
is boundedly invertible. Therefore $\bar{T}_{k}\pm i\lambda$ is boundedly
invertible.%
%TCIMACRO{\TeXButton{End Proof}{\endproof}}%
%BeginExpansion
\endproof
%EndExpansion

The next Lemma gives the basic relationship between the norms on $\mathcal{E}
$ and $\mathcal{L}_{k}^{2}.$ The proof indicates the relationship between $k$
and the growth rate of $\tilde{M}.$

\begin{lemma}
\label{le5.5}Let $u\in C_{c}^{\infty}\left(  \tilde{E}\oplus\tilde{F}\right)
. $ Then for all sufficiently large $k,$ $\left\Vert u\right\Vert _{C^{\ast}%
}\leq K\left\Vert u\right\Vert _{k},$ where $K$ depends only on $k. $
\end{lemma}

%

%TCIMACRO{\TeXButton{Proof}{\proof}}%
%BeginExpansion
\proof
%EndExpansion
The $\mathcal{L}^{1}$ norm on $\mathbb{C}\left(  \pi\right)  $ is $\left\Vert
a\right\Vert _{\mathcal{L}^{1}\left(  \pi\right)  }=\sum_{g\in\pi}\left\vert
a\left(  g\right)  \right\vert .$ It majorizes the $C^{\ast}$norm. Let
$w=k^{\rho\left(  x\right)  }u.$ Then
\begin{align*}
\left\Vert u\right\Vert _{C^{\ast}}^{2}  &  =\left\Vert \sum_{g}\left(
ug,u\right)  g\right\Vert _{C^{\ast}\left(  \pi\right)  }\leq\left\Vert
\sum_{g}\left(  ug,u\right)  g\right\Vert _{\mathcal{L}^{1}\left(  \pi\right)
}\\
&  =\sum_{g}\left\vert \left(  ug,u\right)  \right\vert \leq\sum_{g}\left(
\left\vert \left(  wg,w\right)  \right\vert \sup_{x\in\tilde{M}}k^{-\left(
\rho\left(  x\right)  +\rho\left(  gx\right)  \right)  }\right)  .
\end{align*}

\noindent Since $\left\vert d\left(  x,x_{0}\right)  -\rho\left(  x\right)
\right\vert $ is bounded, there is a $C$ such that $k^{-\left(  \rho\left(
x\right)  +\rho\left(  gx\right)  \right)  }\leq Ck^{-\left(  d\left(
x_{0},x\right)  +d(x_{0},gx\right)  }$ for all $x.$ Then the last expression
above is less than or equal to%

\begin{align*}
C\sum_{g}\left(  \left\vert \left(  wg,w\right)  \right\vert \sup_{x\in
\tilde{M}}k^{-\left(  d\left(  x_{0},x\right)  +d(x_{0},gx\right)  }\right)
&  \leq C\left\Vert w\right\Vert ^{2}\sum_{g}k^{-d\left(  x_{0},gx_{0}\right)
}\\
&  =C\left\Vert u\right\Vert _{k}^{2}\sum_{g}k^{-d\left(  x_{0},gx_{0}\right)
}.
\end{align*}

\noindent The next to last step follows from the Cauchy inequality and the
fact that $d\left(  x_{0},x\right)  +d\left(  x_{0},gx\right)  \geq d\left(
x_{0},gx_{0}\right)  .$ We will show that the series converges for $k$
sufficiently large.

We claim that the number of points $N\left(  r\right)  $ in any orbit of $\pi$
on $\tilde{M}$ lying in a ball $B$ of radius $r$ is bounded by $e^{cr}$ for
some $c.$ From the condition on the injectivity radius, it follows that there
exists $\epsilon>0$ such that $d\left(  x_{1},x_{2}\right)  >2\epsilon$ for
any $x_{1},$ $x_{2}$ in the orbit. For any $\epsilon>0$ there is a minimum
volume $V\left(  \epsilon\right)  $ for balls of radius $\epsilon$ \cite[Lemma
A1.1.3]{shu}. The volume of $B$ satisfies $Vol\left(  B\right)  <e^{mr}$ for
some $m.$ Now $N\left(  r\right)  V\left(  \epsilon\right)  <Vol\left(
B\right)  ,$ so $N\left(  r\right)  <\dfrac{Vol\left(  B\right)  }{V\left(
\epsilon\right)  }<\dfrac{e^{mr}}{V\left(  \epsilon\right)  }.$ We consider
balls of radius $n\in\mathbb{N}$ with center $x_{0}.$ Then
\[
\sum_{g}k^{-d\left(  x_{0},gx_{0}\right)  }\leq\sum_{n=1}^{\infty}k^{-\left(
n-1\right)  }e^{cn}=k\sum e^{\left(  c-\log k\right)  n},
\]

\noindent and the last series converges for $k>e^{c}.$%
%TCIMACRO{\TeXButton{End Proof}{\endproof}}%
%BeginExpansion
\endproof
%EndExpansion

Let $T_{C^{*}}$ be $T$ acting on $C_{c}^{\infty}$ with the inner product
$\left(  \cdot,\cdot\right)  _{C^{*}},$ and $\bar{T}_{C^{*}}$ its closure.

\begin{lemma}
\label{le5.6}For $k$ sufficiently large, $\bar{T}_{C^{\ast}}$ is an extension
of $\bar{T}_{k}.$
\end{lemma}

%

%TCIMACRO{\TeXButton{Proof}{\proof}}%
%BeginExpansion
\proof
%EndExpansion
A bounded operator between normed spaces extends to an operator between their
completions with the same norm. By \ref{le5.5} for $k$ large the identity map
of $C_{c}^{\infty}$ with the $k$- and $C^{\ast}$-norms extends to
$\mathcal{L}_{k}^{2}\rightarrow\mathcal{E.}$ The identity on $C_{c}^{\infty}$
extends to bounded maps $\mathcal{L}_{k}^{2}\rightarrow\mathcal{L}^{2}$ for
any $k>0,$ since $\left\Vert u\right\Vert _{k}\geq\left\Vert u\right\Vert .$
The pointwise inner product $\left\langle u,u\right\rangle $ on $C_{c}%
^{\infty}$ extends to an $\mathcal{L}^{1}$ function of $u\in\mathcal{L}^{2}.$
If $\left(  u,u\right)  _{k}=\int\left\langle u,u\right\rangle k^{2\rho\left(
x\right)  }dx>0,$ Then $\left(  u,u\right)  =\int\left\langle u,u\right\rangle
dx>0.$ Therefore the maps are injective.

The maps $\mathcal{L}_{k}^{2}\rightarrow\mathcal{E}$ are injective. This
follows from a factorization of $\mathcal{L}_{k}^{2}\rightarrow\mathcal{L}%
^{2}$ as $\mathcal{L}_{k}^{2}\rightarrow\mathcal{E}\rightarrow\mathcal{L}%
^{2}.$ There is a bounded trace $Tr:C^{\ast}\left(  \pi\right)  \rightarrow
\mathbb{C}$ which on elements of $\mathbb{C}\left[  \pi\right]  $ is the
coefficient of $e.$ By \ref{5.1} for $u\in C_{c}^{\infty},$ $\left(
u,u\right)  =Tr$ $\left(  u,u\right)  _{C^{\ast}}.$ Then%
\[
\left\Vert u\right\Vert ^{2}=\left(  u,u\right)  =Tr\left(  u,u\right)
_{C^{\ast}}\leq K\left\Vert \left(  u,u\right)  \right\Vert _{C^{\ast}%
}=K\left\Vert u\right\Vert _{C^{\ast}}^{2}.
\]

\allowbreak\noindent This provides the map $\mathcal{E}\rightarrow
\mathcal{L}^{2}.$ It follows directly that $\mathcal{D}\left(  \bar{T}%
_{k}\right)  $ is identified with a subset of $\mathcal{D}\left(  \bar
{T}_{C^{\ast}}\right)  $ and $\bar{T}_{C^{\ast}}=$ $\bar{T}_{k}$ on
$\mathcal{D}\left(  \bar{T}_{k}\right)  .$%
%TCIMACRO{\TeXButton{End Proof}{\endproof}}%
%BeginExpansion
\endproof
%EndExpansion

\noindent In general, $Tr$ isn't faithful on $C^{*}\left(  \pi\right)  ,$ so
$\mathcal{E}$ isn't a subspace of $\mathcal{L}^{2}.$ It is if $C^{*}\left(
\pi\right)  $ is replaced by the reduced algebra.

A regular operator on a Hilbert module is a closed operator $A$ with dense
domain such that $A^{*}$ has dense domain and $A^{*}A+I$ is surjective.

\textit{Proof of Theorem \ref{th5.3}}. Choose $k$ so that $\bar{T}_{C^{\ast}}
$ is an extension of $\bar{T}_{k},$ then $\lambda$ so that $\bar{T}_{k}\pm
i\lambda$ is boundedly invertible, so surjective. Then $\bar{T}_{C^{\ast}}\pm
i\lambda$ has dense range, and is boundedly invertible since $\bar{T}%
_{C^{\ast}}$ is symmetric \cite[Theorem 5.18]{wei}. Henceforth, symbols like
$\bar{T}$ are closures in the $C^{\ast}$-norm. Since $\bar{T}$ is symmetric
and $\bar{T}\pm i\lambda$ is boundedly invertible, $\bar{T}+z$ is boundedly
invertible for all nonreal $z$ \cite[Theorem 5.21]{wei}. $\left(  \bar
{T}+i\right)  \left(  \bar{T}-i\right)  =\bar{T}^{2}+I,$ so $\bar{T}^{2}+I$ is
surjective. $\bar{T}$ is self adjoint \cite[Theorem 5.21]{wei}, so $T^{\ast
}\bar{T}+I$ is surjective.
\[
T^{\ast}\bar{T}=\left[
\begin{array}
[c]{ll}%
\tilde{D}^{\ast}\overline{\tilde{D}} & 0\\
0 & \tilde{D}^{\#\ast}\overline{\tilde{D}^{\#}}%
\end{array}
\right]  ,
\]

\noindent so $\tilde{D}^{*}\overline{\tilde{D}}+I$ is surjective.
$\overline{\tilde{D}}$ is thus a regular operator. By \cite[Proposition
9.9]{lan}, $\tilde{D}^{*}\overline{\tilde{D}}$ is self adjoint, and thus
closed. By \cite[Proposition 2.5]{mil2}, it has spectrum in $[0,\infty).$%
%TCIMACRO{\TeXButton{End Proof}{\endproof}}%
%BeginExpansion
\endproof
%EndExpansion

In the remainder of this section we will consider invariant operators like
$\tilde{D}$ exclusively. For notational convenience the tildes will be omitted.

\subsection{\textbf{\label{ss5.3}}}

We need more information in some special cases. $D_{\mu}=d+\delta_{\mu}$ is
unitarily equivalent to $d+\delta-\left(  dh\wedge+dh\llcorner\right)  $
acting on $\bar{\Omega}.$We use this operator to form $T.$The principal symbol
is given by Clifford multiplication, so $\left\Vert \sigma\left(
x,\cdot\right)  \right\Vert =1.$ $d$ and $\delta_{\mu}$ are handled similarly.
Since $\bar{T}$ is self-adjoint in each case, $\tilde{D}_{\mu}^{\ast
}=\overline{\tilde{D}}_{\mu},$ $\tilde{d}_{\mu}^{\ast}=\overline{\tilde
{\delta}}_{\mu},$ and $\tilde{\delta}_{\mu}^{\ast}=\overline{\tilde{d}}.$ We
suppress the tildes from now on. By Theorem \ref{th5.3}, $\bar{D}_{\mu}%
^{2}=D_{\mu}^{\ast}\bar{D}_{\mu},$ $d_{\mu}^{\ast}\bar{d},$ and $\bar{d}%
d_{\mu}^{\ast}$ are symmetric with spectrum in $[0,\infty).$ When the presence
of weighting makes no difference, we will omit the subscript $\mu.$ Since the
images of $d$ and $\delta$ are orthogonal, it follows that $\bar{D}=\bar
{d}+\bar{\delta}=\bar{d}+d^{\ast}.$ $\bar{D}^{2}=D^{\ast}\bar{D}=\bar
{d}d^{\ast}+d^{\ast}\bar{d},$ since $\operatorname{Im}\bar{d}\subset\ker
\bar{d}$ and $\operatorname{Im}\bar{\delta}\subset\ker\bar{\delta}.$

\begin{lemma}
\label{le5.7}Let $f\left(  t\right)  \in C\left(  \text{Spec}~\left(  d^{\ast
}\bar{d}+I\right)  ^{-1}\right)  $ or $C\left(  \text{Spec}~\left(  \bar
{d}d^{\ast}+I\right)  ^{-1}\right)  $ as is appropriate. Then
\end{lemma}

\begin{enumerate}
\item $f\left(  \left(  d^{*}\bar{d}+I\right)  ^{-1}\right)  \bar{d}=f\left(
1\right)  \bar{d}$

\item $\bar{d}f\left(  \left(  \bar{d}d^{*}+I\right)  ^{-1}\right)  =f\left(
1\right)  \bar{d}.$

\item $\bar{d}f\left(  \left(  \bar{D}^{2}+I\right)  ^{-1}\right)
=\overline{f\left(  \left(  \bar{D}^{2}+I\right)  ^{-1}\right)  \bar{d}}$

\item $d^{*}f\left(  \left(  \bar{D}^{2}+I\right)  ^{-1}\right)
=\overline{f\left(  \left(  \bar{D}^{2}+I\right)  ^{-1}\right)  d^{*}}$
\end{enumerate}

%

%TCIMACRO{\TeXButton{Proof}{\proof}}%
%BeginExpansion
\proof
%EndExpansion
(1) and (2). We prove the first. Since $\left(  d^{*}\bar{d}+I\right)  \bar
{d}=\bar{d},$ $\left(  d^{*}\bar{d}+I\right)  ^{-1}\bar{d}=\bar{d}.$ By
continuity we may assume $f$ smooth and write $f\left(  t\right)  =f\left(
1\right)  +g\left(  t\right)  (t-1).$ Then
\[
f\left(  \left(  d^{*}\bar{d}+I\right)  ^{-1}\right)  \bar{d}=f\left(
1\right)  \bar{d}+g\left(  \left(  d^{*}\bar{d}+I\right)  ^{-1}\right)
\left(  \left(  d^{*}\bar{d}+I\right)  ^{-1}-I\right)  \bar{d}=f\left(
1\right)  \bar{d}.
\]

(3) and (4) are well known. They are proved by approximating $f$ by a sequence
of polynomials and using the relations $\bar{d}\left(  \bar{D}^{2}+I\right)
=\left(  \bar{D}^{2}+I\right)  \bar{d}$ and $d^{*}\left(  \bar{D}%
^{2}+I\right)  =\left(  \bar{D}^{2}+I\right)  d^{*}.$%
%TCIMACRO{\TeXButton{End Proof}{\endproof}}%
%BeginExpansion
\endproof
%EndExpansion

We establish the properties of the complexes $E_{\mu}$ with differentials
$d_{E_{\mu}}=\bar{d}\left(  \bar{D}_{\mu}^{2}+I\right)  ^{-1/2}$ of section
\ref{ss1.4}.

$d_{E}^{2}\subset\bar{d}\bar{d}\left(  \bar{D}_{\mu}^{2}+I\right)
^{-1/2}\left(  \bar{D}_{\mu}^{2}+I\right)  ^{-1/2}=0$ by Lemma \ref{le5.7}(3).

$d_{E}$ is bounded: by \cite[Proposition 2.6]{mil2}, $\mathcal{D}\left(
\left(  \bar{D}_{\mu}^{2}+I\right)  ^{1/2}\right)  =\mathcal{D}\left(  \bar
{D}_{\mu}\right)  ,$ so $\operatorname{Im}\left(  \bar{D}_{\mu}^{2}+I\right)
^{-1/2}$ $\subset\mathcal{D}\left(  \bar{d}\right)  .$ The conclusion follows
from \cite[Exercise 5.6]{wei}. Also
\begin{equation}
d_{E_{\mu}}^{\ast}=\left(  \bar{d}\left(  \bar{D}_{\mu}^{2}+I\right)
^{-1/2}\right)  ^{\ast}=\left(  \left(  \bar{D}_{\mu}^{2}+I\right)
^{-1/2}\bar{d}\right)  ^{\ast}=d_{\mu}^{\ast}\left(  \bar{D}_{\mu}%
^{2}+I\right)  ^{-1/2} \label{5.2}%
\end{equation}

\noindent since $\left(  \bar{D}_{\mu}^{2}+I\right)  ^{-1/2}$ is bounded.

We establish an isomorphism the between the complexes of differential forms
$\left(  \bar{\Omega}_{d,\mu},\bar{d}\right)  $ and $\left(  \bar{\Omega}%
_{\mu},d_{E_{\mu}}\right)  .$ $\left(  d^{*}\bar{d}+I\right)  ^{1/2}$ is a
unitary between $\bar{\Omega}_{d}$ and $\bar{\Omega}:$ by \cite[Proposition
2.6]{mil2}, $\mathcal{D}\left(  \left(  d^{*}\bar{d}+I\right)  ^{1/2}\right)
=\mathcal{D}\left(  \bar{d}\right)  $ and
\[
\left(  u,v\right)  _{d}=\left(  u,v\right)  +\left(  \bar{d}u,\bar
{d}v\right)  =\left(  \left(  d^{*}\bar{d}+I\right)  ^{1/2}u,\left(  d^{*}%
\bar{d}+I\right)  ^{1/2}v\right)  .
\]

\noindent The isomorphism will follow from the fact that $\left(  d^{*}\bar
{d}+I\right)  ^{1/2}$ is a cochain isomorphism, i.e.
\[
\left(  d^{*}\bar{d}+I\right)  ^{-1/2}\bar{d}\left(  \bar{D}^{2}+I\right)
^{-1/2}\left(  d^{*}\bar{d}+I\right)  ^{1/2}=\bar{d}.
\]

\noindent By Lemma \ref{le5.7}(1), the left side is $\bar{d}\left(  \bar
{D}^{2}+I\right)  ^{-1/2}\left(  d^{\ast}\bar{d}+I\right)  ^{1/2}.$ Since
\[
\left(  \bar{D}^{2}+I\right)  ^{-1/2}=\left(  \bar{d}d^{\ast}+I\right)
^{-1/2}\left(  d^{\ast}\bar{d}+I\right)  ^{-1/2},
\]
using Lemma \ref{le5.7}(3) it is
\[
\bar{d}\left(  \bar{d}d^{\ast}+I\right)  ^{-1/2}\left(  d^{\ast}\bar
{d}+I\right)  ^{-1/2}\left(  d^{\ast}\bar{d}+I\right)  ^{1/2}=\bar{d}\left(
d^{\ast}\bar{d}+I\right)  ^{-1/2}\left(  d^{\ast}\bar{d}+I\right)  ^{1/2}%
=\bar{d}.
\]
The last equality holds since $\mathcal{D}\left(  \left(  d^{\ast}\bar
{d}+I\right)  ^{1/2}\right)  =\mathcal{D}\left(  \bar{d}\right)  .$

Now consider the complexes $E_{\mu}$ with the modified differentials
$\beta_{\mu}$ and unitaries $\tau_{\mu}.$ The above shows that $\beta_{\mu}$
is bounded and $\beta_{\mu}^{2}=0.$

\begin{lemma}
\label{le5.8}$\tau_{\mu}\beta_{\mu}\tau_{\mu^{-}}=\beta_{\mu^{-}}^{\ast}.$
\end{lemma}

%

%TCIMACRO{\TeXButton{Proof}{\proof}}%
%BeginExpansion
\proof
%EndExpansion
On $\Omega_{c},$%
\begin{align}
\left(  e^{2h}\ast\right)  d\left(  e^{-2h}\ast\right)   &  =\left(  -\right)
^{nj+n+1}e^{2h}\delta e^{-2h}=\left(  -\right)  ^{nj+n+1}\delta_{\mu^{-}%
},\label{5.3}\\
\left(  e^{2h}\ast\right)  \delta_{\mu}\left(  e^{-2h}\ast\right)   &
=e^{2h}\left(  e^{-2h}\ast\delta\ast e^{2h}\right)  e^{-2h}=\left(  -\right)
^{nj+n}d.\nonumber
\end{align}

\noindent By a standard calculation, $e^{2h}\ast D_{\mu}^{2}=D_{-\mu}%
^{2}e^{2h}\ast,$ so $\tau_{\mu}D_{\mu}^{2}=D_{\mu^{-}}^{2}\tau_{\mu}. $ Then
\[
\tau_{\mu}\left(  \bar{D}_{\mu}^{2}+I\right)  ^{-1}\tau_{\mu^{-}}=\left(
\bar{D}_{\mu^{-}}^{2}+I\right)  ^{-1}.
\]
If $p\left(  t\right)  $ is a polynomial, it follows that $\tau_{\mu}p\left(
\left(  \bar{D}_{\mu}^{2}+I\right)  ^{-1}\right)  \tau_{\mu^{-}}=p\left(
\left(  \bar{D}_{\mu^{-}}^{2}+I\right)  ^{-1}\right)  .$ Therefore $\tau_{\mu
}\left(  \bar{D}_{\mu}^{2}+I\right)  ^{-1/2}\tau_{\mu^{-}}=$ $\left(  \bar
{D}_{\mu^{-}}^{2}+I\right)  ^{-1/2}.$ The conclusion follows from (\ref{5.2}),
(\ref{5.3}), and a check of conventions.%
%TCIMACRO{\TeXButton{End Proof}{\endproof}}%
%BeginExpansion
\endproof
%EndExpansion

\section{\label{s6}Discussion}

The purpose of this section is to explain connections between this paper and
other work on analysis and algebraic topology on manifolds with periodic or
approximately periodic ends. The contents of this paper represent a hybrid of
the two approaches. The main theme is the connection between finite
domination, the Fredholm property, and contractibility of complexes. Results
and notation from the rest of the paper will be used freely. In this section
the $C^{\ast}$-algebra $A$ is $\mathbb{C}$ unless otherwise stated. The main
results aren't known to hold for general $A.$

The fundamental fact concerning index theory on complete manifolds is due to
Anghel \cite{ang}. We state it in its original form. It can be generalized to
complexes. Consider an essentially self-adjoint first order elliptic
differential operator acting on an Hermitian bundle. Let $D$ be its closure, a
bounded operator in the graph norm $\left\|  \cdot\right\|  _{D}.$

\begin{theorem}
\label{th6.1}\cite[Theorem 2.1]{ang} $D$ is Fredholm if and only if there is a
constant $c>0$ and a compact subset $K\subset M$ such that $\left\Vert
Du\right\Vert \geq c\left\Vert u\right\Vert _{D}$ if $u\in\mathcal{D}\left(
D\right)  $ and $Supp\left(  u\right)  \cap K=\emptyset.$
\end{theorem}

The hypothesis of the Theorem is sometimes referred to as invertibility at
infinity. Observe that if $D$ is invariant under a proper isometric action of
$\mathbb{Z}$, then $K$ must be empty. Therefore $D$ is Fredholm if and only if
it is invertible. (This was first proved by Eichhorn.) In earlier work,
versions of this fact were proved. It was applied after an excision argument
to reduce to a periodic situation. (In the present paper, this step
corresponds to Lemma \ref{le3.2}.)

Theorem \ref{th6.1} has been applied to operators which are the sum of a
generalized Dirac operator and a potential. The potentials are vector bundle
maps which are fiberwise strictly positive on the complement of a compact set.
(Most of the relevant papers are in the bibliography of \cite{foha}.) The
operators in the present paper are of the form $d+\delta-\left(  2\log
k\right)  d\rho\llcorner$. Theorem \ref{th0.1} states that if $M$ has finitely
many quasi-periodic ends and finitely generated rational homology, then the
operator is Fredholm for certain values of $k.$ The set of critical points of
$\rho$ can be compact only if $M$ admits a boundary. We have therefore shown
that even if this is not the case, the operator may nonetheless be invertible
at infinity. Section 4.1 contains a relevant example.

The first work related to this paper, by Lockhart and McOwen \cite{lomc} and
Melrose and Mendoza, concerned manifolds with cylindrical ends. However, the
subsequent results of Taubes represent a proper generalization, so we discuss
these first. Let $M$ be a smooth manifold with finitely many periodic ends.
For simplicity, we consider the case of one end. Let $\bar{N}^{+}\subset
\bar{N}$ be the model for the end, where $\bar{N}$ is an infinite cyclic
covering of the compact manifold $N.$ Let $C=\left\{  C_{c}^{\infty}\left(
E_{j}\right)  ,d^{j}\right\}  $ be an elliptic complex on $M$ which is
periodic when restricted to $\bar{N}_{+}.$ The $E_{j}$ are Hermitian vector
bundles. The theory works for differentials $d_{j}$ of any orders, thus in
particular for arbitrary elliptic operators. The operators act on
exponentially weighted Sobolev spaces. The first step is to extend
$C|_{\bar{N}^{+}}$ periodically to all of $\bar{N}.$ Call the result $\bar
{C}.$ Then $\bar{C}$ is Fredholm if and only if $C$ is. Whether $\bar{C}$ is
Fredholm is determined by the cohomology of a family of complexes on $N$
indexed by $\lambda\in\mathbb{C}^{*}.$

We sketch the construction. It is based on Fourier series for an infinite
cyclic covering, generalizing the covering of a point by $\mathbb{Z}.$ We work
in the context of Section \ref{ss5.1}. The transformation $\lambda\tau$ can be
generalized in the case $\pi=$ $\mathbb{Z}.$ We replace the regular
representation on $C^{\ast}\left(  \mathbb{Z}\right)  $ by the nonunitary
representation where $z^{n}$ acts by $k^{-n}z^{n}$ for some $k>0.$ Let
$\psi_{k}$ be the associated flat bundle. Extend the definition of $\tau$ by
$\tau u=\sum_{n}uz^{-n}\otimes k^{-n}z^{n}$ for $u\in C_{c}^{\infty}\left(
\bar{E}_{j}\right)  .$ This is an invariant section of $\bar{E}_{j}\otimes
\bar{\psi}_{k}.$ The weighted $C^{\ast}\left(  \mathbb{Z}\right)  $-inner
products on invariant sections are gotten by replacing $dx$ by $k^{2\rho
\left(  x\right)  }dx.$ The component of $1$ of the induced inner product on
$C_{c}^{\infty}\left(  \bar{E}_{j}\right)  $ is the $k$-inner product. As in
Section \ref{ss5.1}, there is an induced elliptic complex on $N$ with
coefficients in $\psi_{k}.$ Since $C^{\ast}\left(  \mathbb{Z}\right)
=C\left(  S^{1}\right)  $, this corresponds to a family of elliptic complexes
on $N$ parametrized by $\left\{  \lambda|~\left\vert \lambda\right\vert
=k\right\}  .$ This consists of the quotient complex $C_{N}$ of $\bar{C}$ with
coefficients in a family of flat line bundles $\left\{  L_{\lambda}\right\}  $
on $N.$ $L_{\lambda}=\bar{N}\times\mathbb{C}/\left\{  \left(  x,c\right)
=\left(  zx,\lambda c\right)  \right\}  .$ It may be considered as an
unparametrized complex $C_{N\times S^{1}}$ over $N\times S^{1}.$ The Fourier
coefficient of $1$ of the families inner product is the $\mathcal{L}^{2}$
inner product. Thus $\lambda\tau$ induces an isomorphism between $\bar{C}_{k}$
and the $\mathcal{L}^{2}$ completion of $C_{N\times S^{1}}.$ When $N$ is a
point this is the Parseval theorem.

\begin{theorem}
\label{th6.2}\cite[Section 4]{tau} The following are equivalent.
\end{theorem}

\begin{enumerate}
\item $\bar{C}_{k}$ is Fredholm.

\item $\bar{C}_{k}$ is contractible.

\item The cohomology of the family vanishes for all $\lambda$ such that
$\left|  \lambda\right|  =k.$
\end{enumerate}

Under the assumption that the Euler characteristic of $C_{N}$ vanishes, and a
further condition on its symbol, Taubes then shows that $\bar{C}_{k}$ is
Fredholm for all but a discrete set of $k.$ The results also hold if the
differentials are asymptotically periodic in the sense that they converge to
periodic operators in the direction of the end.

The original work of Lockhart and McOwen \cite{lomc} dealt with manifolds with
cylindrical ends of the form $V\times\mathbb{R}_{+}$ and elliptic operators
$D$ invariant on the ends by translation by $\mathbb{R}_{+}.$ In this case $D$
splits as $b\left(  x\right)  \dfrac\partial{\partial t}+A,$ Where $A$ is an
operator on $V$ and $x\in V.$ A family of operators $D_{\lambda}$ on $V$ is
obtained by replacing $\dfrac\partial{\partial t}$ by $i\lambda\in\mathbb{C}.$
It is shown that $D_{k}$ is Fredholm on $\bar{N}$ if and only if $D_{\lambda}$
is invertible for all $\lambda$ such that $\operatorname{Im}\lambda=\log k.$ A
translation to the $\mathbb{Z}$-periodic situation can be accomplished as
follows. The quotient of $V\times\mathbb{R}_{+}$ by $\mathbb{N}$ is $N=V\times
S^{1},$ with the induced operator $D_{N}. $ $D_{N}$ with coefficients in the
family of flat bundles is invertible for exactly the same $k.$ As a result,
all the previously stated results hold. The assumptions used by Taubes to
establish the existence of a large set of Fredholm values of $k$ are automatic
in this case.

Theorem \ref{th6.2} gives another proof (for $A=\mathbb{C}$) that the
operators considered in this paper are Fredholm for the specified values of
$k$. It doesn't seem to be sufficient to compute their indexes.

\begin{proposition}
\label{pr6.3}If $H_{\ast}\left(  M;\mathbb{C}\right)  $ is finitely generated,
the de Rham complex of $N$ with coefficients in a flat line bundle
$L_{\lambda}$ has vanishing cohomology for all $\lambda$ with $\left\vert
\lambda\right\vert >0$ sufficiently small or large.
\end{proposition}

%

%TCIMACRO{\TeXButton{Proof}{\proof}}%
%BeginExpansion
\proof
%EndExpansion
We use the de Rham theorem for closed manifolds and Poincar\'{e} duality. It
is then sufficient to prove that the local coefficient simplicial homology of
$N$ with coefficients in $L_{\lambda}$ is zero for the specified values of
$\lambda.$ Let $\bar{C}$ be the chains of $\bar{N}$. Any $\lambda\in
\mathbb{C}^{\ast}$ determines a homomorphism $e\left(  \lambda\right)
:\mathbb{C}\left[  z,z^{-1}\right]  \rightarrow\mathbb{C}$ by evaluation on
$\lambda.$ Then $\bar{C}\otimes_{e\left(  \lambda\right)  }\mathbb{C}$
computes homology with coefficients in $L_{\lambda}.$ We work in the context
of Section \ref{ss3.4}. Since $H_{\ast}\left(  M\right)  $ is finitely
generated, so is $H_{\ast}\left(  \bar{N}\right)  .$ Let $P$ be a finitely
generated complex equivalent to $\bar{C},$ and $h$ a self-equivalence of $P$
induced from $z^{-1}.$ Let $T$ be the mapping torus of $h.$ It is
$\mathbb{C}\left[  z,z^{-1}\right]  $-module equivalent to $\bar{C}.$ There is
then an equivalence $\bar{C}\otimes_{e\left(  \lambda\right)  }%
\mathbb{C\rightarrow}T\otimes_{e\left(  \lambda\right)  }\mathbb{C}.$ The
latter complex is the mapping cone of $I-\lambda h:P\rightarrow P.$ Since $P$
is finitely generated, $I-\lambda h$ is invertible for $\left\vert
\lambda\right\vert >0$ sufficiently small or large.%
%TCIMACRO{\TeXButton{End Proof}{\endproof}}%
%BeginExpansion
\endproof
%EndExpansion

Hughes and Ranicki \cite{hura} develop topological and algebraic theories in
parallel. We discuss the algebraic. The objects are complexes $\bar{C}$ of
finitely generated free right $A\left[  z,z^{-1}\right]  $-modules, where $A$
is any ring with identity. The relation between finite domination and
contractibility appears in this context as well.

The Novikov rings are $A\left(  \left(  z\right)  \right)  $ and $A\left(
\left(  z^{-1}\right)  \right)  $, which are the formal Laurent series
containing finitely many negative (resp. positive) powers of $z.$

\begin{theorem}
\label{th6.4}\cite[Theorem 1]{ran2} $\bar{C}$ is finitely dominated if and
only if the homology of the complexes $\bar{C}\otimes_{A\left[  z,z^{-1}%
\right]  }A\left(  \left(  z\right)  \right)  $ and $\bar{C}\otimes_{A\left[
z,z^{-1}\right]  }A\left(  \left(  z^{-1}\right)  \right)  $ is zero.
\end{theorem}

\noindent For the local coefficient chains of an infinite cyclic covering of a
compact manifold, the homology of one complex vanishes if and only if that of
the other does. These complexes look like $\bar{C}$ at one end and like
$\bar{C}^{\ell f}$ at the other.

There is an analogy with weighted simplicial chain complexes. If $P$ is a free
$A\left[  z,z^{-1}\right]  $-module, $P\otimes_{A\left[  z,z^{-1}\right]
}A\left(  \left(  z\right)  \right)  $ is isomorphic to $P^{0}\otimes
_{A}A\left(  \left(  z\right)  \right)  ,$ where $P^{0}$ is the module
generated by a set of free generators. Similarly for $A\left(  \left(
z^{-1}\right)  \right)  .$ As in Section \ref{ss3.3}, let $P=P^{0}%
\otimes\mathbb{C}\left[  z,z^{-1}\right]  $ be an extended $A\left[
z,z^{-1}\right]  $-module. Then $P_{\left(  k\right)  }$ is the Hilbert module
tensor product $P^{0}\otimes_{A}A\left[  z,z^{-1}\right]  _{\left(  k\right)
}.$ We may therefore think (heuristically and somewhat incorrectly) of the
chains with coefficients in the Novikov rings as corresponding to the values
$k=\infty$ and $k=0.$

A conjecture of Bueler \cite{bue} is relevant to the present paper. Let $M$ be
complete, oriented, and connected. Suppose that the Ricci curvature is bounded
below. The heat kernel $K_{t}$ for the Laplacian on functions is unique. Let
$d\mu=K_{t}\left(  x_{0},x\right)  dx$ for some fixed $x_{0}$ and $t>0.$ The
conjecture is that the weighted $\mathcal{L}_{\mathcal{\ }}^{2}$ cohomology of
$M$ is isomorphic to the de Rham cohomology. It is shown that in a variety of
situations the weighted Laplacian is Fredholm, although in most the dimension
of its kernel isn't determined. These results have limited contact with the
present paper, since $K_{t}$ tends to decay more rapidly than the weight
functions used here. Carron \cite{car} has given counterexamples to this
conjecture. The method applies only to manifolds with infinitely generated cohomology.

Yeganefar \cite{yeg} has established the equality of the weighted and de Rham
cohomologies in many cases not covered by this paper. This leads to a
topological interpretation of the $\mathcal{L}^{2}$ cohomology of manifolds
with finite volume and sufficiently pinched negative curvature. A standing
hypothesis is that $d\rho\neq0$ outside of a compact set.


\begin{thebibliography}{99}                                                                                               %


\bibitem {ang}N. Anghel, \textit{An abstract index theorem on non-compact
Riemannian manifolds}, Houston J. Math., \textbf{19} (1993), 223-237

\bibitem {botu}R. Bott and L. W. Tu, \textit{Differential Forms in Algebraic
Topology, }Springer-Verlag, 1986

\bibitem {bue}E. L. Bueler, \textit{The heat kernel weighted Hodge Laplacian
on noncompact manifolds,} Trans. Amer. Math. Soc., \textbf{351} (1999), 683-713

\bibitem {car}G. Carron, \textit{A Counter Example to the Bueler's
Conjecture,} www.arxiv.org/abs/math.DG/0509550, Proc. Amer. Math. Soc. to appear.

\bibitem {che}P. R. Chernoff, \textit{Essential self-adjointness of powers of
generators of hyperbolic systems}, J. Funct. Anal., \textbf{12} (1973), 401-414

\bibitem {fac}T. Fack, \textit{Sur la conjecture de Novikov}, in `Index Theory
of Elliptic operators, Foliations, and Operator Algebras', 43-102, Contemp.
Math., \textbf{70}, Amer. Math. Soc., 1988

\bibitem {foha}J. Fox and P. Haskell, \textit{Homology Chern characters of
perturbed Dirac operators}, Houston J. Math., \textbf{27 }(2001), 97-121

\bibitem {gro}A. Grothendieck, \textit{Produits Tensorial Topologiques et
\'{E}spaces Nucl\'{e}aires}, Mem. Amer Math. Soc., \textbf{16}, Amer. Math.
Soc., 1955

\bibitem {hil}M. Hilsum, \textit{Fonctorialit\'{e} en }$K$\textit{-th\'{e}orie
bivariante pour les vari\'{e}t\'{e}s Lipschitziennes}, $K$-Theory, \textbf{3}
(1989), 401-440

\bibitem {hud}J. F. P. Hudson, \textit{Piecewise Linear Topology}, Benjamin, 1969.

\bibitem {hura}C. B. Hughes and A. Ranicki, \textit{Ends of Complexes,
}Cambridge Univ. Press, 1996

\bibitem {kami}J. Kaminker and J. G. Miller, \textit{Homotopy invariance ot
the index of the signature operator over }$C^{\ast}$\textit{-algebras}, J.
Operator Theory, \textbf{14} (1985), 113-127

\bibitem {kari}R. Kadison and J. R. Ringrose, \textit{Fundamentals of the
Theory of Operator Algebras.} \textit{Vol. }1, Academic Press, 1983

\bibitem {kas1}G. G. Kasparov, \textit{An index for invariant elliptic
operators, K-theory, and representations of Lie groups}, Soviet Math. Dokl.,
\textbf{27 }(1983), 105-109

\bibitem {kas2}G. G. Kasparov, \textit{K-theory, group }$C^{\ast}%
$\textit{-algebras, and higher signatures (conspectus), }in\textit{\ `}Novikov
Conjectures, Index Theorems, and Rigidity'. vol. 1 , 101-146, Cambridge Univ.
Press 1995

\bibitem {lan1}E. C. Lance, \textit{K-Theory for certain group }$C^{\ast}%
$\textit{-algebras}, Acta Math., \textbf{151} (1983), 209-230

\bibitem {lan}E. C. Lance, \textit{Hilbert }$C^{\ast}$\textit{-modules},
London Math. Soc., 1995

\bibitem {lef}S. Lefschetz, \textit{Introduction to Topology}, Princeton Univ.
Press, 1949

\bibitem {lomc}R. B. Lockhart and R. C. McOwen, \textit{Elliptic differential
operators on noncompact manifolds}, Ann. Scuola Norm. Sup. Pisa cl. Sci.,
\textbf{12} (1985), 409-447

\bibitem {mel}R. B. Melrose, \textit{The Atiyah-Patodi-Singer Index Theorem},
A. K. Peters, 1993

\bibitem {mil1}J. G. Miller, \textit{Signature operators and surgery groups
over }$C^{\ast}$\textit{-algebras}, $K$-Theory, \textbf{13} (1998), 363-402

\bibitem {mil2}J. G. Miller, \textit{Differential operators over }$C^{\ast}%
$\textit{-algebras}, Rocky Mountain J. Math,. \textbf{29 }(1999), 239-269

\bibitem {mun}J. R. Munkres, \textit{Topology}, 2nd ed., Prentice Hall, 2000

\bibitem {pal}R. S. Palais, \textit{Differential operators on vector bundles},
in `Seminar on the Atiyah-Singer Index Theorem', 51-92, Princeton Univ. Press, 1965

\bibitem {pan1}P. Pansu,\textit{\ Cohomologie }$L^{p}$\textit{: invariance
sous quasiisometries, }www.math.upsud.fr/\symbol{126}pansu, preprint 2004

\bibitem {pan2}P. Pansu, \textit{Introduction to }$L^{2}$\textit{\ Betti
numbers}, in `Riemannian Geometry', 53-86, Amer. Math. Soc., 1996

\bibitem {pas}W. Paschke, \textit{Inner product modules over }$B^{\ast}%
$\textit{-algebras}, Trans. Amer. Math. Soc. \textbf{182 }(1973), 443-468

\bibitem {pivo}M. Pimsner and D. Voiculescu, \textit{Exact sequences for
K-groups and Ext-groups of certain cross-products of }$C^{\ast}$%
\textit{-algebras}, J. Operator Theory, \textbf{4} (1980), 93-118.

\bibitem {ran2}A. Ranicki, \textit{Finite domination and Novikov rings},
Topology, \textbf{34 }(1995), 619-632

\bibitem {roe}J. Roe, \textit{Index theory, Coarse geometry, and the Topology
of Manifolds}, Amer. Math. Soc., 1996

\bibitem {sch}J. A. Schafer, \textit{The Bass conjecture and group von Neumann
algebras}, $K$-Theory, \textbf{19 }(2000), 211-217

\bibitem {seg}G. Segal, \textit{Fredholm complexes}, Q. J. Math., \textbf{21}
(1970), 385-402

\bibitem {shu}M. Shubin, \textit{Spectral theory of elliptic operators on
non-compact manifolds}, Ast\'{e}risque, \textbf{207} (1992), 37-108

\bibitem {sie}L. C. Siebenmann, \textit{The structure of tame ends}, Notices
Amer. Math. Soc., \textbf{66t-G7} (1966), 861

\bibitem {sim}B. Simon, \textit{Representations of Finite and Compact Groups},
Amer. Math. Soc., 1996

\bibitem {ste}N. Steenrod, \textit{The Topology of Fiber Bundles}, Princeton
Univ. Press 1951.

\bibitem {tau}C. H. Taubes, \textit{Gauge theory on asymptotically periodic
4-manifolds}, J. Differential Geom., \textbf{25} (1987), 363-430

\bibitem {wal1}C. T. C. Wall, \textit{Finiteness conditions for }%
$CW$\textit{\ complexes}, Ann. of Math., \textbf{81 }(1965), 56-69

\bibitem {wal2}C. T. C. Wall, \textit{Finiteness conditions for }%
$CW$\textit{\ complexes II}, Proc. Royal Soc. London Ser. A, \textbf{295}
(1965), 129-139

\bibitem {w-o}N. E. Wegge-Olsen, $K$\textit{-Theory and }$C^{\ast}%
$\textit{-algebras}, Oxford Univ. Press, 1993

\bibitem {wei}J. Weidmann, \textit{Linear Operators in Hilbert Spaces},
Springer-Verlag, 1980

\bibitem {whi}H. Whitney, \textit{Geometric Integration Theory}, Princeton
Univ. Press, 1957

\bibitem {yeg}N. Yeganefar, \textit{Sur la }$L^{2}$\textit{-cohomologie des
vari\`{e}t\'{e}s \'{a} courbure n\'{e}gative}, Duke Math. J., \textbf{122}
(2004), 145-180
\end{thebibliography}
\end{document}